\documentclass[12pt]{amsart}
\pdfoutput=1
\usepackage{microtype}
\usepackage[a4paper, margin=1in]{geometry}
\usepackage[all]{xy}
\usepackage{mathtools}
\usepackage{hyperref}
\usepackage{amssymb}
\usepackage{mathrsfs}
\usepackage{tikz}
\usepackage{tikz-cd}
\hypersetup{
    colorlinks=true,
    linkcolor=black,
    filecolor=black,
    citecolor= black,
    urlcolor=black}

\usepackage{cleveref}

\newtheorem*{theorem*}{Theorem}
\newtheorem{theorem}{Theorem}[section]
\newtheorem{lemma}[theorem]{Lemma}
\newtheorem{proposition}[theorem]{Proposition}
\newtheorem{conjecture}[theorem]{Conjecture}
\theoremstyle{definition}
\newtheorem{definition}[theorem]{Definition}
\newtheorem{example}[theorem]{Example}

\newtheorem{corollary}[theorem]{Corollary}
\newtheorem{question}[theorem]{Question}

\theoremstyle{remark}
\newtheorem{remark}[theorem]{Remark}

\newcommand{\rank}{\textup{rank}}

\newcommand{\Tor}{\textup{Tor}}

\newcommand{\Supp}{\textup{Supp}}
\newcommand{\Spec}{\textup{Spec}}
\newcommand{\ann}{\textup{ann}}

\newcommand{\Assh}{\textup{Assh}}
\newcommand{\ceil}[1]{\lceil {#1} \rceil}
\newcommand{\flor}[1]{\lfloor{#1}\rfloor}

\newcommand{\bp}[2]{{#1}^{[p^{#2}]}}
\newcommand{\bpq}[2]{{#1}^{[#2]}}
\newcommand{\dash}{\rule{.4cm}{.1mm}}
\newcommand{\Hom}{\textup{Hom}}

\numberwithin{equation}{section}

\usepackage [backend= biber, style= alphabetic, sorting= nty]{biblatex}
\begin{document}

%%\subjclass[2010]{Primary }

\title{$h$-function, Hilbert-Kunz density function and Frobenius-Poincar\'e function}
\author{Cheng Meng}
\address{Yau Mathematical Sciences Center, Tsinghua University, Beijing 100084, China.}
\email{cheng319000@tsinghua.edu.cn}
\author{Alapan Mukhopadhyay}
\address{Institute of Mathematics, CAG,
EPFL SB MATH
MA A2 383 (Bâtiment MA), 
Station 8,
CH-1015 Lausanne,
Switzerland}
\email{alapan.mathematics@gmail.com}
\date{\today}
\maketitle
\begin{abstract}
Given ideals $I,J$ of a noetherian local ring $(R, \mathfrak m)$ such that $I+J$ is $\mathfrak m$-primary and a finitely generated $R$-module $M$, we associate an invariant of $(M,R,I,J)$ called the $h$-function. Our results on $h$-functions allow extensions of the theories of Frobenius-Poincar\'e functions and Hilbert-Kunz density functions from the known graded case to the local case, answering a question of V.Trivedi. When $J$ is $\mathfrak m$-primary, we describe the support of the corresponding density function in terms of other invariants of $(R, I,J)$. We show that the support captures the $F$-threshold: $c^J(I)$, under mild assumptions, extending results of V. Trivedi and Watanabe. The $h$-function encodes Hilbert-Samuel, Hilbert-Kunz multiplicity and $F$-threshold of the ideal pair involved. Using this feature of $h$-functions, we provide an equivalent formulation of a conjecture of Huneke, Mustaţă, Takagi, Watanabe; recover a result of Smirnov and Betancourt; give a new proof of a result answering Watanabe-Yoshida's question comparing Hilbert-Kunz and Hilbert-Samuel multiplicity and establish lower bounds on $F$-thresholds. We also point out that a conjecture of Smirnov-Betancourt as stated is false and suggest a correction which we relate to the conjecture of Huneke et al.

We develop the theory of $h$-functions in a more general setting which yields a density function for $F$-signature. A key to many results on $h$-functions is a `convexity technique' that we introduce, which in particular proves differentiability of Hilbert-Kunz density functions almost everywhere on $(0,\infty)$, thus contributing to another question of Trivedi.
\end{abstract}
\tableofcontents
\section{Introduction}\label{se: introduction}

Hilbert-Kunz multiplicity and $F$-signature are numerical invariants extensively appearing in commutative algebra and algebraic geometry in prime characteristics; see  \cite{HanMonsky93, SmithVdb97, BuchweitzChen97, SecCoef, Vraciu08, Hanes, Brrationality, BrennerMonsky10, Ilya16, perfectoidsignature23, AberbachSteib24} for instance. These quantify severity of singularities at a point of a variety and also relate to other invariants, such as the cardinality of the local \'etale fundamental group of the punctured spectrum of a strongly $F$-regular local ring; see \cite{AberbachEnescu}, \cite{HunekeYao02}, \cite{BE}, \cite{JavierSchwedeTucker} and the references in \Cref{se: background material}. The theory of Hilbert-Kunz multiplicity in the graded setup has witnessed two new natural generalizations in recent years: The Hilbert-Kunz density function and the Frobenius-Poincar\'e function. Let $S$ be a standard graded ring in prime characteristics and $\mathfrak a$ be a homogeneous ideal of finite colength. When $\dim(S) \geq 2$, Trivedi proves that the Hilbert-Kunz density function $g_{S, \mathfrak a}$ is a compactly supported real valued continuous function  of a real variable, whose integral recovers the Hilbert-Kunz multiplicity $e_{HK}(\mathfrak a, S)$; see \cite{TriExist}, \Cref{se: background material} for details. For the graded pair $(S, \mathfrak a)$, where $\dim(S)$ is not necessarily at least two, the associated Frobenius-Poincar\'e function is shown to be an entire function in one complex variable, whose value at the origin is the Hilbert-Kunz multiplicity $e_{HK}(\mathfrak a, S)$; see \cite{AlapanExist}, \Cref{se: background material}. These two functions not only encode finer (graded)invariants of $(S, \mathfrak a)$ than the Hilbert-Kunz multiplicity but also allow application of geometric tools, such as sheaf cohomology on $\text{Proj}(S)$, and tools from homological algebra. Successful applications of the Hilbert-Kunz density functions have resolved Watanabe and Yoshida's conjecture on the values of Hilbert-Kunz multiplicity of quadric hypersurfaces in large enough characteristics, rationality of Hilbert-Kunz multiplicities and $F$-thresholds of two dimensional normal rings among other results; see \cite{TrivediQuadric}, \cite{TriFthreshold}, \cite{TriCurves}, \cite{TrivediCharZero}.\\

Building extensions of these two theories to the setting of a noetherian local ring is a natural question; see Trivedi's question \cite[Question 1.3]{TriExist}. In this article, we extend the theories of Hilbert-Kunz density functions and Frobenius-Poincar\'e functions to the local setting. In the local setting, the role of a grading is replaced by an ideal adic filtration. However, the techniques used in the graded setup do not extend. We handle the additional challenges in this greater generality by developing a theory of, what we call the $h$-function. The theory of $h$-functions that we develop proves properties of density functions previously unknown even in the classical graded setup, \Cref{pr: monotonicity result for the desnity function}, \Cref{co: HK density is strictly positive in the interior} et cetera. Moreover, our results on $h$-functions yield applications to problems in multiplicity theories appearing in \cite{WY00}, \cite{HMTW}, \cite{BetancortSmirnov} which are a priori unrelated to Hilbert-Kunz density functions and Frobenius-Poincar\'e functions; see \Cref{se: applications}. A brief overview follows:\\

Fix a noetherian local domain  $(R, \mathfrak m)$ of prime characteristic $p>0$ and  Krull dimension $d$, where the Frobenius endomorphism is a finite map. Fix two ideals $I,J$ of $R$ such that $I+J$ is $\mathfrak m$-primary. We prove:\\

\noindent \textbf{Theorem A:} Consider the sequence of functions of a real variable
\[h_{n,I,J}(s)= l(\frac{R}{(I^{\ceil {sp^n}}+ \bp{J}{n})R}),\]
where $\bp{J}{n}$ is the ideal generated by $\{f^{p^n} \, |\, f \in J\};$ and $l(\dash)$ is the length function.

\begin{enumerate}
    \item (\Cref{th: existence of h-function for a domain}, \Cref{th: uniform convergence of h function for modules}) There is a real-valued function of a real variable denoted by $h_{I,J}(s)$ such that given an interval $[s_1, s_2] \subseteq \mathbb R$, there is a constant $C$ depending only on $s_1, s_2$ satisfying
    \[|h_{I,J}(s)-\frac{h_{n,I,J}(s)}{p^{nd}}| \leq \frac{C}{p^n}, \, \text{for all}\, s \in [s_1, s_2] \, \text{and}\, n \in \mathbb{N}.\]

    \noindent Consequently, the sequence of functions $\frac{h_{n,I,J}(s)}{p^{nd}}$ converges to $h_{I,J}(s)$ and the convergence is uniform on every compact subset of $\mathbb R$.

    \item (\Cref{pr: cotinuity of $h$ functions}, \Cref{th: Lipschitz continuity for a family}) Given real numbers $s_2> s_1>0$, there is a constant $C'$- depending only on $s_1, s_2$ such that for $x, y \in [s_1, s_2]$,

    \[|h_{I,J}(x)- h_{I,J}(y)| \leq C'|x-y| .\]
     That is, away from zero, $h_{I,J}$ is locally Lipschitz continuous.
\end{enumerate}
The function $h_{I,J}$ is called the \textit{h}-function associated to the pair $(I,J)$.
\vspace{.3cm}
 In fact we prove a version of the above theorem for an ideal $I$ and a family of ideals $J_{\bullet}$ satisfying what we call \textbf{Condition C},  allowing for applications to other numerical invariants such as $F$-signature; see \Cref{th: existence of h-function for a domain}.\\

Special instances of this $h$-function have been considered by different authors: in \cite{TaylorInterpolation} when both $I$ and $J$ are $\mathfrak{m}$-primary; in \cite{BlicleSchwedeTuckerFractal} when $R$ is regular, $I$ is principal and $J= \mathfrak m$ to study $F$-signature of a pair and in \cite{kosuke2017function} in the same set up but in a different context. More recently, the function $\phi_J(R,z^t)$ featured in \cite{AberbachSteib24}\footnote{This appeared after our preprint was posted.} and played a crucial role in the resolution of Watanabe-Yoshida conjecture up to dimension seven. This $\phi_J(R,z^t)$ also is another special instance of $h$-functions. Our theory of $h$-functions generalizes these results.  Moreover, the techniques involved in our proofs yield uniform convergence which is crucial for us.\\

\noindent \textbf{Local Frobenius-Poincar\'e functions:} In \Cref{th: bounds on h function special case}, we prove that there is a polynomial $P_1(s)$ of degree $\dim(R/J)$ such that $h_{I,J}(s) \leq P_1(s)$ for all $s$. Using this polynomial bound, we prove existence and holomorphicity of a function $F_{R,I,J}(y)$ on the open lower half complex plane, which we call the Frobenius-Poincar\'e function in the local setting; see \Cref{th: exitence of Frob-Poincare for general ideals}. We moreover show: 
\[F_{R,I,J}(y)= \underset{\mathbb R}{\int} h_{I,J}(t)e^{-ity}iydt. \]
When $J$ is $\mathfrak{m}$-primary, we prove $F_{R,I,J}(y)$ is entire. When \textit{$(R, \mathfrak m, J)$ comes from a graded pair $(S, \mathfrak a)$}, i.e. $(R,\mathfrak{m})$ is the localization of a standard graded ring $S$ at the homogeneous maximal ideal, $I$ is the homogeneous maximal ideal and $J$ comes from a homogeneous ideal of finite colength $\mathfrak a$, $F_{R,I,J}(y)$ coincides with the Frobenius-Poincar\'e function of the pair $(S, \mathfrak a)$ as developed in \cite{AlapanExist}; see \Cref{pr: the graded Frob-Poincare coincides wth the local Frob Poincare},(3). Unlike \cite{AlapanExist}, our treatment allows us to consider Frobenius-Poincar\'e function of 
a graded pair $(S, \mathfrak a)$, where $\mathfrak a$ need not have finite colength; see \Cref{pr: the graded Frob-Poincare coincides wth the local Frob Poincare}, (2).\\

\noindent \textbf{Local Hilbert-Kunz density functions:} Extending the theory of Hilbert-Kunz density functions is more involved. Set

\[f_n(s)= h_{n,I,J}(s+\frac{1}{p^n})- h_{n,I,J}(s).\]

When $(R,\mathfrak{m},J)$ comes from a graded pair $(S, \mathfrak a)$, as is meant above, where $\dim(S) \geq 2$, we point out that the sequence of functions
\[\frac{f_n(s)}{(p^n)^{d-1}}\]

\noindent converges uniformly to the Hilbert-Kunz density function of $(S, \mathfrak a)$ defined by Trivedi using the underlying graded structure in \cite{TriExist}; see \Cref{th: h function differentiable when J is m primary}. So, for a local ring $(R, \mathfrak{m})$, we take the limit of $\frac{f_n(s)}{(p^n)^{d-1}}$, whenever it exists, as the value of the Hilbert-Kunz density function associated to $(R,I,J)$ at the point $s$. For arbitrary ideals $I,J$ of a local ring $(R, \mathfrak m)$, the pointwise convergence of $f_n(s)/(p^n)^{d-1}$ at every point is not clear. In fact when $I=0$ the sequence does not converge at $s=0$, see \Cref{eg: where the limit defining the desnity function does not exist}. We handle the convergence problem for $f_n(s)/(p^n)^{d-1}$ by relating it to the differentiability of $h_{I,J}$ at $s$. In \Cref{th: differentiability of h function implies existence of density function} we prove, 
\begin{center}
\textit{If $h_{I,J}$(s) if differentiable at $s$, $f_n(s)/(p^n)^{d-1}$ converges to $h'_{I,J}(s)$}.
\end{center}

\noindent We warn that differentiability of $h_{I,J}$ is not necessary for the above convergence, as \Cref{eg: density functions exist but is not continuous} shows. In the direction of differentiability of $h$ and its relation with the local density function, we prove:\\

\noindent \textbf{Theorem B:} Let $h_{I,J}$ be as before.
\begin{enumerate}
    \item The left and right hand derivative of $h_{I,J}$ exist at each nonzero real number; see \Cref{th: one sided differentiability of h-functions}, (3).
    
    \item Outside a countable subset of $(0,\infty)$, $h_{I,J}$ is continuously differentiable; see \Cref{th: one sided differentiability of h-functions}, (3).

    \item The set of points where $f_n(s)/(p^n)^{d-1}$ converges and the limiting function is continuous coincides with the set of points where $h_{I,J}$ is continuously differentiable; see \Cref{th: continuous differentiability vs continuity of the density function}, \Cref{de: continuity of density function}, \Cref{de: continuous differentiability of a function}.
\end{enumerate}

\noindent \textbf{Theorem B}, (2) along with \Cref{th: differentiability of h function implies existence of density function} imply that if $(R,\mathfrak{m})$ is local, then for any $R$-ideals $I,J$ as above, $f_n(s)/(p^n)^{d-1}$ converges outside a countable subset of $\mathbb R$ and coincides with the derivative of $h_{I,J}(s)$. Thus outside this countable set, the limiting function $f_n(s)/(p^n)^{d-1}$ yields a well-defined notion of density function. When $J$ is $\mathfrak{m}$-primary and $I$ is non-zero, the Lebesgue integral of the  density function is indeed the Hilbert-Kunz multiplicity $e_{HK}(J,R)$; see \Cref{pr: integral of the density function}.\\

\noindent \textbf{Advantages of the $h$-function approach:} One advantage of our approach via $h$-functions, over the earlier approaches in the graded setup is that many of our results (e.g. \Cref{th: differentiability of h function implies existence of density function}, \textbf{Theorem B}) on density functions hold more generally for a family satisfying \textbf{Condition C}. This generalization in particular yields a theory of density function for $F$-signature which is absent even in the graded setup. We focus on this $F$-signature density function in a forthcoming work. Note that, the monotonicity property of the density function (\Cref{pr: monotonicity result for the desnity function}) is another consequence of our approach, which was previously unknown and perhaps a priori unexpected even in the graded setup of \cite{TriExist}.\\ 

\noindent \textbf{Recovering the graded case:} When $(R,\mathfrak{m},J)$ comes from a graded pair $(S, \mathfrak a)$ with $\dim(S) \geq 2$, we prove that the corresponding $h$-function $h_{R,\mathfrak{m},J}$ is continuously differentiable and the derivative coincides with the Hilbert-Kunz density function that Trivedi defines; see \Cref{th: h function differentiable when J is m primary}. We moreover generalize existence and continuity of the graded density function to the case when the homogeneous ideal $\mathfrak a$ does not have finite colength; see \Cref{th: HK density when J is homogeneous but not of finite colength}. These two results essentially stem from \textbf{Theorem B}, (3) which in turn follows from the integral formula (\Cref{th: h is the integral of f}) and  existence of one sided limit for the density function and the one sided derivatives of $h_{I,J}$ (\Cref{pr: one sided limits exists}). Outside the graded setup, our results \Cref{th: h function and integration} and \Cref{eg: examples of differentiable density function} allow constructions of $h$-function and density function with a given order of smoothness.\\

\noindent \textbf{Choices of definitions for the local Hilbert-Kunz density function:} We now comment on other possible choices for the density function in the local setup. One could take either the left or the right derivative of $h_{I,J}$ as the density function corresponding to $I,J$: Both of these are integrable with integral taking the value $e_{HK}(J,R)$, when $J$ is $\mathfrak{m}$-primary and $I$ is nonzero; see \Cref{pr: integral of the density function}. Such a choice would produce density function defined everywhere on $\mathbb{R}-\{0 \}$ and everywhere on $\mathbb{R}$, when $I$ is nonzero. Moreover, for a graded pair $(S, \mathfrak{a})$ as above, with $\dim(S) \geq 2$, either of the choices coincides with Trivedi's density function defined in \cite{TriExist}. Still, we choose to work with $f_n(s)$, because of its interpretation as the length of $\ceil{sp^n}$-th graded piece of the $I$-adic filtration on $R/\bp{J}{n}$. In general, the left derivative of $h_{I,J}$ may not coincide with the limit of $(f_n(s)/(p^n)^{d-1})_n$; see \Cref{eg: h function is not differentiable at a point} and the examples in \Cref{sse: examples} for subtleties with alternate definitions of density functions. We do not know of an example where the right derivative of $h_{I,J}$ differs from the limit of $(f_n(s)/(p^n)^{d-1})_n$ on $\mathbb{R}_{>0}$; see \Cref{qe: continuity of density function}.\\

\noindent \textbf{A convexity technique:} Central to many results in this article, including \textbf{Theorem B}, is a `convexity technique' that we introduce. In \Cref{th: uniform convergence to the convex functional}, the `convexity technique' manifests into the convexity of the following function on $[s_0, \infty)$, given any $s_0>0$,
\[s \mapsto H(s,s_0)= h_{I,J}(s)/c(s)-h_{I,J}(s_0)/c(s_0)+\int^{s}_{s_0} h_{I,J}(t)c'(t)/c^2(t)dt,\]

\noindent where $c(s)= s^{\mu-1}/(\mu-1)!$, $\mu$ being the cardinality of a set of generators of $I$. \textbf{Theorem B} then follows from general properties of convex functions. \Cref{th: uniform convergence to the convex functional} has numerous applications: The left and right limit of the density function exists at every point on $\mathbb{R}_{>0}$, even though density function is not known to be defined everywhere; see \Cref{pr: one sided limits exists}. Moreover, we prove that the density function is positive in the interior of the support, this property was unknown even in the graded setup; see \Cref{c: minimal stable point is the support}, \Cref{co: HK density is strictly positive in the interior}. Another consequence is the twice differentiablility of the $h$-function outside a set of measure zero, contributing to Trivedi's question about the order of differentiability of the Hilbert-Kunz density function; see \cite[Question 1]{TrivediQuadric}, \Cref{re: density function is differentiable outside a set of measure zero}. The underlying idea of the same convexity argument is used to prove Lipschitz continuity of $h$-functions stated in \textbf{Theorem A}.\\

\noindent\textbf{Invariants encoded by $h$-functions:} The behaviour of $h_{I,J}$ near zero essentially encodes the Hilbert-Kunz multiplicity of $R/I$ with respect to $J$:\\

\noindent \textbf{Theorem C:}(\Cref{th: asymptotic behaviour of h function near zero}) Suppose $\dim(R/I)=d'$. Denote the set of minimal primes of $R/I$ of dimension $d'$ by $\text{Assh}(R/I)$. Then
    \[ \underset{s \to 0+}{\lim}\frac{h_{I,J}(s)}{s^{d-d'}}=\frac{1}{(d-d')!}\underset{P \in \Assh(R/I)}{\sum}e_{HK}(J,R/P)e(I,R_P),\]
where $e(I, \dash)$ denotes the Hilbert-Samuel multiplicity with respect to $I$. In particular, the order of vanishing $h_{I,J}(s)$ at $s=0$ is $d-d'$ and $h_{I,J}$ is continuous at zero if and only if $I$ is nonzero. \textbf{Theorem C} extends part of \cite[Theorem 4.6]{BlicleSchwedeTuckerFractal}, where $R$ is regular, $I$ is a principal ideal and $J= \mathfrak m$.\\

The $h$-function treats different numerical invariants of $(R,I,J)$ on an equal footing. When $I$ is $\mathfrak{m}$-primary, for $s>0$ and close to zero $h_{I,J}(s)= e(I,R)\frac{s^d}{d!}$; see \Cref{le: h function is stable near boundaries}.
When $J$ is $\mathfrak{m}$-primary, for large $s$, $h_{I,J}(s)= e_{HK}(R,J)$. Since $h_{I,J}(s)$ is a increasing function there is a smallest point after which $h_{I,J}(s)= e_{HK}(R,J)$. We describe this `minimal stable point' of $h_{I,J}(s)$ in terms another numerical invariant. \\ 

\noindent \textbf{Theorem D:}(\Cref{th: maximal stable point of h function for J}, \Cref{le: F threshold upto tight closure}) Suppose $J$ is $\mathfrak m$-primary. Let $\alpha_{R,I,J}= \text{sup}\{s \in \mathbb R \, |\, s>0 \, ,h_{I,J}(s) \neq e_{HK}(J,R) \}.$ Consider the sequence of numbers,
\[r^J_I(n)=\max\{t \in \mathbb{N}|I^t \nsubseteq (J^{[p^n]})^*\},\]
where $(J^{[p^n]})^*$ denotes the tight closure of the ideal $(J^{[p^n]})$; see \Cref{de: tight closure}. Then 
$(r^J_I(n)/p^n)_n$ is an increasing sequence converging to $\alpha_{R,I,J}$.\\
The above description of $\alpha_{R,I,J}$ resembles that of the well known $F$-threshold $c^J(I)$; see \Cref{de: definition of F-limbus and F-threshold}. Recall that $F$-threshold is an invariant extensively studied in prime characteristic singularity theory; see \cite{HMTW}, \cite{MTW} and is closely related log canonical threshold via reduction modulo $p$; see \cite{TW}, \cite{HW}. In general, $\alpha_{R,I,J}$ is bounded above by $c^J(I)$. We prove, under suitable hypothesis, for example, strong $F$-regularity at every point of $\text{Spec}(R)-\mathfrak \{\mathfrak{m}\}$, $\alpha_{R,I,J}$ coincides with $c^J(I)$; see \Cref{th: when the maximal support coincides with F-threshold}. \Cref{th: maximal stable point of h function for J}, \Cref{th: when the maximal support coincides with F-threshold}, \Cref{c: minimal stable point is the support} extend Trivedi and Watanabe's description of the support of Hilbert-Kunz density function when $R$ is graded and strongly $F$-regular on the punctured spectrum; see \Cref{re: maximal support of h and relation to Trivedi;s result}, \cite[Theorem 4.9]{TriFthreshold}. We do not know whether $\alpha_{R,I,J}$ always coincides with $c^J(I)$; see \Cref{qe: minnimal stable point vs F-threshold}. We prove that on $(0, \alpha_{R,I,J})$, the function $h_{I,J}$ is in fact strictly increasing; see \Cref{c: minimal stable point is the support}.\\

\noindent \textbf{Applications:} Our applications of the theory of $h$-functions in \Cref{se: applications} highlight its feature of capturing different numerical invariants such as $F$-threshold, Hilbert-Kunz and Hilbert-Samuel multiplicity simultaneously. In \Cref{sse: Watanabe's question}, we answer Watanabe's question comparing Hilbert-Samuel and Hilbert-Kunz multiplicity. Although Watanabe's question had been affirmatively confirmed by Hanes, the motivation behind this question was obscure- we do not find any mention of the motivation in the existing literature. We think that our approach using $h$-function, which is different from Hanes', naturally yields the question itself. Indeed, in \Cref{pr: algebraic version of Watanabe's question}, we first prove an equivalence between Watanabe's question and whether the minimal stable point $\alpha_{R,I,I}$ of $h_{I,I}$ is strictly greater than one, provided $\dim(R)$ is at least two. We prove that the latter question is equivalent to showing $I^{p^n+1} \nsubseteq (\bpq{I}{p^n})^*$ for large enough $n$; see \Cref{pr: algebraic version of Watanabe's question}. The last question appears to be more plausible to us, because of the metaprinciple: `Tight closure is small'.\footnote{this is the reason for assigning the attribute `tight' to the closure operation. Aside, see \Cref{re: continuity of density function and Watanabe's inequality} for another potential motivation for Watanabe-Yoshida question.} We affirmatively answer this containment problem in \Cref{th: comparison of powers and powers of tight closure}. Our result also implies that the $F$-threshold $c^I(I)>1$, when $\dim(R)$ is at least two; see \Cref{co: lower bound on F threshold} for a more general statement.  As further applications, we establish comparisons between multiplicities, which are similar to those appearing in the equimultiplicity theories; see \Cref{pr: comparison coming from behaviour of h function at infinity}. The main tools in \Cref{pr: comparison coming from behaviour of h function at infinity} are our result on the asymptotic behaviour of $h_{I,J}(s)$ as $s$ becomes large (\Cref{th: behaviour of h at infinity}) and the monotonicity property of $h$ functions; see \Cref{th: decreasingness of normalized h function}. As a consequence, we show that $h_{I,J}(s)= e(I+J, R)s^r/r!$, when $I$ and $J$ are generated by parts of a full system of parameters consisting of $ 0 \leq r \leq d$ and $d-r$ elements respectively; see \Cref{pr: h functions when the ideals split a full system of parameters}. 

In a different direction, we show that even a coarse approximation of an $h$-function recovers Smirnov and Betancourt's result comparing Hilbert-Kunz multiplicity of a ring and its quotient by part of a system of parameters; see \Cref{th: Hilbert-Kunz multiplicity and F-thresold}. However, we point out that their conjecture motivating their aforesaid result is false; see \Cref{pr: Smirnov Betancourt conjecture is false}. So we propose a corrected version in \Cref{co: corrected F-threshold and Hilbert-Kunz}. Indeed \Cref{pr: equivalence of conjectures} shows that a special case of our proposed version is equivalent to another conjecture of Huneke, Mustaţă, Takagi and Watanabe; the latter has been already verified when the rings and ideals involved are graded. We show that the general case of Huneke et al.'s conjecture is equivalent to a question about $h$-functions; see \Cref{pr: HMTW conjecture in terms of the h function}. These applications are facilitated by \Cref{th: decreasingness of normalized h function} and \Cref{pr: monotonicity result for the desnity function}, which in turn relies on our `convexity technique'.\\

Some questions regarding $h$-functions and the resulting density function are listed in \Cref{se: questions}.\\

\noindent \textbf{Notation and conventions:} All rings are commutative and noetherian. The symbol $p$ denotes a positive prime number. Unless otherwise said, the pair $(R, \mathfrak m)$ denotes a noetherian local ring $R$- not necessarily a domain- with maximal ideal $\mathfrak m$. By saying $(R, \mathfrak m)$ is graded, we mean $R$ is a standard graded ring with homogeneous maximal ideal $\mathfrak{m}$. When $(R,\mathfrak{m})$ is assumed to be graded, $R$-modules and ideals are always assumed to be $\mathbb Z$-graded. The local or graded ring denoted by the symbol $R$ is assumed to have prime characteristic $p>0$ and $R$ is assumed to be \textit{F}-finite, i.e. the Frobenius endomorphism of $R$ is finite. We index the sequences of numbers and functions by nonnegative integers $n$. Whenever the letter $q$ appears in such a sequence, $q$ denotes $p^{n}$. For an ideal $I \subset R$, $\bp{I}{n}$ or $\bpq{I}{q}$ denotes the ideal generated by $\{f^{q} \, |\, f \in I\}$ and is called the $q$ or $p^n$-th \textit{Frobenius power} of $I$. The operator $l_R(\dash)$ or simply $l(\dash)$ denotes the length function. For an $R$-module $M$, $F^n_*M$ denotes the $R$-module whose underlying abelian group is $M$, but the $R$-action comes from restriction scalars through the iterated Frobenius morphism $F^n: R \rightarrow R$. We use the attributes increasing and non-decreasing and similarly, decreasing and non-increasing interchangeably. 

\section{Background material}\label{se: background material}
We recall some results on Hilbert-Kunz multiplicity, Hilbert-Kunz density function and Frobenius-Poincar\'e function. In this section, let $(R,\mathfrak m)$ be a noetherian local or graded ring, $J$ be an $\mathfrak m$-primary ideal, $M$ be a finitely generated $R$-module.

Although the germ of Hilbert-Kunz multiplicity was present in Kunz's seminal article \cite{Kunz}, its existence was not proven until Monsky's work. 

\begin{theorem}(see \cite{MonExist})\label{th: existence of HK multiplicity}
There is a real number denoted by $e_{HK}(J,M)$ such that,
\[l(\frac{M}{\bp{J}{n}M})= e_{HK}(J,M)(p^n)^{\textup{dim}(M)}+ O((p^n)^{\textup{dim}(M)-1}).\]   

\noindent The number $e_{HK}(J,M)$ is called the \textit{Hilbert-Kunz multiplicity} of $M$ with respect to $J$.
\end{theorem}

The Hilbert-Kunz multiplicity $e_{HK}(\mathfrak{m}, R)$ is always at least one. Smaller values of $e_{HK}(\mathfrak{m},R)$ predicts milder singularity of $(R, \mathfrak m)$; see for e.g. \cite[Cor 3.6]{AberbachEnescu}, \cite{BE}, \cite{HuYao}. It is imperative to consider Hilbert-Kunz multiplicity with respect to arbitrary ideals, for e.g. to realize $F$-signature (see \Cref{eg: when J is big})- an invariant characterizing strong $F$-regularity of $(R,\mathfrak m)$- in terms of Hilbert-Kunz multiplicity; see \cite[Cor 6.5]{TuckerPolstra}. We refer the readers to \cite{HunekeExp}, \cite[Chapter 2]{AlapanThesis} and the references there in for surveying the state of art.\\

When $(R, \mathfrak m)$ is graded, Trivedi's Hilbert-Kunz density function refines the notion of Hilbert-Kunz multiplicity:

\begin{theorem}\label{th: Hilbert-Kunz denisty after Trivedi}(see \cite{TriExist})
    Let $(R, \mathfrak m)$ be graded, $J$ be a finite colength homogeneous ideal, $M$ be a finitely generated $\mathbb Z$-graded $R$-module. Consider the sequence of functions of a real variable $s$, 
    \[\tilde g _{n,M, J}(s)=  l ([\frac{M}{\bpq{J}{q}M}]_{\flor {sq}}).\]
\begin{enumerate}
    \item There is a compact subset of $\mathbb R$ containing the supports of all $\tilde g _n$'s.
    
    \item If $\dim(M)\geq 1$, there is a function-denoted by $\tilde g _{M, J}$- such that $(\frac{1}{q})^{\dim(M)-1}\tilde g _{n,M, J}(s)$ converges pointwise to $\tilde g _{M, J}(s)$ for all $s \in \mathbb R $.
    
    \item When $\dim(M) \geq 2$, the above convergence is uniform and $\tilde g _{M,J}$ is continuous.
    \item \[e_{HK}(J, M)= \int \limits_{0}^{\infty}\tilde g _{M, J}(s)ds. \]
\end{enumerate}
\end{theorem}

\begin{definition} \label{de: HK desnity after Trivedi}
    The function $\tilde g_{M, J }$ is called the \textit{Hilbert-Kunz density function} of $(M, J)$.
\end{definition}

For a graded ring $(R,\mathfrak m)$, the \textit{Frobenius-Poincar\'e} function produces another refinement of the Hilbert-Kunz multiplicity. Frobenius-Poincar\'e functions are  essentialy a limiting function of the Hilbert series of $\frac{M}{\bpq{J}{q}M}$ in the variable $e^{-iy}$, see \cite[Rmk 3.6]{AlapanExist}.

\begin{theorem}\label{th: Frobenius-Poincare}(see \cite{AlapanExist})
    Let $(R,\mathfrak{m})$ be graded, $M$ be a finitely generated $\mathbb Z$-graded $R$-module, $J$ be a finite colength homogeneous ideal. Consider the sequence of entire functions on $\mathbb C$
    \[G_{n,M,J}(y)= (\frac{1}{q})^{\dim(M)}l([\frac{M}{\bpq{J}{q}M}]_j)e^{-iyj/q}.\]
    \begin{enumerate}
        \item The sequence of functions $G_{n,M,J}(y)$ converges to an entire function $G_{M,J}(y)$ \footnote{Note the difference in notation from \cite{AlapanExist}.} on $\mathbb C$. The convergence is uniform on every compact subset of $\mathbb C$.

        \item \[G_{M,J}(0)= e_{HK}(J,M).\]
    \end{enumerate}
\end{theorem}
The last theorem holds for any graded ring which are not necessarily standard graded. For the notion of Hilbert-Kunz density function in the non-standard graded setting, see \cite{TrivediWatanabeDomain}. By \cite[Theorem 8.3.2]{AlapanThesis}, for a standard graded $(R,\mathfrak m)$ of Krull dimension at least one, the holomorphic Fourier transform of $\tilde g_{M,J}$ is $G_{M,J}$, i.e.
\[G_{M,J}(y)= \int\limits_{0}^{\infty}\tilde g_{M,J}(s)e^{-iys}ds.\] 
Thus when $\dim(M) \geq 2$, the Hilbert-Kunz density function and the Frobenius-Poincar\'e function determine each other; see \cite[Rmk 8.2.4]{AlapanThesis}. Both Hilbert-Kunz density function and Frobenius-Poincar\'e function capture more subtle graded invariants of $(M,J)$ than the Hilbert-Kunz multiplicity. For example, when $R$ is two dimensional and normal, $J$ is generated by forms of the same degree, $\tilde g_{R,J}$ and $G_{R,J}$ determine and are determined by slopes and ranks of factors in the Harder-Narasimhan filtration of a syzygy bundle associated to $J$ on $\text{Proj}(R)$; see \cite{TriCurves}, \cite{Brrationality}, \cite[Example 3.3]{TriExist}, \cite[Chap 6]{AlapanExist}. For other results on Hilbert-Kunz density functions and Frobenius-Poincar\'e functions, see the reference section of \cite{AlapanThesis}. These two functions and the Hilbert-Kunz multiplicity of $(R,J)$ detects $J$ up to its tight closure. Recall:

\begin{definition}\label{de: tight closure}(\cite[Def 3.1]{HH})
Let $A$ be a ring of characteristic $p>0$. We say $x \in A$ is in the tight closure of an ideal $I$ if there is a $c$ not in any minimal primes of $A$ such that $c x^{p^n} \in \bp{I}{n}$ for all large $n$. The elements in the tight closure of $I$ form an ideal, denoted by $I^*$. 
\end{definition}

\begin{theorem}(\cite[Prop 5.4, Thm 5.5]{HunekeExp}, \cite[Thm 8.17]{HH})\label{th: tightclosure vs multiplicity}
Let $I \subseteq J$ be two ideals in $(R,\mathfrak m)$.
\begin{enumerate}
    \item If $I^*= J^*$, $e_{HK}(I,R)= e_{HK}(J,R)$. 
    
    \item Conversely, when $R$ is formally equidimensional, i.e. all the minimal primes of the completion $\hat R$ have the same dimension, $e_{HK}(I,R)= e_{HK}(J,R)$ implies $I^*= J^*$. Therefore, when $(R,\mathfrak{m})$ is a graded ring where all the minimal primes have the same dimension, $\tilde g_{I,R}= \tilde g_{J,R}$ or $G_{I,R}= G_{J,R}$ implies $I^*=J^*$.
\end{enumerate}
\end{theorem}

\section{$h$-function}

Given ideals $I,J$ of a local ring $(R, \mathfrak m)$ such that $I+J$ is $\mathfrak m$-primary and a finitely generated $R$-module $M$, we assign a real-valued function $h_{M,I,J}$ of a real variable, which we refer to as the corresponding \textit{$h$-function}. The existence and continuity of $h_{M,I,J}$ is proven in \Cref{sse: reduction to the domain case}. When $R$ is additionally a domain and $M=R$, given an ideal $I$ and a family of ideals $\{ J_{n}\}_{n \in \mathbb N}$- satisfying what we call \textbf{Condition C} below- in \Cref{sse: h function for families}, we associate a corresponding $h$-function which is continuous on $\mathbb R_{>0}$. Moreover, we prove basic properties of $h$-functions, which are used later in the article and provide examples. 

\subsection{$h$-functions of a domain}\label{sse: h function for families}

\begin{definition}\label{de: weak family}
 Let $I_ \bullet=\{I_n\}_{n \in \mathbb N}$ be a family of ideals of the $F$-finite local ring $R$.   

\begin{enumerate}

\item $I_ \bullet$ is called a \textit{weak $p$-family} if there exists $c \in R$ not contained in any minimal primes of maximal dimension of $R$ such that $cI_n^{[p]} \in I_{n+1}$.

\item $I_ \bullet$ is called a \textit{weak $p^{-1}$-family} if exists a nonzero $\phi \in \Hom_R(F_*R,R)$ such that $\phi(F_*I_{n+1}) \subset I_n$.

\item A  \textit{big $p$-family (resp. big $p^{-1}$-family)} is a weak $p$ (resp. $p^{-1}$)-family $I_ \bullet$ such that there is an $\alpha \in \mathbb N$ for which $\bp{\mathfrak m}{n+ \alpha} \subseteq I_n$ for all $n$.
\end{enumerate}
 
\end{definition}

A family of ideals where (1) holds with $c=1$ and $\mathfrak{m}^{[p^n]} \subseteq I_n$, has been called a $p$-family of ideals; see \cite{hernandez2018local}. Notions of $p$ and $p^{-1}$-families provide an abstract framework for proving existence of asymptotic numerical invariants:

\begin{theorem}\label{th: TP conditions}(see \cite[Theorem 4.3]{TuckerPolstra})\label{l: Tucker Polstra conditions} Let $(R,\mathfrak{m},k)$ be an $F$-finite local domain of dimension $d$, $\{I_n\}_{n \in \mathbb{N}}$ a sequence of ideals such that $\mathfrak{m}^{[p^n]} \subset I_n$ for all $n \in \mathbb{N}$.

\begin{enumerate}

\item If there exists a nonzero $c \in R$ such that $cI_n^{[p]} \subset I_{n+1}$ for all $n \in \mathbb{N}$, then $\eta=\lim_{e \to \infty}1/p^{nd}l_R(R/I_n)$ exists, and there exists a positive constant $C$ that only depends on $c$ such that $\eta-1/p^{nd}l_R(R/I_n) \leq C/p^n$ for all $n \in \mathbb{N}$.

\item If there exists a nonzero $\phi \in \Hom_R(F_*R,R)$ such that $\phi(F_*I_{n+1}) \subset I_n$ for all $e \in \mathbb{N}$, then $\eta=\lim_{n \to \infty}1/p^{nd}l_R(R/I_n)$ exists, and there exists a positive constant $C$ that only depends on $\phi$ such that $1/p^{nd}l_R(R/I_n)-\eta \leq C/p^n$ for all $n \in \mathbb{N}$.

\item If the conditions in (1) and (2) are both satisfied then there exists a constant $C$ that only depends on $c$ and $\phi$ such that $\lvert 1/p^{nd}l_R(R/I_n)-\eta\lvert \leq C/p^n$.
\end{enumerate}
\end{theorem}

\begin{lemma}\label{le: sum of weak p- families is a weak p-family} Let $(R, \mathfrak m)$ be a local domain.
Let $I_n$, $J_n$ be two weak $p$-families, then so is the family $I_n+J_n$. If $I_n$, $J_n$ are two weak $p^{-1}$-families, then so is the family $I_n+J_n$. When one of the families are big ($p$ or $p^{-1}$), then so is their sum.
\end{lemma}

\begin{proof}
Suppose there are nonzero elements $c_1,c_2$ such that $c_1I^{[p]}_n \subset I_{n+1}$ and $c_2J^{[p]}_n \subset J_{n+1}$, then $c=c_1c_2$ is still nonzero and satisfies $cI^{[p]}_n \subset I_{n+1}$, $cJ^{[p]}_n \subset J_{n+1}$. So $c(I_n+J_n)^{[p]} \subset I_{n+1}+J_{n+1}$. Suppose there are nonzero elements $ \phi_1,  \phi_2 \in \Hom_R(F_*R,R)$, such that $\phi_1(F_*I_{n+1}) \subset I_n$ and $\phi_2(F_*J_{n+1}) \subset J_n$. For $\phi \in \Hom_R(F_*R,R)$ and $r \in R$, define $F_*r\cdot\phi \in \text{Hom}_R(F_*R,R)$ by the formula $(F_*r \cdot\phi) (F_*s)= \phi(F_*(rs))$. This puts an $F_*R$-module structure on $\text{Hom}_R(F_*R,R)$, which turns out to be a torsion free module of rank one. So the $F_*R$-submodules of $\Hom_R(F_*R,R)$ generated by $\phi_1$ and $\phi_2$ have a nonzero intersection, or in other words, there exist nonzero $c_1,c_2 \in R$ and a nonzero element $\phi \in \Hom_R(F_*R,R)$ such that $\phi=\phi_1(F_*(c_1\cdot))=\phi_2(F_*(c_2\cdot))$. Thus, $\phi(F_*I_{n+1}) \subset I_n$ and $\phi(F_*J_{n+1}) \subset J_n$. So $\phi(F_*(I_{n+1}+J_{n+1})) \subset I_n+J_n$.

To prove the `big'ness, assume  that there is an $\alpha$ such that $\mathfrak{m}^{[p^{n+\alpha}]} \subseteq I_{n}$. Then we have $\mathfrak{m}^{[p^{n+ \alpha}]}  \subseteq I_{n}+ J_n$.
\end{proof}

\textbf{Condition C}, defined below, provides the right framework where we can prove existence of $h$-functions; see \Cref{th: existence of h-function for a domain}.

\noindent \textbf{Condition C:} Let $(R, \mathfrak m)$ be an $F$-finite local ring, $I$ is an ideal and $J_{\bullet}=\{J_n\}_{n \in \mathbb N}$ be a family of ideals in $R$. We say $I, J_ \bullet$ satisfies \textbf{Condition C} if 

\begin{enumerate}
    \item The family $J_{\bullet}$ is weakly $p$ and also weakly $p^{-1}$.
    \item For each real number $t$, there is an $\alpha$ such that $\bp{\mathfrak m}{\alpha+n} \subseteq I^{\ceil{tq}}+ J_n$ for all $n$.
\end{enumerate}

\noindent

\begin{definition}\label{de: h_n for a family}
Let $(R, \mathfrak m)$ be a local  ring of characteristic $p>0$. Let $I$ be an ideal and $J_{\bullet}= \{J_n\}_{n \in \mathbb N}$ be a family of ideals in $R$, such that $I+J_n$ is $\mathfrak m$-primary for all $n$. For a finitely generated $R$-module $M$ and $s \in \mathbb R$, set
\begin{enumerate}
    \item $h_{n, M, I, J_\bullet}(s)= l(\frac{M}{(I^{\ceil {sq}}+ J_n)M})$, here a nonpositive power of an ideal is the unit ideal by convention.
    \item For an integer $d$, set
    \[h_{n, M, I, J_{\bullet},d}(s)= \frac{1}{q^d}l(\frac{M}{(I^{\ceil {sq}}+ J_n)M}).\]
    \item We denote the limit of the sequence of numbers $h_{n, M, I, J_\bullet,d}(s)$, whenever it exists, by $h_{M,I, J_{\bullet}, d}(s)$. 
\end{enumerate}
Whenever one or more of the parameters $M,I , J_{\bullet}$ is clear from the context, we suppress those from $h_{n,M,I, J_\bullet}(s)$, $h_{n,M,I, J_\bullet,d}(s)$ or $h_{M,I, J_\bullet, d}(s)$. In the absence of an explicit $d$, it should be understood that $d= \dim(M).$ When $J_n= \bp{J}{n}$ for some ideal $J$, $h_{n,M,I,J}, h_{n,M,I,J,d}, h_{M,I,J}$ stand for $h_{n,M,I, J_\bullet}$, $h_{n,M,I, J_\bullet,d}$ and $h_{M,I, J_\bullet, d}$ respectively.
\end{definition}

\begin{remark}\label{re: conflicting notation}
(1) With the notational conventions and suppression of parameters declared above, $h_{n,M,I,J}$ stands for both $l(\frac{M}{(I^{\ceil {sq}}+ J_n)M})$ and $\frac{1}{q^{\dim(M)}}l(\frac{M}{(I^{\ceil {sq}}+ J_n)M})$. But in the article, it is always clear from the context what $h_{n,M,I,J}$ denotes. So we do not introduce further conventions.

\noindent (2) When $(R, \mathfrak m)$ is graded, $M,I$ and $J_{\bullet}$ are homogeneous, $h_{n,M,I,J}= h_{n,M_{\mathfrak m}, IR_{\mathfrak m}, J_{\mathfrak m}}$. So once we prove statements involving $h_n$'s in the local setting, the corresponding statements in the graded setting follow. 
\end{remark}

The following comparison between ordinary powers and Frobenius powers is used throughout this article:

\begin{lemma}\label{l: comparison of ususal and Frobenius powers}
Let $R$ be a ring of characteristic $p>0$, $J$ be a nonzero $R$-ideal generated by $\mu$ elements, $k \in \mathbb{N}$, and $q=p^n$ is a power of $p$. Then $J^{q(\mu+k-1)} \subset (J^{[q]})^k \subset J^{qk}$.
\end{lemma}

\begin{proof}The second containment is trivial. We prove the first containment. Let $J=(a_1,...,a_{\mu})$, then $J^{q(\mu+k-1)}$ is generated by $a_1^{u_1}...a_{\mu}^{u_{\mu}}$ where $\sum u_i=q(\mu+k-1)$. Let $a=a_1^{u_1}...a_{\mu}^{u_{\mu}}$, $v_i=\lfloor u_i/q \rfloor$ and $b=a_1^{v_1}...a_{\mu}^{v_{\mu}}$, then since $qv_i \leq u_i$, $b^q$ divides $a$. Now $qv_i \geq u_i-q+1$, so $\sum qv_i \geq q(\mu+k-1)+(-q+1)\mu=q(k-1)+\mu >q(k-1)$, so $\sum v_i \geq k$. This means $b \in J^k$ and $a \in J^{k[q]}=J^{[q]k}$.
\end{proof}

\begin{theorem}\label{th: existence of h-function for a domain}
Let $(R,\mathfrak{m},k)$ be an $F$-finite local domain of dimension $d$. Let $J_\bullet$ be a family of ideals such that there is a nonzero $c \in R$ and $\phi \in \textup{Hom}_R(F_*R, R)$ satisfying $cJ_n^{[p]} \subseteq J_{n+1}$ and $\phi(F_*J_{n+1})\subseteq J_n$. Let $I$ be an ideal such that for each $s \in \mathbb R$, there is an integer $\alpha$ such that $\bp{\mathfrak{m}}{n+\alpha} \subseteq I^{\ceil{sq}}+J_n$ for all $n$ (i.e. $I,J_{\bullet}$ satisfy \textbf{Condition C}). Set $I_n(s)= I^{\ceil{sq}}+J_n$. 

\begin{enumerate}
\item Fix $t \in \mathbb{R}$. Choose $\alpha \in \mathbb N$ such that $\bp{\mathfrak{m}}{n+ \alpha} \subseteq I^{\ceil{tq}}+ J_n$ for all $n$. Then there exists a positive constant $C$ depending only on $c, \phi, I, \alpha$ and independent of the specific choice of $J_{\bullet}$, such that for any $s \in (-\infty,t]$, 
$$h_{R,I, J_{\bullet},d}(s)=\lim_{n \to \infty}1/p^{nd}l_R(R/I_n(s)) \, \text{exists, and}$$ 

\begin{equation}\label{eq: rate of convergence of h function}
 |1/p^{nd}l_R(R/I_n(s))-h_{R,I, J_{\bullet},d}(s)| \leq C/p^n \, \text{for all}\, n \in \mathbb{N}.   
\end{equation}

\item Given choices $I, J_\bullet$ and $t \in \mathbb R$, one can choose $C$ depending only on $t$, such that \Cref{eq: rate of convergence of h function} holds for  $s \in [0, t]$.

\item On every bounded subset of $\mathbb{R}$, the sequence of functions $h_{n, I,J_{\bullet}, d}(s)$ converges uniformly to  $h_{R,I,J_ {\bullet}}(s).$ 
\end{enumerate}
\end{theorem}

\begin{proof}
(1) When $I=0$, $I_n(s)= J_n$, so everything follows from \Cref{th: TP conditions}.

We assume $I$ is nonzero for the rest of the proof. Note $cI_n(s)^{[p]} =cI^{\lceil sq \rceil[p]}+cJ_n^{[p]}\subseteq cI^{\lceil sq \rceil p}+cJ_n^{[p]}\subseteq I^{\lceil spq \rceil}+J_{n+1}$ as $\ceil{sq}p \geq \ceil{sqp}$. So 
\begin{equation}\label{eq: the family satisfied TP1}
cI_n(s)^{[p]} \subseteq I_{n+1}(s) \ .    
\end{equation} 
Suppose $I$ is generated by $\mu$-many elements. Then
$$I^{\lceil spq \rceil} \subseteq I^{\lceil sq \rceil p-p} \subseteq I^{[p](\lceil sq \rceil-\mu)};\,\,\,\,  \text{see \Cref{l: comparison of ususal and Frobenius powers}}.$$
Fix a nonzero $r \in (I^\mu)^{[p]}$. Then the last containment implies,

\begin{equation}\label{eq: TP2 is satisfied}
\phi(F_*r F_* I_{n+1}(s))= \phi(F_*(rI^{\ceil{spq}}))+ \phi(F_* (rJ_{n+1})) \subseteq \phi(F_*(I^{\ceil{sq}[p]}))+J_n \subseteq I_n(s) \, \text{for all}\, s \in \mathbb R \ .    
\end{equation}
\Cref{eq: the family satisfied TP1} and \Cref{eq: TP2 is satisfied} imply that, for all $s$, the nonzero elements $c \in R$ and $\phi(F_*r\cdot\dash ) \in \text{Hom}_R(F_*R, R)$ endow $I_n(s)$ with weakly $p$ and $p^{-1}$-family structures, respectively. The ideal $\bp{\mathfrak{m}}{n+ \alpha}$ is contained in $I_n(t)$ and hence in $I_n(s)$ for $s \leq t$. The rest follows by applying \Cref{th: TP conditions} to the family $I_{n+{\alpha}}(s)$ for every $s\leq t$. The feasibility of choosing $C$ depending only on $c, \phi, \alpha$ and $r$ also follows from \Cref{th: TP conditions}. Since $r \in \bpq{(I^\mu)}{p}$ can be chosen depending only on $I$, the choice of $C$ depends only on $c, \phi, \alpha$ and $I$.\\

\noindent (2) Once $I,J_n$ satisfying the hypothesis is given and $t \in \mathbb R$ is given, $c, \phi, \alpha$ can be chosen depending only on $I,J_n, t$. \\

\noindent (3) Every bounded subset of $\mathbb{R}$ is contained in some interval $(-\infty, t]$. The dependence of $C$ only on $I,J_n$ and $t$ implies (3).
\end{proof}

The domain assumption is made in the above theorem just so that we can apply \Cref{th: TP conditions}.

\begin{lemma}\label{le: compactly supported h-function}
 Suppose $I$ and $J_\bullet$ satisfy the hypothesis of \Cref{th: existence of h-function for a domain}. Suppose there is an integer $r$ such that $I^{rp^n} \subseteq J_n$ for all $n$. Then $h_{n,I,J_{\bullet}}(s)$ and $h_{I,J_{\bullet},d}$ are constant on $[r, \infty)$. 
\end{lemma}

The next two propositions produce examples of an ideal $I$ and ideal family $J_ \bullet$ satisfying \textbf{Condition C}. For specific choices of $J_{\bullet}$ and $I$, the corresponding functions $h_{I, J_\bullet,d}$ encode widely studied invariants of a prime characteristic ring such as Hilbert-Kunz multiplicity, $F$-signature, and $F$-threshold. The ring $R$ is not necessarily a domain in the next two examples.

\begin{proposition}\label{eg: J_n is big}
Let $J_\bullet$ be a family of ideals which is a big $p$ and also $p^{-1}$-family. For any ideal $I$, the pair $I$, $J_\bullet$ satisfies \textbf{Condition C}.
\end{proposition}

\begin{proof}
Since $J_{\bullet}$ is big, there is an $\alpha$ such that $\bp{\mathfrak m}{n + \alpha} \subseteq J_n$. Thus for every $s \in \mathbb R$, $\bp{\mathfrak m}{n + \alpha} \subseteq I^{\ceil {sq}}+J_n$.
\end{proof}

When $R$ is a domain, a big $p$, $p^{-1}$-family $J_{\bullet}$ thus produces an $h$-function. Thanks to \Cref{le: compactly supported h-function} such an $h_{I,J_\bullet}$ is eventually constant. 

\begin{example}\label{eg: when J is big}
 Examples of $J_\bullet$ which are both big $p$ and also $p^{-1}$ include $J_n=J^{[p^n]}$, where $J$ is an $\mathfrak m$-primary ideal. Another example of interest is when $J_n$ is the sequence of ideals defining $F$-signature of $(R,\mathfrak m)$, which we now recall. Set $p^{\alpha}= [k:k^p]$. Take  
 \[J_n= \{ x \in R \,|\, \phi(x) \in \mathfrak m, \, \text{for all}\, \phi \in \text{Hom}_R(F^n_*R,R) \}.\]
 Then $p^{\alpha n}l(R/J_n)$ coincides with the free rank of $F^n_*R$: the maximal rank of a free module $M$ such that there is an $R$-module surjection $F^n_*R \rightarrow M$; see \cite[Prop 4.5]{KevinFsignature}. The family $J_n$ is both weakly $p$ and $p^{-1}$, and $J_n$ contains $\bp{\mathfrak m}{n}$. Thanks to \Cref{th: TP conditions}, the limit
 \[s(R):= \lim_{n \to \infty}(\frac{1}{q})^{\dim(R)}l(\frac{R}{J_n})\]
exists. The number $s(R)$ measuring the asymptotic growth of the free rank of $F^n_*R$ is called the $F$-signature of $R$. The ring $(R,\mathfrak m)$ is strongly $F$-regular if and only if $s(R)$ is positive; see \cite[Theorem 0.2]{AberbachLeuschke}. When $R$ is a domain, for any nonzero ideal $I$, we have $h_{I, J_ \bullet}(s)=s(R)$ for large $s$. The continuity, left-right differentiability of such $h_{I,J_\bullet}$ are consequences of \Cref{th: one sided differentiability of h-functions}.
\end{example}

\noindent The examples of $h$-functions produced by the result below are central to extending theories of Frobenius-Poincar\'e and Hilbert-Kunz density functions to the local setting.

\begin{proposition}\label{eg: I+J is m-primary}
For any pair of ideals $I, J$ such that $I+J$ is $\mathfrak{m}$-primary, the ideal $I$ and the family $J_n= \bp{J}{n}$ satisfies \textbf{Condition C}.    
\end{proposition}

\begin{proof}
For any nonzero $c \in R$ and nonzero $\phi \in \Hom_R(F_*R, R)$, $c(J^{[p^n]})^{[p]} \subseteq \bp{J}{n+1}$ and $\phi(F_*\bp{J}{n+1}) \subset \bp{J}{n}$. So the family $J_n= \bp{J}{n}$ is weakly $p$ and $p^{-1}$  Since $I+J$ is $\mathfrak m$-primary, given a real number $s$, $\bp{\mathfrak m}{\alpha} \subseteq I^{\ceil{s}}+J$ for some $\alpha$. Then $\bp{\mathfrak{m}}{\alpha+n} \subseteq (I^{\ceil{s}}+ J)^{[p^n]} \subseteq I^{\ceil{sq}}+ \bpq{J}{q}$. So $I^{\ceil{sq}}+ \bpq{J}{q}$ is a big $p$ and $p^{-1}$-family.  
\end{proof}

For two $\mathfrak m$-primary ideals $I,J$, in \cite{TaylorInterpolation} Taylor considers $s$-multiplicity (function) which, up to multiplication by a positive number depending on $s$, coincides with the corresponding $h_{I,J}$. When $J_n= \bpq{J}{q}$, our proof of the existence of $h$-function in \Cref{th: existence of h-function for a domain} is not only different from the proof of Theorem 2.1 of \cite{TaylorInterpolation}, but also still valid when both $I$ and $J$ are not necessarily $\mathfrak{m}$-primary.  Moreover, in \Cref{th: existence of h-function for a domain}, the flexibility of choosing $C$ depending only on $\phi$ and $c$ is a byproduct of our proof; this flexibility is crucial in \Cref{th: m-adic continuity} and later.

\subsection{Growth and $\mathfrak{m}$-adic continuity of $h$-function}

Next, we examine how $h_{n,I, J_{\bullet}}(s)$ changes when the $I$ or $J_{\bullet}$ is replaced by another ideal or ideal family which is $\mathfrak m$-adically close to the initial one. The results we prove are used later in the proof of \Cref{th: m-adic continuity of h function} and in \Cref{se: relation among different function}, for example, to prove continuity of Hilbert-Kunz density function $\tilde{g}_{M,J}$ for a non $\mathfrak m$-primary ideal $J$; see \Cref{th: HK density when J is homogeneous but not of finite colength}.

\begin{lemma}Let $R$ be a noetherian local ring, $I,J$ be two $R$-ideals such that $I+J$ is $\mathfrak{m}$-primary. Let $I'$, $J'$ be two ideals such that $I \subset I'$, $J \subset J'$. Then $h_{n,M,I,J}(s) \geq h_{n,M,I',J'}(s)$.
\end{lemma}

\begin{proof}If $I \subset I'$, $J \subset J'$ then $(I^{\lceil sp \rceil}+J^{[p]})M \subset (I'^{\lceil sp \rceil}+J'^{[p]})M$, so $l(M/(I^{\lceil sp \rceil}+J^{[p]})M) \geq l(M/(I'^{\lceil sp \rceil}+J'^{[p]})M)$, which just means $h_{n,M,I,J}(s) \geq h_{n,M,I',J'}(s)$.
\end{proof}

\begin{theorem}\label{th: m-adic continuity}
Let $(R,\mathfrak{m})$ be a noetherian local ring. Assume $I,J_\bullet$  satisfies \textbf{Condition C}.

\item (1) Fix $s_0 \in \mathbb R$. We can choose a positive integer $t$ depending only on $I,J_\bullet,s_0$ such that for any ideals $J \subset \mathfrak{m}^{t}$ and $I'$ containing $I$, and all $n$, $$h_{n,M,I',J_\bullet}(s)=h_{n,M,I',J_{\bullet}+\bp{J}{n}}(s)\,\, \text{for}\,\, s \leq s_0,$$
where $J_{\bullet}+\bp{J}{n}$ denotes the family $n \mapsto J_{n}+\bp{J}{n}$

\item (2) Assume $J_{\bullet}$ is both big $p$ and $p^{-1}$-family. There exists a constant $c$ such that for any ideals $I' \subset \mathfrak{m}^t$, $t \in \mathbb{N}$ and $s \in \mathbb{R}$, $$h_{n,M,I,J_\bullet}(s-c/t) \leq h_{n,M,I+I',J_\bullet}(s) \leq h_{n,M,I,J_\bullet}(s) \leq h_{n,M,I+I_t,J_\bullet}(s+c/t).$$

\item (3) Fix $s_0>0$. There exists a positive integer $t_0$ and a positive real number $c$, both only depending on $s_0, I, J_ \bullet$ such that for any $t \geq t_0$, and any ideal $I_t \subseteq \mathfrak{m}^t$, $$h_{n,M,I,J_ \bullet}(s-c/t) \leq h_{n,M,I+I_t,J_ \bullet}(s) \leq h_{n,M,I,J_ \bullet}(s) \leq h_{n,M,I+I_t,J_ \bullet}(s+c/t),$$
for any $s \leq s_0$ and $n \in \mathbb{N}$. 
\end{theorem}

\begin{proof}\noindent (1) Let $t$ be the smallest integer such that $\mathfrak{m}^{t[q]} \subset I^{\ceil{s_0q}}+J_n$ for all $n$. By the previous lemma, it suffices to consider the case where $J=\mathfrak{m}^t$. So for $I \subseteq I'$, 
$$I'^{\lceil sq \rceil}+J_n=I'^{\lceil sq \rceil}+J_n+\mathfrak{m}^{t[q]} \, \text{for}\, s \leq s_0 \, \text{and all}\, n \in \mathbb{N},$$
proving the desired statement.\\

\item (2) Since $J_ \bullet$ is a big family, we can choose $t_0$ such that $\mathfrak{m}^{t_0[q]} \subseteq J_n$ for all $n$. We may also assume $I'=\mathfrak{m}^t$. Let $\mathfrak{m}$ be generated by $\mu$-elements, set $\epsilon_t=t_0\mu/t$. Then $\mathfrak{m}^{t\lceil \epsilon_tq \rceil} \subseteq \mathfrak{m}^{t_0\mu q} \subseteq \mathfrak{m}^{t_0[q]}  \subset J_n$ for all $n$. So $$(I+\mathfrak{m}^t)^{\lceil sq \rceil}=\sum_{0 \leq j \leq \lceil sq \rceil}I^{\lceil sq \rceil-j}\mathfrak{m}^{tj} \subset I^{\lceil sq \rceil-\lceil \epsilon_tq \rceil}+\mathfrak{m}^{t\lceil \epsilon_tq \rceil} \subset I^{\lceil sq \rceil-\lceil \epsilon_tq \rceil}+J_n \subseteq I^{\ceil{(s-t_0\mu/t)q}}+J_n$$ 
Thus we have

$$l(M/(I^{\lceil (s- t_0\mu/t)q \rceil}+J_n)M)\leq l(M/((I+\mathfrak{m}^t)^{\lceil sq \rceil}+J_n)M) \leq l(M/(I^{\lceil sq \rceil}+J_n)M).$$
So taking $c=t_0\mu$ verifies the first two inequalities. These equalities are independent of $s$, so we may replace $s$ by $s+c/t$ to get the third inequality.\\

\item (3) By (1) we can choose $t_1$ depending on $s_0,I,J_ \bullet$ so that $h_{n,M,I',J_{\bullet}+\mathfrak{m}^{t_1[q]}}(s)=h_{n,M,I',J_{\bullet}}(s)$ whenever $I \subset I'$ and $s \leq s_0 +1$. By (2), we can choose $c$ depending on $J_{\bullet}$ and $\mathfrak{m}^{t_1}$ such that
%...........
$$h_{n,M,I,J_ \bullet +\mathfrak{m}^{t_1[q]}}(s-\frac{c}{t}) \leq h_{n,M,I+I_t,J_{\bullet}+\mathfrak{m}^{t_1[q]}}(s) \leq h_{n,M,I,J_{\bullet}+\mathfrak{m}^{t_1[q]}}(s) \leq h_{n,M,I+I_t,J_{\bullet}+\mathfrak{m}^{t_1[q]}}(s+\frac{c}{t}),$$
%.........................................
for $I_t \subseteq m^t$. Take $t_0=c$. Since for $t \geq t_0$ and $s \leq s_0$, $s+ \frac{c}{t} \leq s_0+1$, the above chain of inequalities imply 
$$h_{n,M,I,J_{\bullet}}(s-c/t) \leq h_{n,M,I+I_t,J_{\bullet}}(s) \leq h_{n,M,I,J_{\bullet}}(s) \leq h_{n,M,I+I_t,J_{\bullet}}(s+c/t).$$
\end{proof}

Assertion (1) of the theorem above allows us to replace $J_{\bullet}$ by a big $p$ and $p^{-1}$-family in questions involving local structure of $h$-functions. This observation is repeatedly used later; see \Cref{th: HK density when J is homogeneous but not of finite colength}.

Next we prove that the sequence $h_{n,I, J_\bullet ,d}(s)$ is uniformly bounded on every compact subset. When $J_\bullet= \bp{J}{n}$ for some $J$, we refine the bound to show that $h_{n,I, J_\bullet, d}(s)$ is bounded above by a polynomial of degree $\dim(R/J)$ in  \Cref{th: bounds on h function special case}. The uniform (in $n$) polynomial bound on $h_n$ is used in the extension of the theory of Frobenius-Poincar\'e functions in \Cref{le: absolute convergence of the terms of sequence} and \Cref{th: exitence of Frob-Poincare for general ideals}.

\begin{lemma}\label{le: bound on the h function for a family}
In a local ring $(R, \mathfrak{m})$, let $I, J_\bullet$ satisfy \textbf{Condition C}. Let $M$ be a finitely generated $R$-module. Given $s_0 \in \mathbb{R}$, there is a constant $C$ depending only on $s_0$ such that
    $$h_{n,M,I, J_\bullet}(s) \leq Cq^d$$
for all $s \leq s_0$ and all $n \in \mathbb{N}$.
\end{lemma}

\begin{proof}
    Choose $\alpha$ such that $\bp{\mathfrak m}{n + \alpha} \subseteq I^{\ceil {s_0q}}+ J_n$. So for $s \leq s_0$,
\[h_{n,M,I, J_\bullet}(s) \leq l(\frac{M}{\bp{\mathfrak m}{n + \alpha}M}) \leq Cq^d.\]
The last ineuqality is a consequence of \cite{MonExist}.
\end{proof}

\begin{remark}\label{re: changing the residue field}
Given a noetherian local ring $(R, \mathfrak m, k)$ containing $\mathbb{F}_p$ and a field extension $k \subseteq L$, denote by $S$ the $\mathfrak m$-adic completion of $L \otimes_k \hat{R}$. Here $\hat R$ is the $\mathfrak m$-adic completion of $R$ which can be treated as a $k$-algebra thanks to the existence of a coefficient field of $\hat R$; see \cite[tag 0323]{stacks-project}. The residue field of the local ring $S$ is isomorphic to $L$. The natural map $R \rightarrow S$ is faithfully flat. Now given a finite length $R$-module $M$, $l_R(M)= l_S(S \otimes_R M)$. We use this observation to make simplifying assumption on the residue field of $R$.
\end{remark}

\begin{theorem}\label{th: bounds on h function special case}Let $(R,\mathfrak{m})$ be a noetherian local ring of dimension $d$, $I,J$ be two $R$-ideals such that $I+J$ is $\mathfrak{m}$-primary. Assume $I$ is generated by $\mu$ elements, $M$ is generated by $\nu$ elements, and $d'=\dim R/J$. Then:
\begin{enumerate}
\item There exist a polynomial $P_1(s)$ of degree $d'$  such that for any $s \geq 0$,
$$\frac{l(M/(I^{\lceil sq \rceil}+J^{[q]})M)}{l(R/\mathfrak{m}^{[q]})} \leq P_1(s).$$
Moreover if $d'>0$, the leading coefficient of $P_1$ can be taken to be $\frac{\nu e(I,R/J)}{d'!}$ 
\item There exist a polynomial $P_2(s)$ such that
$$\frac{l(M/(I^{\lceil sq \rceil}+J^{[q]})M)}{q^d} \leq P_2(s).$$
In other words, $h_{n,M,d}(s)/q^d \leq P_2(s)$.
\item There exists a polynomial $P_3$ of degree $d'$ and leading coefficient $\frac{\nu e(I,R/J)e_{HK}(R)}{d'!}$ such that for any $s \geq 0$,
$$\overline{\lim}_{n \to \infty}\frac{l(M/(I^{\lceil sq \rceil}+J^{[q]})M)}{q^d} \leq P_3(s).$$
\end{enumerate}
\end{theorem}
\begin{proof} We may assume that the residue field is perfect using \Cref{re: changing the residue field}.

\begin{enumerate}
\item Suppose $M$ is generated by $\nu$ many elements. Then
\begin{align*}
l(M/(I^{\ceil{sq}}+J^{[q]})M) \leq \nu l(R/I^{\ceil{sq}}+J^{[q]})\\
\leq \nu l(R/(I^{\ceil s })^{[q]}+J^{[q]})\\
\leq \nu l(F^n_*R/(I^{\ceil s}+J)F^n_*R)\\
\leq \nu \mu_R(F^n_*R)l(R/I^{\ceil s }+J)
\end{align*}
Let $P_0$ be the  Hilbert-Samuel polynomial for the $I$-adic filtration on $R/J$; $P_0$ has degree $d'$ and leading coefficient $\frac{\nu e(I,R/J)}{d'!}>0$. Since $P_0$ is a polynomial with positive leading coefficient, we can fix $s_0$ such that for $s \geq s_0$, $l(R/I^{\ceil s}+J)=P_0(\ceil s)$ and $P_0$ is increasing. Thus for $s\geq s_0$,
$$l(R/I^{\lfloor s \rfloor}+J)\leq P_0(s+1).$$
When $R/J$ has Krull dimension zero, we can take the constant polynomial $P_1(s)=l(R/J)<\infty$, then $l(R/I^{\ceil s}+J)\leq P_1(s)$ for all $s$. When $R/J$ has positive Krull dimension, we can add a suitable positive constant to $P_0(s+1)$ to get a $P_1$ so that $l(R/I^{\lfloor s \rfloor}+J)\leq P_1(s)$ on $[0,s_0+1]$ and thus on $\mathbb R$.

\item Since $\lim_{n \to \infty}l(R/\mathfrak{m}^{[q]})/q^d=e_{HK}(R)$ exists,
$$C=\sup_n l(R/\mathfrak{m}^{[q]})/q^d$$
exists. So for any $n$, $l(R/\mathfrak{m}^{[q]})/q^d \leq C$, and $P_2=CP_1$ satisfies (2).

\item
\begin{align*}
\overline{\lim}_{n \to \infty}\frac{l(M/(I^{\lceil sq \rceil}+J^{[q]})M)}{q^d}\\
\leq \overline{\lim}_{n \to \infty}\frac{l(M/(I^{\lceil sq \rceil}+J^{[q]})M)}{l(R/\mathfrak{m}^{[q]})}\overline{\lim}_{n \to
 \infty}\frac{l(R/\mathfrak{m}^{[q]})}{q^d}\\
 \leq e_{HK}(R)P_1(s).
\end{align*}
So $P_3=e_{HK}(R)P_1$ works.
\end{enumerate}
\end{proof}

\subsection{Lipschitz continuity of $h$-functions: Application of a `convexity technique'}

Proving continuity of $h_{R,I, J_{\bullet}}$- when $R$ is a domain is more involved than proving its existence. In this subsection, we develop results aiding the proof of local Lipschitz continuity of $h_{R,I, J_{\bullet}}$ proven in \Cref{th: Lipschitz continuity for a family}. When $J_n= \bpq{J}{q}$, these results are used to prove existence and continuity of the $h$-function of a finitely generated module in \Cref{th: uniform convergence of h function for modules}, by reducing the problem to the case where $R$ is reduced. The key result which allows these reduction steps is \Cref{th: bound on the density function}. We prove \Cref{th: bound on the density function} by utilizing the monotonicity of a certain numerical function. This technique of using the monotonicity- which we call the `convexity technique'- is repeatedly used later for example to prove left and right differentiability of the $h$-function among other properties. The required monotonicity result appears in \Cref{l: Boij's theorem strengthened}. This is an adaptation and generalization of Boij-Smith's result in \cite{BoijGregory15}  which is suitable for our purpose.\\

\begin{lemma}\label{l: Boij's theorem strengthened}
Let $(R, \mathfrak m)$ be a noetherian local ring, $I$ be an $\mathfrak{m}$-primary ideal generated by $\mu$ elements, $M$ be a finitely generated $R$-module, $S$ be the polynomial ring of $\mu$-variables over $R/\mathfrak m$. Then the function $i \to l(I^iM/I^{i+1}M)/l(S_i)$ is decreasing for $i \geq 0$.
\end{lemma}

\begin{proof}Consider the associated graded ring $\textup{gr}_I(R)$. Since $I$ is generated by a set of $\mu$ elements, as a graded ring $\textup{gr}_I(R)$ is a quotient of the standard graded polynomial ring $R/I[T_1,...,T_\mu]$ over $R/I$. Recall $S=\frac{R}{\mathfrak m}[T_1,...,T_\mu]$. Since $M/IM$ is Artinian, there exists a filtration $$0=N_0 \subset N_1 \subset ... \subset N_l=M/IM,\,\text{such that} \, N_{j+1}/N_j=R/\mathfrak{m}\, \, \text{for}\, 0 \leq j \leq l-1.$$ 
Let $M_j$ be the $\textup{gr}_I(R)$-submodule of $\textup{gr}_I(M)$ spanned by $N_j$. Then $M_{j+1}/M_j$ is annihilated by $\mathfrak{m}\textup{gr}_I(R)$. So it is naturally a $\textup{gr}_I(R)/\mathfrak{m}\textup{gr}_I(R)$-module, hence is an $S$-module, and it is generated in degree $0$. So by Theorem 1.1 of \cite{BoijGregory15}, for any $i \geq 0$, $$l(M_{j+1}/M_j)_i/l(S_i) \geq l(M_{j+1}/M_j)_{i+1}/l(S_{i+1}).$$ Since truncation at degree $i$ is an exact functor from $\textup{gr}_I(R)$-modules to $R$-modules, taking sum over $0 \leq j \leq l-1$ we get $l(M_l)_i/l(S_i) \geq l(M_l)_{i+1}/l(S_{i+1})$. Since $M_l=\textup{gr}_I(R)N_l=\textup{gr}_I(M)$, we are done.
\end{proof}

\noindent When $I$ is a principal ideal, the above lemma manifests into the following easily verifiable result. 

\begin{example}Let $R$ be a noetherian local ring, $f$ be an element in $R$ such that $R/fR$ has finite length. Then for any $j
\geq i \geq 0$, $l(f^iR/f^{i+1}R) \geq l(f^jR/f^{j+1}R)$. This means that the function $i \to l(R/f^iR)$ is convex on $\mathbb N$; see \Cref{de: convex function}.
\end{example}

\begin{theorem}\label{th: bound on the density function}Let $R$ be a noetherian local ring, $M$ be a finitely generated module of dimension $d$. Suppose $I,J_\bullet$ satisfies \textbf{Condition C}. Fix $0<s_1<s_2<\infty \in \mathbb{R}$. Then there is a constant $C$ and a power $q_0=p^{n_0}$ that depend on $s_1,s_2$, but independent of $n$ such that for any $s_1 \leq s \leq s_2-1/q$ and $q \geq q_0$
$$l(\frac{(I^{\lceil sq \rceil}+J_n)M}{(I^{\lceil sq \rceil+1}+J_n)M}) \leq Cq^{d-1}$$
In other words, whenever $s_1 \leq s \leq s_2-1/q$ and $q \geq q_0$,
$$|h_{n,M}(s+1/q)-h_{n,M}(s)| \leq Cq^{d-1}.$$
\end{theorem}

\begin{proof}We may assume $s_1,s_2 \in \mathbb{Z}[1/p]$. Otherwise, since $\mathbb{Z}[1/p]$ is dense in $\mathbb{R}$, we can choose $s'_1 \in (0,s_1) \cap\mathbb{Z}[1/p]$, $s'_2 \in (s_2,\infty) \cap\mathbb{Z}[1/p]$ and replace $s_1,s_2$ by $s'_1,s'_2$. Choose $s_3 \in \mathbb{Z}[1/p]$ such that $0<s_3<s_1$ and choose $q_0$ such that $s_1q_0,s_2q_0,s_3q_0 \in \mathbb{Z}$. Let $I$ be generated by a set of $\mu$ many elements. Applying \Cref{l: Boij's theorem strengthened} to the module $M/J_nM$ we know for any $0 \leq t \leq \lceil sq \rceil$,

$$\frac{l(\frac{I^{\ceil{sq}}(M/J_nM)}{I^{\ceil{sq}+1}(M/J_nM)})}{{{\mu+ \lceil sq \rceil -1}\choose{\mu-1}}} \leq \frac{l(\frac{I^{t}(M/J_nM)}{I^{t+1}(M/J_nM)})}{{{\mu+t-1}\choose{\mu-1}}}.$$

\noindent Rewritten, the above inequality yields

$$\frac{l(\frac{(I^{\lceil sq \rceil}+J_n)M}{(I^{\lceil sq \rceil+1}+J_n)M})}{{{\mu+ \lceil sq \rceil -1}\choose{\mu-1}}} \leq \frac{l(\frac{(I^t+J_n)M}{(I^{t+1}+J_n)M})}{{{\mu+t-1}\choose{\mu-1}}}.$$
Thus for $s_1 \leq s \leq s_2- \frac{1}{q}$ and $q \geq q_0$,

\begin{equation*}
\begin{split}
(\lceil sq \rceil-s_3q)l(\frac{(I^{\lceil sq \rceil}+J_n)M}{(I^{\lceil sq \rceil+1}+J_n)M}) &\leq {{\mu+ \ceil{sq} -1}\choose{\mu-1}}\sum^{\lceil sq \rceil-1}_{t=s_3q}\frac{l(\frac{(I^t+J_n)M}{(I^{t+1}+J_n)M})}{{{\mu+t-1}\choose{\mu-1}}} \\
&\leq \frac{{{\mu+\ceil{sq}-1}\choose{\mu-1}}}{{{\mu+s_3q-1}\choose{\mu-1}}}l(\frac{(I^{\lceil sq \rceil}+J_n)M}{(I^{s_3q}+J_n)M}) \\
&\leq \frac{{{\mu+\ceil{sq}-1}\choose{\mu-1}}}{{{\mu+s_3q-1}\choose{\mu-1}}}[l(\frac{M}{(I^{\lceil sq \rceil}+J_n)M})-l(\frac{M}{(I^{s_3q}+J_n)M})] \\
&\leq \frac{{{\mu+s_2q-1}\choose{\mu-1}}}{{{\mu+s_3q-1}\choose{\mu-1}}}[l(\frac{M}{(I^{s_2q}+J_n)M})-l(\frac{M}{(I^{s_3q}+J_n)M})].
\end{split}
\end{equation*}

\noindent Therefore for $s_1 \leq s \leq s_2- \frac{1}{q}$ and $q \geq q_0$,

$$l(\frac{(I^{\lceil sq \rceil}+J_n)M}{(I^{\lceil sq \rceil+1}+J_n)M}) \leq \frac{1}{s_1q-s_3q}\frac{{{\mu+s_2q-1}\choose{\mu-1}}}{{{\mu+s_3q-1}\choose{\mu-1}}}[l(\frac{M}{(I^{s_2q}+J_n)M})-l(\frac{M}{(I^{s_3q}+J_n)M})].$$
By \Cref{le: bound on the h function for a family}, we can choose a constant $C'$ depending only on $s_2$ such that for $s \leq s_2$,
\[l(\frac{M}{(I^{sq}+J_n)M}) \leq C'q^d .\]
Since ${{\mu+s_2q-1}\choose{\mu-1}}/{{\mu+s_3q-1}\choose{\mu-1}}$ is bounded above by a constant depending on $s_1, s_2, s_3$ and $s_3$ depends only on $s_2$, we can choose $C$ depending only on $s_1, s_2$ such that for all $n$ and $q \geq q_0,$

\[l(\frac{(I^{\lceil sq \rceil}+J_n)M}{(I^{\lceil sq \rceil+1}+J_n)M}) \leq Cq^{d-1}.\]

Here $C$ is a constant only depending on $s_1,s_2,s_3$, and $s_3$ depends only on $s_1$.
\end{proof}

\noindent Therefore, whenever $h_{M,I, J_{\bullet}}$ exists, it is locally Lipschitz continuous away from zero. 

\begin{theorem}\label{th: Lipschitz continuity for a family}
Let $I$ be an ideal and $J_\bullet$ be a family of ideals satisfying \textbf{Condition C} in a domain $(R, \mathfrak m)$ of Krull dimension $d$. Given real numbers $0<s_1< s_2$, there is a constant $C$ depending only in $s_1,s_2$ such that for any $x, y \in [s_1, s_2]$,
\[|h_R(x)-h_R(y)| \leq C|x-y|\]
\end{theorem}

\begin{proof}
 Given $s_1,s_2$ as above and $x,y$ in $[s_1,s_2]$, by \Cref{th: bound on the density function}, we can choose a constant $C$ depending only on $s_1, s_2$ such that
\[|h_{n,R}(x)- h_{n,R}(y)|= |h_{n,R}( \frac{\ceil {qx}}{q})- h_{n,R}(\frac{\ceil {qy}}{q})| \leq C|\frac{\ceil {qx}}{q} - \frac{\ceil {qy}}{q}|q^{d} \, \text{for all}\, n .\]
Divide both sides by $q^{d}$ and take limit as $n$ approaches infinity. Since for any real number $s$, $\frac{h_n(s)}{q^d}$ and $\ceil {qs} /q$ converge to $h_R(s)$ and $s$ respectively, 
\[|h_R(x)-h_R(y)| \leq C|x-y|.\]   \end{proof}

\begin{lemma}\label{le: reduction to the reduced case}
Assume the residue field of $R$ is perfect and $M$ is a module of dimension $d$. For each integer $n_0 \geq 0$ and fixed $0<s_1<s_2<\infty \in \mathbb{R}$, there is a constant $C$ independent of $n$ such that
$$|h_{n+n_0,M,I,J}(s)-h_{n,F^{n_0}_*M,I,J}(s)| \leq Cq^{d-1}$$
for any $s_1 \leq s \leq s_2$.
\end{lemma}

\begin{proof}For any $q_0$, $\lceil sqq_0 \rceil \leq \lceil sq \rceil q_0 \leq \lceil sqq_0 \rceil+q_0$. We have,

\begin{gather*}
|h_{n+n_0,M,I,J}(s)-h_{n,F^{n_0}_*M,I,J,d}(s)|\\
=|l(M/(I^{\lceil sqq_0 \rceil}+J^{[qq_0]})M)-l(F^{n_0}_*M/(I^{\lceil sq \rceil}+J^{[q]})F^{n_0}_*M)|\\
=|l(M/(I^{\lceil sqq_0 \rceil}+J^{[qq_0]})M)-l(M/(I^{\lceil sq\rceil[q_0]}+J^{[qq_0]})M)|\\
=(l(I^{\lceil sqq_0 \rceil}+J^{[qq_0]})M/(I^{\lceil sq\rceil q_0}+J^{[qq_0]})M)+l(I^{\lceil sq \rceil q_0}+J^{[qq_0]})M/(I^{\lceil sq\rceil[q_0]}+J^{[qq_0]})M)) \ .
\end{gather*}

\noindent Note that $1/q_0\ceil{sqq_0} \geq sq \geq \ceil{sq}-1$, so $\ceil{sq}q_0 \leq \ceil{sqq_0}+q_0$, so $I^{\ceil{sqq_0}+q_0} \subset I^{\ceil{sq}q_0}$. 

\noindent Suppose $I$ is generated by $\mu$ elements, then by \Cref{l: comparison of ususal and Frobenius powers}, $I^{\lceil sq \rceil q_0} \subset I^{(\lceil sq \rceil-\mu+1) [q_0]}$. 
Now by \Cref{th: bound on the density function}, we can choose a constant $C$ depending only on $s_1, s_2$ but independent of $q$ such that for all $s \in [s_1, s_2]$, 

\begin{gather*}
l(\frac{(I^{\lceil sqq_0 \rceil}+J^{[qq_0]})M}{(I^{\lceil sq\rceil q_0}+J^{[qq_0]})M})+l(\frac{(I^{\lceil sq \rceil q_0}+J^{[qq_0]})M}{(I^{\lceil sq\rceil[q_0]}+J^{[qq_0]})M})\\
\leq l(\frac{(I^{\lceil sqq_0 \rceil}+J^{[qq_0]})M}{(I^{\lceil sqq_0\rceil+q_0}+J^{[qq_0]})M})+ l(\frac{(I^{(\lceil sq \rceil-\mu+1) [q_0]}+J^{[qq_0]})M}{(I^{\lceil sq\rceil[q_0]}+J^{[qq_0]})M}) \leq Cq^{d-1} \ .
\end{gather*}
\end{proof}
The lemma above allows us to replace $M$ by $F^{n_0}_*M$ in the proof of the existence of $h_{M,I,J,d}$. Since we may replace $R$ by $R/\text{ann}(F^{n_0}_*M)$ and for large enough $n_0$, $\text{ann}(F^{n_0}_*M)$ contains the nilradical of $R$, we may assume $R$ is reduced while proving the existence of $h_{M,I,J}$.

\begin{corollary}\label{c: replacing by frobenius pushforward} Assume the residue field of $R$ is perfect. For each $n_0 \geq 0$, $h_{M,I,J,d}(s)$ exists if and only if $h_{F^{n_0}_*M,I,J,d}(s)$ exists, and if they both exist then
$$q_0^dh_{M,I,J,d}(s)=h_{F^{n_0}_*M,I,J,d}(s).$$
\end{corollary}

\subsection{$h$-function of a module}\label{sse: reduction to the domain case}
For a noetherian local ring $(R, \mathfrak m)$, $R$-ideals $I,J$ such that $I+J$ is $\mathfrak m$-primary and a finitely generated $R$-module $M$, we prove the existence of $h_{M,I,J}$ in \Cref{th: uniform convergence of h function for modules} and prove the local Lipschitz continuity of $h_{M,I,J}$ in \Cref{pr: cotinuity of $h$ functions}. First, we prove preparatory results to reduce the problem of the existence of $h_{M,I,J}$ to the situation where $M=R$ and $R$ is a domain. Recall:
\begin{definition}\label{de: Assh}Set $\Assh R=\{P \in \Spec R: \dim R=\dim R/P\}$.
\end{definition}

\begin{lemma}\label{le: generic dependence}\cite[Proof of Lemma 1.3]{MonExist}If $M,N$ are two $R$-modules such that $M_P \cong N_P, \forall P \in \Assh R$. Then there is an exact sequence
$$0 \to N_1 \to M \to N \to N_2 \to 0$$
such that $\dim N_1, \dim N_2 \leq \dim(R)-1$. Moreover it breaks up into two short exact sequences:
$$0 \to N_1 \to M \to N_3 \to 0$$
$$0 \to N_3 \to N \to N_2 \to 0.$$
\end{lemma}
\begin{lemma}\label{le: length inequality after quotient modulo pricipal ideal}
Let $N \subset M$ be two $R$-modules of finite length, and take $a \in R$, then $l(M/aM) \geq l(N/aN)$.
\end{lemma}
\begin{proof}Consider the commutative diagram,

\begin{center}
  \begin{tikzcd}
0 \arrow[r] & 0:_Na \arrow[r] \arrow[d] & N \arrow[r, "a"] \arrow[d] & N \arrow[r] \arrow[d] & \frac{N}{aN} \arrow[r] \arrow[d] & 0 \\
0 \arrow[r] & 0:_Ma \arrow[r]           & M \arrow[r, "a"]           & M \arrow[r]           & \frac{M}{aM} \arrow[r]           & 0
\end{tikzcd}  
\end{center}

We see the map $0:_Na \to 0:_Ma$ is injective. By the additivity of length on exact sequences we see $l(M/aM)=l(0:_Ma) \geq l(0:_Na)=l(N/aN)$.
\end{proof}

\begin{lemma}\label{l: growth Tor_1 improved version} Let $(R,\mathfrak{m},k)$ be a local ring of dimension $d$. Suppose $I,J_{\bullet}$ satisfies \textbf{Condition C}, and $M$ is a module of dimension $d' \leq d-1$. Fix $s_0 \in \mathbb{R}$. Then there are constants $C_1, C_2$ depending on $s_0$ but independent of $n$ such that $l(\Tor^R_0(R/(I^{\lceil sq \rceil}+J_n),M)) \leq C_1q^{d-1}$ and $l(\Tor^R_1(R/(I^{\lceil sq \rceil}+J_n),M)) \leq C_2q^{d-1}$ for any $s \leq s_0$. Moreover if $J_\bullet$ is big, $C_1,C_2$ can be chosen independently of $s$. 
\end{lemma}

\begin{proof}
Since $I, J_\bullet$ satisfies \textbf{Condition C}, we can find an $\mathfrak{m}$-primary ideal $J$ such that for $s \leq s_0$, $\bpq{J}{q} \subseteq I^{\ceil{sq}}+ J_n$ for all $n$.
As $M/J^{[q]}M$ surjects onto $\Tor^R_0(R/(I^{\lceil sq \rceil}+J_n),M)$, 
 and we can find a constant $C_1$, such that $l(M/J^{[q]}M) \leq C_1q^{\dim M}$, $l(\text{Tor}_0^R(R/I^{\ceil{sq}}+\bpq{J}{q},M)) \leq C_1q^{d-1}$.

To see the bound on $\textup{Tor}_1$, for a fixed $s \leq s_0$, consider the exact sequence:
$$0 \to (I^{\lceil sq \rceil}+J_n)/J^{[q]} \to R/J^{[q]} \to R/(I^{\lceil sq \rceil}+J^{[q]}) \to 0$$
So by the long exact sequence of Tor, it suffices to show that we can choose $C_2$ satisfying $$l(\Tor^R_1(R/J^{[q]},M)) \leq C_2q^{d-1} \, \text{and} \, l(\frac{I^{\lceil sq \rceil}+J_n}{J^{[q]}} \otimes M) \leq C_2q^{d-1}.$$ 
Choosing a $C_2$ satisfying the first inequality is possible thanks to \cite[Lemma 1.1]{SecCoef}. For the remaining inequality, by taking a prime cyclic filtration of $M$, we may assume $M=R/P$ for some $P \in \Spec(R)$ with $\dim M \leq \dim R-1$. In this case, $P \notin \Assh(R)$. So we can choose $b \in P$ such that $\dim R/bR \leq \dim R-1$. Apply \Cref{le: length inequality after quotient modulo pricipal ideal} to $I^{\lceil sq \rceil}+J^{[q]}/J^{[q]} \subset R/J^{[q]}$, we see that we can enlarge $C_2$ independently of $s$ and $q$ so that
\begin{gather*}
  l(\frac{I^{\lceil sq \rceil}+J_n}{J^{[q]}} \otimes_R R/P) \leq l(\frac{I^{\lceil sq \rceil}+J_n}{J^{[q]}} \otimes_R R/bR)\\ \leq l(R/J^{[q]}\otimes_R R/bR)=l(R/bR+J^{[q]}) \leq C_2q^{d-1}.  
\end{gather*}
If $J_{\bullet}$ is big then for large $s$, $I^{\ceil{sq}} \subset J_n$ for all $n$. So we can find constant $D$ such that for every $n$, $l(\Tor^R_0(R/(I^{\lceil sq \rceil}+J_n),M))$ and $l(\Tor^R_1(R/(I^{\lceil sq \rceil}+J_n),M))$ are constant for $s \geq D$. So $C_1,C_2$ can be chosen independent of $s \in \mathbb{R}$.
So we are done.
\end{proof}

\begin{lemma} \label{le: h is determinde by generic behaviour}Let $M,N$ be two finitely generated $R$-modules that are isomorphic at $P \in \Assh R$. Then for any $t>0$, there is a constant $C$ depending on $M,I,J,t$ but independent of $n$ such that for any $s<t$
$$|h_{n,M,d}(s)-h_{n,N,d}(s)| \leq C/q$$
Moreover if $J$ is $\mathfrak{m}$-primary, then $C$ can be chosen independently of $t$.
\end{lemma}
\begin{proof}
By \Cref{le: generic dependence}, there is an exact sequence
$$0 \to N_1 \to M \to N \to N_2 \to 0$$
such that $\dim N_1, \dim N_2 \leq d-1$. And it breaks up into two short exact sequences:
$$0 \to N_1 \to M \to N_3 \to 0$$
$$0 \to N_3 \to N \to N_2 \to 0$$
Now by the long exact sequence of Tor we get
$$|l(M/(I^{\lceil sq \rceil}+J^{[q]})M)-l(N_3/(I^{\lceil sq \rceil}+J^{[q]})N_3)| \leq l(N_1/(I^{\lceil sq \rceil}+J^{[q]})N_1),$$

\begin{equation*}
    \begin{split}
        & |l(\frac{N_3}{(I^{\lceil sq \rceil}+J^{[q]})N_3})-l(\frac{N}{(I^{\lceil sq \rceil}+J^{[q]})N})|\\
        & \leq l(\frac{N_2}{(I^{\lceil sq \rceil}+J^{[q]})N_2})+l(\Tor_1^R(\frac{R}{I^{\lceil sq \rceil}+J^{[q]}},N_2)) \ .
    \end{split}
\end{equation*} 

Thus by \Cref{l: growth Tor_1 improved version}, there is a constant $C$ such that
$$|l(M/(I^{\lceil sq \rceil}+J^{[q]})M)-l(N/(I^{\lceil sq \rceil}+J^{[q]})N)| \leq Cq^{d-1},$$
and divide both sides by $q^d$ to get the conclusion.
\end{proof}

\begin{lemma}\label{le: uniform Cauchyness}
Let $(R,\mathfrak{m},k)$ be a local ring, $I,J$ be two ideals such that $I+J$ is $\mathfrak{m}$-primary, and $M$ be a finitely generated $R$-module of dimension $d$. For any $0<s_1<s_2<\infty$, there is a constant $C$ depending on $M,I,J,s_1,s_2$ but independent of $n$ such that for any $s_1 \leq s \leq s_2$
$$|h_{n+1,M,d}(s)-h_{n,M,d}(s)| \leq C/q$$
\end{lemma}
\begin{proof}
We may assume that the residue field is perfect using \Cref{re: changing the residue field}. Choose sufficiently large $n_0$ such that $R/\ann(F^{n_0}_*M)$ is reduced. The positive constants $C_1, C_2, C_3$ chosen below depends only on $M,I,J,s_1,s_2$ and is independent of $n$. By \Cref{le: reduction to the reduced case},
$$|h_{n+n_0,M,I,J}(s)-h_{n,F^{n_0}_*M,I,J}(s)| \leq C_1q^{d-1}$$
and
$$|h_{n+n_0+1,M,I,J}(s)-h_{n+1,F^{n_0}_*M,I,J}(s)| \leq C_1q^{d-1}$$
So it suffices to prove existence of a suitable $C$ such that
$$|h_{n+1,F^{n_0}_*M,d}(s)-h_{n,F^{n_0}_*M,d}(s)| \leq C/q.$$
Replacing $M$ by $F^{n_0}_*M$ and $R$ by $R/\ann(F^{n_0}_*M)$, so we may assume $R$ is reduced. In this case it suffices to prove
$$|h_{n+1,M,I,J}(s)-h_{n,F_*M,I,J}(s)| \leq C_2q^{d-1}.$$
Thanks to the reducedness of $R$, $R_P=R_P/PR_P$ is a field and $M_P$ is free over $R_P$ for any $P \in \Assh(R)$. Applying equation 2.2 of \cite{Kunz76} to the domain $R/P$, we see the localizations of $M^{\oplus p^d}$ and $F_*M$ are isomorphic at all $P \in \Assh R$. So by \Cref{le: h is determinde by generic behaviour}, 
$$|h_{n,F_*M,I,J}(s)-p^dh_{n,M,I,J}(s)| \leq C_3q^{d-1}.$$
Thus one can choose a $C'$ which depends only on $M,I,J, s_1, s_2$ such that for all $s \in [s_1,s_2]$ and $n \in \mathbb N$,
\[|h_{n+1,M,I,J}(s)-p^dh_{n,M,I,J}(s)| \leq C'q^{d-1}.\]
Dividing by $(pq)^d$ and letting $C=C'/p^d$, we get
$$|h_{n+1,M,I,J,d}(s)-h_{n,M,I,J,d}(s)| \leq C/q.$$
\end{proof}

\begin{theorem}\label{th: uniform convergence of h function for modules}Let $(R,\mathfrak{m},k)$ be a noetherian local ring, $I,J$ be two $R$-ideals such that $I+J$ is $\mathfrak{m}$-primary, and $M$ is a finitely generated $R$-module. Then for every $s \in \mathbb R$, $$\lim_{n \to \infty}\frac{1}{q^{\dim(M)}}h_{n,M,I,J}(s)=h_{M,I,J}(s)$$ exists. Moreover the convergence is uniform on $[s_1,s_2]$ for any $0<s_1<s_2<\infty$.
\end{theorem}

\begin{proof}
    By replacing $R$ by $R/ \text{ann}(M)$, we may assume $\dim(M)= \dim(R)$.  Given $s_1, s_2$ as in the statement, it follows from \Cref{le: uniform Cauchyness} that $h_{n,M,I,J}(s)/q^{\dim(M)}$ is uniformly Cauchy on $[s_1, s_2]$. So the theorem follows.
\end{proof}

We also have:

\begin{theorem}\label{pr: cotinuity of $h$ functions}
Let $(R,\mathfrak{m},k)$ be a local ring of dimension $d$, $I$, $J$ be two $R$-ideals such that $I+J$ is $\mathfrak{m}$-primary, and $M$ be a finitely generated $R$-module. Then:

\begin{enumerate}
\item $h_M(s)$ is Lipschitz continuous on $[s_1,s_2]$ for any $0<s_1<s_2<\infty$. Consequently, it is continuous on $(0,\infty)$.

\item $h_M(s)$ is increasing \footnote{In fact, $h_M(s)$ is strictly increasing until a point and then possibly becomes constant; see \Cref{c: minimal stable point is the support}}. It is $0$ on $(-\infty,0]$. It is continuous if and only if it is continuous at $0$, if and only if $\lim_{s\to 0^+}h_M(s)=0$. The limit $\lim_{s\to 0^+}h_M(s)$ always exists and is nonnegative.

\item Assume $J$ is $\mathfrak{m}$-primary. Then for $s >>0$, $h_{n,M}(s)=e_{HK}(J,M)$ is a constant. 

\item There is a polynomial $P(s)$ of degree $\dim R/J$ such that $h_M(s) \leq P(s)$ on $\mathbb R$.
\end{enumerate}
\end{theorem}

\begin{proof}
\begin{enumerate}
\item An argument similar to that in the proof of \Cref{th: Lipschitz continuity for a family} with $R$ replaced by $M$ and $J_n= \bpq{J}{q}$ yields a proof. The difference is that when $J_n= \bpq{J}{q}$, we know the existence of $h_{M,I,J}$.   

\item If $s_1 \leq s_2$, then $\ceil{s_1q} \leq \ceil{s_2q}$, so $I^{\ceil{s_2q}} \subset I^{\ceil{s_1q}}$. This implies 
$$l(M/(I^{\ceil {s_1q}}+J^{[q]})M) \leq l(M/(I^{\ceil {s_2q}}+J^{[q]})M),$$
which is just
$$h_{n,M}(s_1) \leq h_{n,M}(s_2).$$
So after dividing $p^{n\dim M}$ and let $n \to \infty$, we get $h_M(s_1) \leq h_M(s_2)$. This implies $h_M(s)$ is increasing; so in particular the limit $\lim_{s\to 0^+}h_M(s)$ always exists and is at least $h_M(0)$. If $s \leq 0$, then $\ceil{sq} \leq 0$, so $I^{\ceil{sq}}=R$. Thus $M/(I^{\ceil {s_1q}}+J^{[q]})M=0$ and $h_{n,M}(s)=0$ for any $n$, so $h_M(s)=0$. So $h_M(s)$ is continuous on $(-\infty,0)$ and $(0,\infty)$, and $\lim_{s\to 0^-}h_M(s)=0=h_M(0)$, so we get (2).

\item Let $J$ be generated by $\mu$ elements. For $s>>0$, $I^{\lfloor s/\mu \rfloor} \subset J$. So $I^{\lceil sq \rceil} \subset I^{\lfloor s/\mu \rfloor q\mu} \subset J^{q\mu} \subset J^{[q]}$, so $h_{n,M}(s)=l(M/J^{[q]}M)$ and $h_M(s)=\lim_{n \to \infty}\frac{l(M/J^{[q]}M)}{q^d}=e_{HK}(J,M)$.

\item This is a corollary of \Cref{th: bounds on h function special case} and \Cref{th: uniform convergence of h function for modules}.
\end{enumerate}
\end{proof}

Now we prove that the map sending ideal pairs $(I,J)$ to $h_{I,J}$ is $\mathfrak{m}$-adically continuous in $I,J$. This result is used to transform questions about $h_{I,J}$, where $I,J$ are not necessarily $\mathfrak{m}$-primary to questions about $h_{I,J}$, where $I,J$ are $\mathfrak{m}$-primary; see \Cref{th: HK density when J is homogeneous but not of finite colength} for example.

\begin{theorem} \label{th: m-adic continuity of h function} Let $(I_t)_t$, $(J'_t)_t$ be two sequences of ideals of $(R, \mathfrak{m})$ such that $I_t+J'_t \subset\mathfrak{m}^t$ for all $t$. Let $M$ be a finitely generated $R$-module. Assume, either $R$ is a domain, $M=R$ and $I, J_{\bullet}$ satisfies \textbf{Condition C} or $R$ is not necessarily a domain, the family $J_n= \bp{J}{n}$ for an ideal $J$ such that $I+J$ is $\mathfrak{m}$-primary. Then for any $s$, $$\lim_{t \to \infty}h_{M,I+I_t,J_{\bullet}+\bpq{J'_t}{q}}(s)=h_{M,I,J_{\bullet}}(s),$$ 
    
\noindent where for a fixed $t$, $J_{\bullet}+\bpq{J'_t}{q}$ denotes the family of ideals $n \mapsto J_n+ \bpq{J'_t}{p^n}$. Moreover, on any compact set in $(0,\infty)$, this convergence is uniform with respect to $s$.

\end{theorem}

\begin{proof}
If $s \leq 0$ then both sides are 0, so there is nothing to prove. Fix $0<s_1<s_2<\infty$, it suffices to prove the uniform convergences on $[s_1,s_2]$.
By \Cref{th: Lipschitz continuity for a family}, (1), there exists a positive integer $t_0$ such that for all $t \geq t_0$,
%.............................
$$h_{n,M,I+I_t,J_\bullet}(s)=h_{n,M,I+I_t,J_{\bullet}+\bp{J'_t}{n}}(s),\,\, \forall s \in [s_1, s_2]\, \text{and}\, \forall n \in \mathbb{N}.$$

By \Cref{th: Lipschitz continuity for a family}, (3), there exist a positive integer $t_1$ and a positive real number $c$, such that for all $t \geq t_1$, $n \in \mathbb{N}$ and $s \in [s_1,s_2]$,
$h_{n,M,I+I_t,J_\bullet}(s) \geq h_{n,M,I,J_\bullet}(s-c/t)$. The last two comparisons imply, for any $t \geq t_0+ t_1$ and $s \in [s_1,s_2]$, 
$$|h_{n,M,I,J_{\bullet}}(s)-h_{n,M,I+I_t,J_{\bullet}+\bp{J'_t}{n}}(s)| \leq h_{n,M,I,J_{\bullet}}(s)- h_{n,M,I,J_\bullet}(s-c/t).$$
Taking limit as $n$ approaches infinity, the last comparison gives
$$|h_{M,I,J_{\bullet}}(s)-h_{M,I+I_t,J_{\bullet}+\bp{J'_t}{n}}(s)| \leq h_{M,I,J_{\bullet}}(s)- h_{M,I,J_\bullet}(s-c/t).$$
Here we used our hypothesis to conclude the convergences of all the sequences involved. The deemed uniform convergence now follows from the local Lipschitz continuity of $h_{M,I,J_{\bullet}}$ on $\mathbb{R}_{>0}$; see \Cref{th: Lipschitz continuity for a family}, \Cref{pr: cotinuity of $h$ functions}.
   
\end{proof}

\noindent We record the associativity formula and the additivity property for $h$-function.

\begin{proposition} (Associativity formula)\label{pr: associativity formula for h function}
Let $M$ be a $d$-dimensional finitely generated $R$-module. Let $P_1, P_2, \ldots,$ $P_t$ be the $d$-dimensional minimal primes in the support of $M$. Then,
\[h_{M,I,J,d}(s)= \sum \limits_{j=1}^{t} l_{R_{P_j}}(M_{P_j})h_{R/P_j, IR/P_j, JR/P_j, d}(s).\]
\end{proposition}

\begin{proof}
 By replacing $R$ by $R/\textup{ann}(M)$, we can assume $\dim(R)=d$. We can always assume $R$ is reduced. Indeed, since $R$ is noetherian, we can choose $e_0$ such that the image of the nilradical of $R$ under the $e_0$-th iteration of the Frobenius is zero. Now by \Cref{c: replacing by frobenius pushforward}, we can replace $M$ by some $F^{e_0}_*M$ and pass to the reduced case. Once $R$ is reduced, the two modules

 \[M \, \, \textup{and} \,\, \bigoplus \limits_{j=1}^{t}(\frac{R}{P_j})^{l_{R_{P_j}}(M_{P_j})}\]

 are isomorphic after localizing at each of the primes $P_j$'s. So the result follows from \Cref{le: h is determinde by generic behaviour}.
\end{proof}

\begin{theorem}(Additivity property)\label{th: additivity property of h function}
Let $0 \rightarrow M' \rightarrow M \rightarrow M'' \rightarrow 0$ be a short exact sequence of finitely generated $R$-modules, $I,J$ be two ideals such that $I+J$ is $\mathfrak{m}$-primary. Then $h_{M,I,J, \dim(M)}= h_{M',I,J, \dim(M)}+  h_{M'',I,J, \dim(M)}$.   
\end{theorem}

\begin{proof}
    The proof follows by using \Cref{pr: associativity formula for h function} and \Cref{th: uniform convergence of h function for modules}.
\end{proof}

\begin{remark}\label{re: associativity and additivity of the other functions}
    Given ideals $I,J$ of a local ring $(R, \mathfrak{m})$ and a finitely generated $R$-module, we associate a Frobenius-Poincar\'e function in \Cref{sse: Frob-Poincare in the local setting} and a Hilbert-Kunz density function in \Cref{se: differentiability of h function}. Using \Cref{re: the sequence defining Frobenius Poincare when J is big}, one can formulate and prove an analogous associativity formula and additivity property of the corresponding Frobenius-Poincar\'e function. At the points where the relevant $h$-functions are differentiable, using \Cref{th: differentiability of h function implies existence of density function}, (2), one can formulate and prove a pointwise associativity formula and an additivity property of the density function. In the graded case, these associativity formulas and the additivity properties coincide with \cite[Thm 4.3, Prop 4.4]{AlapanExist}, \cite[Prop 2.14]{TriExist}, see \Cref{se: relation among different function}. 
\end{remark}

Next, we describe the $h$-function of one dimensional rings. This description and \Cref{th: h function and integration} can be combined to construct other explicit examples of $h$-functions. 

\begin{proposition}\label{pr: h-function of one dimensional ring}
Let $(R,\mathfrak{m})$ be a one dimensional ring. Let $I,J$ be proper ideals such that $I+J$ is nonzero and $M$ be a finitely generated $R$-module. 

\begin{enumerate}
    \item For $s \leq 0$, $h_{M,I,J}(s)=0$. If $I=0$, then $h_{M,I,J}(s)= e_{HK}(J,M)$, for positive $s$. If $J=0$, then $h_{M,I,J}(s)= e(I,M)s$ for nonnegative $s$, where $ e(I,M)$ denotes the Hilbert-Samuel multiplicity. 

    \item Assume $R$ is a discrete valuation ring with a uniformizer $\pi$, $I= (\pi^a)$, $J= (\pi^b)$ for some positive integers $a,b$. Then
    $$h_{M,I,J}(s)= \begin{cases}
        \rank_R{(M)}as, \,\,\,\, \text{if}\,\, s \leq b/a,\\
        \rank_R{(M)}b, \,\,\,\,\text{if}\,\, s \geq b/a ,
    \end{cases}$$
    Where $\rank_R(M)$ is the generic rank of $M$.
    
    \item Assume both $I,J$ are nonzero. Then $h_{M,R,I,J}$ is a positive rational linear combination of piecewise linear functions whose slopes are integers. 
\end{enumerate}
\end{proposition}

\begin{proof}
    The claims in (1) follow from a direct calculation. For (2) and (3), we can assume $R$ is a domain and $M=R$ by \Cref{pr: associativity formula for h function}. Now, (2) follows from a direct calculation once we note for $s>0$,
    $l(R/(\pi^{a \ceil{sq}}R+ \pi^{bq}R)= \text{min}\{a \ceil{sq}, bq\}.$ To deduce (3), denote the normalization of the domain $R$ by $S$ and the maximal ideals of $S$ by $\mathfrak{m}_1, \ldots, \mathfrak{m}_r$. Since each $S_{\mathfrak{m}_{j}}$ is a discrete valuation ring, 
 (3) follows from our claim:
    $$h_{R,I,J}= \sum \limits_{j=1}^{r}h_{S_{\mathfrak{m}_{j}}, IS_{\mathfrak{m}_{j}}, JS_{\mathfrak{m}_{j}}} \,\,\,\text{for}\,\, s>0.$$
    To that end, by \Cref{pr: associativity formula for h function}, $h_{R,I,J}= h_{S,I,J}/\rank_R{S}$, where $S$ is considered as an $R$-module.
    For any $s>0$, and $q$, since 
    $$l(\frac{S}{I^{\ceil{sq}}S+\bpq{J}{q}S})= \sum \limits_{j=1}^{r}l(\frac{S_{\mathfrak{m}_{j}}}{I^{\ceil{sq}}S_{\mathfrak{m}_{j}}+\bpq{J}{q}S_{\mathfrak{m}_{j}}}),
    $$
    $h_{S,I,J}= \sum \limits_{j=1}^{r}h_{S_{\mathfrak{m}_{j}}, IS_{\mathfrak{m}_{j}}, JS_{\mathfrak{m}_{j}}}$. Thus the claim follows.
\end{proof}

We analyze how $h_{R,I,J}$ depends on different closure operations of ideals. The next result is used in \Cref{pr: monotonicity result for the desnity function} and implicitly in \Cref{se: applications}.  We refer to \cite{HunekeSwanson} for results on integral closure of ideals. 

\begin{proposition}\label{pr: invariance under closure operations}
    Let $(R, \mathfrak m)$ be a noetherian local ring of dimension $d$. \begin{enumerate}
        \item Let $I, J$ be ideals such that $I+J$ is $\mathfrak m$-primary. Let $J^*$ be the tight closure of $J$. Then $h_{R,I,J}= h_{R, I, J^*}$.
        
        \item Assume $R$ is a domain. Let $I, J_\bullet$ satisfy \textbf{Condition C}\footnote{We do not need the domain assumption when $J_n= \bpq{J}{q}$.} . Let $\overline I$ be the integral closure of $I$. Then $h_{R,I,J_\bullet}= h_{R, \overline I, J_\bullet}$.
    \end{enumerate} 
\end{proposition}

\begin{proof}
   We first prove (1).   \begin{equation}\label{eq: invariance under closure}
   \frac{1}{q^d}h_{n, R,I, J}(s)-  \frac{1}{q^d}h_{n, R,I, J^*}(s)= \frac{1}{q^d}l(\frac{I^{\ceil{sq}}+ \bpq{(J^{*})}{q}}{I^{\ceil{sq}}+ \bpq{(J)}{q}})\leq \frac{1}{q^d}l(\frac{\bpq{(J^{*})}{q}}{\bpq{J}{q}})\ .
   \end{equation}
    Since $R$ is noetherian there is an element $c \in R$ which is not any minimal primes of $R$ such that $c\bpq{(J^{*})}{q} \subseteq \bpq{J}{q}$ for all $q$. Fix a choice of $r$-many generators of $J^*$, the $q$-th powers of these generate $\bpq{(J^*)}{q}$. Thus the length of $\bpq{(J^{*})}{q}/\bpq{J}{q}$ is bounded above by $rl(R/(c, \bpq{J}{q}))$. Since the Krull dimension of $\frac{R}{cR}$ is at most $d-1$, $l(R/(c, \bpq{J}{q}))= O(q^{d-1})$ by \Cref{le: bound on the h function for a family}. Thus taking limit as $q$ approaches infinity in \Cref{eq: invariance under closure}, we conclude $h_{R,I, J}(s)= h_{R,I, J^*}(s)$.

    For (2), recall that there is a natural number $n_0$ such that
    $$I^{n+n_0} \subseteq \overline I ^{n+n_0} \subseteq I^n, $$ for all natural numbers $n$.
    Thus for a positive real number $s$ and $q$ large enough,
    $$\frac{1}{q^d}h_{n, R,I, J_\bullet}(s-\frac{n_0}{q}) \leq \frac{1}{q^d}h_{n, R,\overline I, J_\bullet}(s) \leq \frac{1}{q^d}h_{n, R,I, J_\bullet}(s).$$
    Thus for a positive real number $s$, $h_{R,\overline I, J_\bullet}(s)= h_{R, I, J_\bullet}(s)$; see \Cref{th: Lipschitz continuity for a family}. The desired equality at zero follows from definition.
\end{proof}

\section{Frobenius-Poincar\'e function in the local setting}\label{sse: Frob-Poincare in the local setting}
Given two ideals $I, J$ of $(R, \mathfrak{m})$, such that $I+J$ is $\mathfrak{m}$-primary and a finitely generated $R$-module $M$ we associate a local version of Frobenius-Poincar\'e function $F_{M,I,J}$ and relate it to the $h$-function; see \Cref{th: exitence of Frob-Poincare for general ideals}. The role of underlying graded structures in the classical setting of \cite{AlapanExist} is replaced by the $I$-adic filtration in the local setting. The local version of Frobenius-Poincar\'e function recovers the classical graded one; see \Cref{pr: the graded Frob-Poincare coincides wth the local Frob Poincare}. Unlike \cite{AlapanExist}, the local version $F_{M,I,J}$ is defined even when $J$ is not $\mathfrak{m}$-primary, at the cost of being holomorphic just on the open lower half plane instead of being an entire function.

Given an ideal $I$ and a family $J_{\bullet}$ and a finitely generated $R$-module $M$, set 
\[f_{n,M,I, J_{\bullet}}(s)= h_{{n,M,I, J_{\bullet}}}(s+\frac{1}{q})- h_{{n,M,I, J_{\bullet}}}(s).\]
    When $J_n= \bpq{J}{q}$, $f_{n,M,I,J}(s)$ denotes $f_{n,M,I, J_{\bullet}}(s)$. We drop one or more parameters in $f_{n,M,I, J_{\bullet}}$ when there is no resulting confusion. For the rest of this article, we denote the imaginary part a complex number $y$ by $\Im y$ and the open lower half complex plane by $\Omega$, i.e. $\Omega= \{ y \in \mathbb{C} \, |\, \Im y<0\}.$

\begin{lemma}\label{le: absolute convergence of the terms of sequence}
Let $(R,\mathfrak{m},k)$ be a local ring of dimension $d$, $I$, $J$ be two $R$-ideal, $I+J$ is $\mathfrak{m}$-primary, and $M$ be a finitely generated $R$-module. Consider the function defined by the infinite series
\[F_{n,M,I,J}(y):=\sum_{j=0}^{\infty} f_{n,M,I,J}(j/q)e^{-iyj/q}\] 
Then $F_{n,M,I,J}(y)$ defines a holomorphic function on $\Omega$. We often drop one or more parameters in $F_{n,M,I,J}$ when there is no chance of confusion.
\end{lemma}
\begin{proof}
There is a polynomial $P$ such that $f_{n,M}(s) \leq h_{n,M}(s+1) \leq P(s)$ for any $s$; see \Cref{th: bounds on h function special case}, \Cref{pr: cotinuity of $h$ functions}, assertion (2). Thus
$$|f_{n,M,R,I,J}(j/q)e^{-iyj/q}| \leq P(j/q)e^{j\Im y/q}.$$
Since for fixed $\epsilon>0$, the series $\sum_{0 \leq j < \infty}P(j/q)e^{-j\epsilon/q}$ converges, on the region where $\Im y<-\epsilon$, the sequence of functions $\sum_{j=0}^{\infty} f_{n,M,R,I,J}(j/q)e^{-iyj/q}$ converges uniformly. The limit function is thus holomorphic \cite[Theorem 1, Chap 5]{Ahlfors}. Taking union over all $\epsilon>0$, we see $F_{n,M}(y)$ exists and is holomorphic on $\Omega$.
\end{proof}

\begin{remark}\label{re: the sequence defining Frobenius Poincare when J is big}
    For a big $p$, $p^{-1}$-family $J_{\bullet}$, the analogous $F_{n,M,I, J_\bullet}(y)$ defined using $f_{n,M,I, J_{\bullet}}$ is entire since the corresponding sum is a finite sum of entire functions. 
\end{remark}

Now, we want to check the convergence of $(F_{n,M,I,J}(y)/q^{\dim(M)})_n$ whenever it exists. We will be repeatedly using the dominated convergence: If a sequence of measurable functions $f_n$ converges to $f$ pointwise on a measurable set $\Sigma$ and there is a measurable function $g$ such that $|f_n| \leq g$ on $\Sigma$ for any $n$ and $\int_\Sigma |g|<\infty$, then $\int_\Sigma |f_n-f|$ converges to $0$, so in particular $\int_\Sigma f_n$ converges to $\int_\Sigma f$.

\begin{theorem}\label{th: exitence of Frob-Poincare for general ideals}
Let $(R,\mathfrak{m},k)$ be a local ring, $I$, $J$ be two $R$-ideals such that $I+J$ is $\mathfrak{m}$-primary, and $M$ be a finitely generated $R$-module of dimension $d$.
\item (1) Assume $J$ is $\mathfrak{m}$-primary. Then $F_{M,I,J}(y)=\lim_{n \to \infty}F_{n,M}(y)/p^{nd}$ exists for all $y \in \mathbb{C}$. This convergence is uniform on any compact set of $\mathbb{C}$. Suppose $h_M(s)$ is constant for $s \geq C$, then $F_{M,I,J}(y)=\int^C_0 h_M(t)iye^{-iyt}dt+h_M(C)e^{-iyC}$.

\item (2) Assume $J$ is not necessarily $\mathfrak{m}$-primary. Then for every $y \in \Omega$, $F_{n,M}(y)/p^{nd}$ converges to
\[F_{M,I,J}(y)= \int\limits_{0}^{\infty}h_{M}(t)e^{-iyt}iydt.\]
Moreover, this convergence is uniform on any compact subset of $\Omega$ and $F_{M}(y):= F_{M,I,J}(y)$ is holomorphic on $\Omega$.
\end{theorem}

\begin{proof}
\item (1) Since $J$ is $\mathfrak{m}$-primary, then $h_M(s)=h_M(C)$ for some fixed $C>0$ and any $s \geq C$; see \Cref{le: compactly supported h-function} and \Cref{pr: associativity formula for h function}. Then,
\begin{equation*}
\begin{split}
F_{n,M}(y)&=\sum_{j=0}^{\infty} f_{n,M}(j/q)e^{-iyj/q}\\
&=\sum_{j=0}^{\infty} (h_{n,M}((j+1)/q)-h_{n,M}(j/q))e^{-iyj/q}\\
&=\sum_{j=0}^{Cq-1} (h_{n,M}((j+1)/q)-h_{n,M}(j/q))e^{-iyj/q}\\
&=\sum_{j=0}^{Cq-1} h_{n,M}(j/q)(e^{-iy(j-1)/q}-e^{-iy(j)/q})+h_{n,M}(C)e^{-iy(C-\frac{1}{q})}\\
&=\sum_{j=0}^{Cq-1} h_{n,M}(j/q)e^{-iyj/q}(e^{iy/q}-1)+h_{n,M}(C)e^{-iy(C-\frac{1}{q})}\\
&=\int_{0}^{C} h_{n,M}(t)e^{-iy\ceil {tq}/q}q(e^{iy/q}-1)dt+h_{n,M}(C)e^{-iy(C-\frac{1}{q})} \ .
\end{split}
\end{equation*}

\noindent Fix a compact subset $K$ of $\mathbb C$. Given $\delta>0$, choose $b>0$ such that for all $y \in K$, $t \in \mathbb R$ and $n \in \mathbb N$

\[ \int_{0}^{b}(\frac{1}{q^d}|h_{n,M}(t)e^{-iy\ceil {tq}/q}q(e^{iy/q}-1)|+ |h_{M}(t)e^{-iyt}(iy)|)dt \leq \frac{\delta}{2}.\]

We have

\begin{gather*}
  |\frac{1}{q^d}F_{n,M}(y) - \int_{0}^{C} h_{M}(t)e^{-iyt}(iy)dt-h_{M}(C)e^{-iyC}| \\
  \leq \int_{0}^{C} |h_{n,M,d}(t)e^{-iy\ceil {tq}/q}q(e^{iy/q}-1)- h(y)iye^{-iyt}|dt+|h_{n,M}(C)e^{-iy(C-\frac{1}{q})}-h_M(C)e^{-iyC}|\\
  \leq \int_{0}^{b} (|h_{n,M,d}(t)e^{-iy\ceil {tq}/q}q(e^{iy/q}-1)|+ |h_{M}(t)e^{-iyt}(iy)|)dt\\
  + \int_{b}^{C} |h_{n,M,d}(t)e^{-iy\ceil {tq}/q}q(e^{iy/q}-1)- h(y)iye^{-iyt}|dt+ |h_{n,M}(C)e^{-iy(C-\frac{1}{q})}-h_M(C)e^{-iyC}| . 
\end{gather*}

Moreover for $y \in K$, there is a constant $C'$ independent of $n$ such that for all $t \in [b,C]$
$$|h_{n,M,d}(\lceil tq \rceil/q)-h_M(t)| \leq C'/q \, \text{and}\, |e^{-iy\lceil tq \rceil/q}q(e^{iy/q}-1)-e^{iyt}(iy)| \leq C'/q.$$ 

Thus we can choose $N_0$ such that for all $n \geq N_0$ and $y \in K$,

\[|\frac{1}{q^d}F_{n,M}(y) - \int_{0}^{C} h_{M}(t)e^{-iyt}(iy)dt-h_{M}(C)e^{-iyC}| \leq \delta.\]

This proves the desired uniform convergence.

\item (2) We prove uniform convergence of $F_{n,M}/q^{d}$ to the integral on every compact subset of $\Omega$; the holomorphicity of $F_M$ is then a consequence of \cite[Theorem1, Chap 5]{Ahlfors}. We have

\begin{equation*}
\begin{split}
F_{n,M}(y)&=\sum_{j=0}^{\infty} f_{n,M}(j/q)e^{-iyj/q}\\
&=\sum_{j=0}^{\infty} (h_{n,M}((j+1)/q)-h_{n,M}(j/q))e^{-iyj/q}\\
&=\sum_{j=0}^{\infty} h_{n,M}(j/q)(e^{-iy(j-1)/q}-e^{-iy(j)/q})\\
&=\sum_{j=0}^{\infty} h_{n,M}(j/q)e^{-iyj/q}(e^{iy/q}-1)\\
&=\int_{0}^{\infty} h_{n,M}(t)e^{-iy\ceil{tq}/q}q(e^{iy/q}-1)dt \ .
\end{split}
\end{equation*}

\noindent The rearrangements leading to the second and third equality are possible thanks to the absolute convergences implied by \Cref{th: bounds on h function special case}. 
Fix any compact $K \subseteq \Omega$. Using triangle inequality, we get

\begin{gather*}
   |h_{n,d}(t)e^{-iy\frac{\ceil {tq}}{q}}q(e^{iy/q}-1)- h(t)e^{-iyt}(iy)|\\
   \leq |h_{n,d}(t)-h(t)||e^{-iy\frac{\ceil{tq}}{q}}q(e^{iy/q}-1)|+ |h(t)||e^{-iy\frac{\ceil{tq}}{q}}-e^{-iyt}||q(e^{iy/q}-1)|\\+ |h(t)||e^{-iyt}||q(e^{iy/q}-1)-iy|\\
   = |h_{n,d}(t)-h(t)||e^{-iy\frac{\ceil{tq}}{q}}q(e^{iy/q}-1)|+ |h(t)e^{-iyt}||e^{-iy(\frac{\ceil{tq}}{q}-t)}-1||q(e^{iy/q}-1)|\\
   + |h(t)e^{-iyt}||q(e^{iy/q}-1)-iy| \ .
\end{gather*}

It follows from the power series expansion of $e^z$ at zero and the boundedness of $K$ that there are constants $C_1$, $C_2$ such that for all $y \in K$ , $t \in \mathbb R$ and $n \in \mathbb N$

\[|q(e^{iy/q}-1)| \leq C_1|y|, \, |q(e^{iy/q}-1)-iy| \leq C_2\frac{|y|^2}{q},\, |e^{-iy(\frac{\ceil{tq}}{q}-t)}-1| \leq C_1|y(\frac{\ceil{tq}}{q}-t)|.\]

Choose $\epsilon >0$ such that $K\subseteq \{y \in \mathbb C \, | \, \Im y< -\epsilon\}$. Using the comparisons above, we get for all $y \in K$ , $t \in \mathbb R$ and $n \in \mathbb N$,

\begin{gather*}
    |h_{n,d}(t)e^{-iy\frac{\ceil {tq}}{q}}q(e^{iy/q}-1)- h(t)e^{-iyt}(iy)|\\
    \leq |h_{n,d}(t)-h(t)|e^{-\epsilon t}C_1|y|+ |h(t)e^{-\epsilon t}|C_1^2|y|^2|\frac{\ceil{tq}}{q}-t|+ |h(t)e^{\epsilon t}|C_2\frac{|y|^2}{q}\\
    \leq |h_{n,d}(t)-h(t)|e^{-\epsilon t}C_1|y|+ |h(t)e^{-\epsilon t}|C_1^2\frac{|y|^2}{q}+ |h(t)e^{-\epsilon t}|C_2\frac{|y|^2}{q} \ .
\end{gather*}

Taking integral on $\mathbb R_ {\geq 0}$, we get for $y \in K$ and all $n \in \mathbb N$
\begin{gather*}
|\frac{1}{q^d}F_{n,M}(y)-F_{M,I,J}(y)|\\
\leq C_1|y|\int \limits_{0}^{\infty}|h_{n,d}(t)-h(t)|e^{-\epsilon t}dt+ (C_1^2+C_2)\frac{|y|^2}{q}\int \limits_{0}^{\infty}|h(t)|e^{-\epsilon t}dt \ .
\end{gather*}
\noindent Thanks to \Cref{th: bounds on h function special case}, (2), we can choose a polynomial $P_2 \in \mathbb R [t]$ such that  $|h_{n,d}(t)|\leq |P_2(t)|$ for all $n$ and $t \in \mathbb R$. Since $|P_2(t)e^{-\epsilon t}|$ is integrable on $\mathbb R_{\geq 0}$, by dominated convergence

\[\lim_{n \to \infty}\int \limits_{0}^{\infty}|h_{n,d}(t)-h(t)|e^{-\epsilon t}dt=0.\]

\noindent Using this in the last inequality implies uniform convergence of $\frac{1}{q^d}F_{n,M}(y)$ to $F_{M,I,J}(y)$ on $K$.
\end{proof}

\begin{remark}
Suppose $h_M(y)$ is constant for $y \geq C$. Since for $y \in \Omega$, $h_M(C)e^{-iyC}$ converges to zero as $C$ approaches infinity, the two descriptions of $h_M$ in this case match on $\Omega$. When $J_{\bullet}$ is both big $p$ and $p^{-1}$, our argument actually produces a corresponding entire function $F_{M,I,J_{\bullet}}(y)$.  
\end{remark}

\begin{definition}\label{de: Frobenius-Poincare general version}
 Let $I,J$ be two ideals in $(R, \mathfrak m)$ such that $I+J$ is $\mathfrak m$-primary. For a finitely generated $R$-module $M$, the function $F_{M,I,J}(y)$ is called the \textit{Frobenius-Poincar\'e function} of $(M,I,J)$.

  We drop one or more parameters from $F_{M,I,J}$ when there is no possible source of confusion.
\end{definition}

The next result directly follows from \Cref{pr: associativity formula for h function}.

\begin{corollary}Let $M,N$ be two $R$-modules such that $\dim M=\dim N=\dim R$ and their localizations are isomorphic at all $P \in \Assh R$. Then $F_M(y)=F_N(y)$.
\end{corollary}
\begin{proof}
This is true because $h_M(s)=h_N(s)$.
\end{proof}

Next we prove that the map sending the pair $I,J$ to $F_{M,I,J}$ is $\mathfrak{m}$-adically continuous in both $I$ and $J$:

\begin{proposition} \label{pr: m adic continuity of the Frobenius Poincare function}Let $t \in \mathbb{N}$, $I_t$, $J_t$ be two sequences of ideals such that $I_t+J_t \subset\mathfrak{m}^t$. Then for any $y \in \Omega$: the open lower half complex plane, $\lim_{t \to \infty}F_{M,I+I_t,J+J_t}(y)=F_{M,I,J}(y)$. If $J$ is $\mathfrak{m}$-primary, then the above holds for $y \in \mathbb{C}$. In either case, the convergence is uniform on a compact subset of $\Omega$ or $\mathbb{C}$.
\end{proposition}

\begin{proof}
Fix a compact subset $K$ of $\Omega$. Choose $\epsilon >0$ such that $\Im y<- \epsilon$ for all $y \in K$. Recall from \Cref{th: bounds on h function special case}, that there is a polynomial $P \in \mathbb{R}[t]$ such that $h_{n,M,I,J}(s) \leq P(s)$ for all $s \in \mathbb R$ and all $n$; so $h_{M,I+I_t,J+J_t}(s)\leq P(s)$ for all $s$. Notice $|P(s)e^{-\epsilon s}|$ is integrable on $\mathbb R_{\geq 0}$ and the sequence $h_{M,I+I_t,J+J_t}$ converges to $h_{M,I,J}$; the convergence is uniform on every compact subset of $(0, \infty)$; see \Cref{th: m-adic continuity of h function}. Say the absolute values of elements of $K$ is bounded above by $D$. Given $\delta>0$, the observations above allows us to choose an interval $[a,b] \subseteq (0, \infty)$ and $t_0 \in \mathbb N$ such that,\\
\item[(a)] $2\int_{0}^{a} |P(s)|e^{-\epsilon s}ds+ 2\int_{b}^{\infty} |P(s)|e^{-\epsilon s}ds \leq \frac{\delta}{2D}.$\\
\item[(b)] $|h_{M,I+I_t, J+J_t}(x)- h_{M,I,J}(x)| \leq \frac{\delta}{2D\int_{a}^{b}e^{-\epsilon s}ds} \, \text{for all}\, t\geq t_0\, \text{and all}\, s \in [a,b].$

\noindent
Therefore by using \Cref{th: exitence of Frob-Poincare for general ideals}, for $y \in K$ and all $t \geq t_0$
\begin{equation*}
\begin{split}
    |F_{M,I+I_t, J+J_t}(y)- F_{M,I,J}(y)|&\leq \int_{0}^{\infty}|y||h_{M,I+I_t,J+J_t}(s)-h_{M,I,J}(s)|e^{-\epsilon s}ds\\
    & \leq D [2\int_{0}^{a} |P(s)|e^{-\epsilon s}ds+ 2\int_{b}^{\infty} |P(s)|e^{-\epsilon s}ds\\&+ \int_{a}^{b}|h_{M,I+I_t, J+J_t}(s)- h_{M,I,J}(s)|e^{-\epsilon s}ds ]  \\
    & \leq \delta \ .
\end{split}
\end{equation*}
This proves uniform convergence of $(F_{M,I+I_t, J+J_t}(y))_t$ to $F_{M,I,J}(y)$ on every compact subset of $\Omega$. The assertion for $\mathfrak m$-primary $J$ follows from a similar argument.
\end{proof}

\section{Existence and properties of density functions in the local setting}\label{se: differentiability of h function}
In this section, we associate a density function in the local setting to a pair $(I, J_{\bullet})$ satisfying \textbf{Condition C}. This extends the classical Hilbert-Kunz density function in the graded setting (\cite{TriExist}) to the local setting when $J_n= \bp{J}{n}$ for an $\mathfrak{m}$-primary ideal $J$; see \Cref{th: h function differentiable when J is m primary}. We prove various properties of the resulting density function and the $h$-function such as the monotonicity property: \Cref{pr: monotonicity result for the desnity function}, and \Cref{pr: one sided limits exists}, which are used in the later sections.

\begin{definition}\label{de: HIlbert-Kunz density}
Let $I$ be an ideal and $J_{\bullet}$ be a family of ideals in $(R,\mathfrak m)$ satisfying \textbf{Condition C}. For a finitely generated $R$-module $M$ and $s \in \mathbb R$, recall

\[f_{n,M,I,J_{\bullet}}(s)= h_{n,M,I,J_\bullet}(s+\frac{1}{q})- h_{n,M,I,J_{\bullet}}(s)= l(\frac{(I^{\ceil{sq}}+J_n)M}{(I^{\ceil{sq}+1}+J_n)M}). \]

\noindent Whenever $((\frac{1}{p^n})^{\dim(M)-1}f_{n,M,I,J_{\bullet}}(s))_n$ converges, we call the limit the \textit{density function} of $(M,I, J_{\bullet})$ at $s$ and denote the limit by $f_{M,I,J_{\bullet}}(s)$. Whenever $f_{M,I,J_{\bullet}}(s)$ exists for all $s \in \mathbb R$, the resulting function $f_{M,I,J_{\bullet}}$ is called the \textit{density function} of $(M,I,J_{\bullet})$. 

We often drop one or more parameters from $f_{n,M,I,J_{\bullet}}(s), f_{M,I,J_{\bullet}}(s), f_{M,I,J_{\bullet}}$ whenever those are clear from the context.
\end{definition}

\Cref{th: differentiability of h function implies existence of density function} shows that, in general, the $f_{M,I,J_{\bullet}}(s)$ exists outside a countable subset of $\mathbb{R}$. In the sequel, we deal with the situation where $f_{M,I,J_{\bullet}}(s)$ is possibly not defined at all points of $\mathbb R$. So we clarify what we mean by the continuity of the density function $f_{M,I,J_{\bullet}}$ 
at a point.

\begin{definition}\label{de: continuity of density function}
Let $\Gamma \subseteq \mathbb{R}$ be the set of all points $s$ such that $((\frac{1}{p^n})^{\dim(M)-1}f_{n,M,I,J_{\bullet}}(s))_n$ converges, endowed with the subspace topology. For $s_0 \in \Gamma$,  we say that $f_{M,I,J_{\bullet}}$ is continuous at $s_0$, if the function $f_{M,I,J_{\bullet}}: \Gamma \rightarrow \mathbb{R}$ is continuous at $s_0$.
\end{definition}

In \Cref{th: differentiability of h function implies existence of density function}, we relate the existence of $f_{M,I,J_{\bullet}}(s)$ to the differentiability of $h_{M,I,J_\bullet}$ at $s$-whenever the derivative of $h_{M,I,J_\bullet}$ at $s$ exists. We show that $h_{M,I,J_\bullet}$ is always left and right differentiable everywhere on the real line. The key to obtaining these results is \Cref{th: uniform convergence to the convex functional}, which proves convexity of a function constructed from the $h$-function. Note that, the $h$-function being Lipschitz continuous is automatically differentiable outside a set of measure zero. But our method shows that the $h$-function is continuously differentiable outside a countable set. Recall: 

\begin{definition}\label{de: convex function}
Let $S$ be a subset of $\mathbb R$. We call a function $\lambda: S \rightarrow \mathbb R$ to be \textit{convex} if for elements of $S$, $s_1< s_2\leq t_1< t_2$,

\[\frac{\lambda(s_2)- \lambda(s_1)}{s_2 -s_1} \geq \frac{\lambda(t_2)- \lambda(t_1)}{t_2 -t_1}.\]
\end{definition}

Convexity is a notion that appears naturally in mathematical analysis. For references on convex functions, see \cite{Convex2006}.

\subsection{Convexity property of $h$-functions} \label{sse: construction of convex function}Let $I,J_{\bullet}, M$ be as above. Now we lay the groundwork for the construction of the convex function $\mathcal{H}(s,s_0)$ in \Cref{th: uniform convergence to the convex functional}. Fix $\mu$ such that $I$ is generated by $\mu$-many elements. Set $M_q= M/J_nM$ and $S$ to be the polynomial ring in $\mu$ many variables over $R/ \mathfrak m$. Given a compact interval $[a,b] \subseteq (0, \infty)$, thanks to \Cref{th: bound on the density function} we can choose $C$ such that for all $x \in [a,b]$ and $n \in \mathbb N$

\[\frac{I^{\ceil{xq}}M_q}{I^{\ceil{xq}+1}M_q}= h_{n}(x+\frac{1}{q})- h_{n}(x) \leq Cq^{\dim(M)-1}.\]
Recall,
\[l(S_{\ceil {xq}})={{\mu+\ceil{xq}-1}\choose{\mu-1}}=1/(\mu-1)!(\ceil{xq})^{\mu-1}+O(\ceil{xq}^{\mu-2}).\]
Fix $s_0 \in \mathbb R$. Taking cues from these two estimates, for $s >s_0$ we define 
\begin{equation}\label{eq: the sequence defining the convex H}
\mathcal{H}_n(s,s_0)=\sum^{\ceil{sq}-1}_{j=\ceil{s_0q}}q^{\mu-\dim(M)-1}l(I^jM_q/I^{j+1}M_q)/l(S_j)\ .  
\end{equation}

\begin{theorem}\label{th: uniform convergence to the convex functional}
Let $I,J_\bullet$ in the local ring $(R,\mathfrak m)$ satisfy \textbf{Condition C}, $M$ be a finitely generated $R$-module of Krull dimension $d$, $I$ be generated by a set of $\mu$ elements. Set $M_q= M/J_nM$, fix $s_0 \in \mathbb R_{>0}$. Assume $M,I,J_{\bullet}$ satisfy either (A) or (B) below:

\begin{enumerate}
    \item[(A)] $R$ is a domain and $M=R$.
    \item[(B)] $J_n= \bpq{J}{q}$ for some ideal $J$ such that $I+J$ is $\mathfrak m$-primary and $M$ is any finitely generated $R$-module.
\end{enumerate}

Set $c(s)= \frac{s^{\mu-1}}{(\mu-1)!}$, where $\mu$ is such that $I$ admits a set of generators consisting of $\mu$ elements. In the context of (A) or (B)\footnote{$h_{M,I,J_{\bullet}}$ exists in the context of (A) or (B)}, set
      \[\mathcal{H}(s,s_0)=h_{M,I, J_{\bullet}}(s)/c(s)-h_{M,I, J_{\bullet}}(s_0)/c(s_0)+\int^{s}_{s_0} h_{M,I, J_{\bullet}}(t)c'(t)/c^2(t)dt.\]
\begin{enumerate}
    \item 
    On any compact subset of $(s_0, \infty)$, $\mathcal{H}_n(s,s_0)$ uniformly converges to $\mathcal{H}(s,s_0)$.

    \item The function $\mathcal{H}(s,s_0)$ is a convex function on $(s_0, \infty)$.
\end{enumerate}
\end{theorem}

\begin{proof}
(1) Let $\mathcal H_n(s,s_0)$ be as in \Cref{eq: the sequence defining the convex H}. We have
\begin{gather*}
\mathcal{H}_n(s,s_0)=\sum^{\ceil{sq}-1}_{j=\ceil{s_0q}}q^{\mu-d-1}l(I^jM_q/I^{j+1}M_q)/l(S_j)\\
=\sum^{\ceil{sq}-1}_{j=\ceil{s_0q}}q^{\mu-d-1}(l(M_q/I^{j+1}M_q)-l(M_q/I^jM_q))/l(S_j)\\
=q^{\mu-d-1}l(M_q/I^{\ceil{sq}}M_q)/l(S_{\ceil{sq}-1})-q^{\mu-d-1}l(M_q/I^{\ceil{s_0q}}M_q)/l(S_{\ceil{s_0q}})\\
+\sum^{\ceil{sq}-1}_{j=\ceil{s_0q}+1}q^{\mu-d-1}l(M_q/I^jM_q)(1/l(S_{j-1})-1/l(S_j)) \ .
\end{gather*}

\noindent Since we are in the context of (A) or (B), $q^{\mu-d-1}l(M_q/I^{\ceil{sq}}M_q)/l(S_{\ceil{sq}-1})$ converges to $h(s)/c(s)$ and $q^{\mu-d-1}l(M_q/I^{\ceil{s_0q}}M_q)/l(S_{\ceil{s_0q}})$ converges to $ h(s_0)/c(s_0)$. Also,

\begin{gather*}
\sum^{\ceil{sq}-1}_{j=\ceil{s_0q}+1}q^{\mu-d-1}l(M_q/I^jM_q)(1/l(S_{j-1})-1/l(S_j))\\
=\int^{s-1/q}_{s_0}\frac{l(M_q/I^{\ceil {tq}}M_q)}{q^d}(\frac{1}{l(S_{\ceil {tq}-1})}-\frac{1}{l(S_{\ceil {tq}})})(q^{\mu})dt \ .
\end{gather*}

\noindent When $q$ approaches infinity, $\frac{l(M_q/I^{\ceil {tq} }M_q)}{q^d}$ converges to $h_M(t)$, and $(\frac{1}{l(S_{\ceil {tq} -1})}-\frac{1}{l(S_{\ceil{tq} })})(q^{\mu})$ converges to $c'(t)/c^2(t)$. Also, all these convergence are uniform on any compact subset of $(0, \infty)$. So we get a uniform convergence (uniform on $s$) on any compact subset of $(s_0, \infty)$:
\begin{align*}
\int^{s-1/q}_{s_0}\frac{l(M_q/I^{\lfloor tq \rfloor}M_q)}{q^d}(\frac{1}{l(S_{\lfloor tq \rfloor-1})}-\frac{1}{l(S_{\lfloor tq \rfloor})})(q^{\mu})dt\\
\to \int^{s}_{s_0} h(t)c'(t)/c^2(t)dt.
\end{align*}
This proves that $\mathcal H_n(s,s_0)$ converges to $\mathcal H(s,s_0)$ and the convergence is uniform on any compact subset of $(s_0, \infty)$.\\

\noindent(2) We claim $\mathcal{H}_n$ is convex on $1/p^n\mathbb{Z} \cap (s_0, \infty)$. To this end, it suffices to show \[\mathcal{H}_n(\frac{i+1}{p^n}, s_0)-\mathcal{H}_n(\frac{i}{p^n}, s_0) \geq \mathcal{H}_n(\frac{i+2}{p^n}, s_0)-\mathcal{H}_n(\frac{i+1}{p^n}, s_0).\] By definition, this is equivalent to showing $$l(I^iM_q/I^{i+1}M_q)/l(S_i) \geq l(I^{i+1}M_q/I^{i+2}M_q)/l(S_{i+1}),$$
which follows from \Cref{l: Boij's theorem strengthened}. This convexity of $\mathcal H_n(s,s_0)$ implies
the convexity of the limit function $\mathcal H(s,s_0)$ on $(s_0, \infty) \cap \mathbb{Z}[1/p]$. Therefore for $s_1<s_2\leq t_1< t_2$ in $(s_0, \infty) \cap \mathbb{Z}[1/p]$,
\[\frac{H(s_2,s_0)-H(s_1,s_0)}{s_2-s_1} \geq \frac{H(t_2,s_0)-H(t_1,s_0)}{t_2-t_1}.\]
Since $\mathcal{H}(s,s_0)$ is continuous on $(s_0, \infty)$,
$(s,t) \to (H(t,s_0)-H(s,s_0))/(t-s)$ is continuous. Moreover as $\mathbb Z[1/p]\cap (s_0, \infty)$ is dense in $(s_0, \infty)$, the slope inequality defining a convex function (see \Cref{de: convex function}) holds for $\mathcal{H}(s,s_0)$ for points in $(s_0, \infty)$. 
\end{proof}

\begin{theorem}\label{th: one sided differentiability of h-functions}
With notations set in the statement of \Cref{th: uniform convergence to the convex functional}, set $\mathcal H(s)= \mathcal H (s,s_0)$. Denote the left and right derivative of a function $\lambda$ at $s \in \mathbb R$ by $\lambda'_{-}(s)$ and $\lambda'_+(s)$ respectively. In the context of situation (A) or (B) stated in \Cref{th: uniform convergence to the convex functional}, 

\begin{enumerate}
\item Outside a countable subset of $(s_0,\infty)$, the derivative of $\mathcal{H}$ exists and is also continuous on the complement of this countable subset. The left and right derivative of $\mathcal{H}$ exists everywhere on $(s_0, \infty)$. The second derivative of $\mathcal{H}$ exists almost everywhere, i.e. outside a subset of Lebesgue measure zero of $(s_0, \infty)$.

\item The left and right derivatives of $\mathcal H$ are both decreasing in terms of $s$. We have $\mathcal{H}'_+(s) \leq \mathcal{H}'_-(s)$, and if $s_1<s_2$, $\mathcal{H}'_-(s_2) \leq \mathcal{H}'_+(s_1)$.

\item Outside a countable subset of $(0,\infty)$, the derivative of $h$ exists and is also continuous. The left and right derivative of $h$ exists everywhere on $(0, \infty)$. The second derivative of $h$ exists almost everywhere on $(0, \infty)$.

\item On $(s_0,\infty)$, $h'_+(s)=\mathcal{H}'_+(s)c(s)$, $h'_-(s)=\mathcal{H}'_-(s)c(s)$ exists, and $h'_+(s) \leq h'_-(s)$ for any $s \in (0, \infty)$. 
\end{enumerate}
\end{theorem}

\begin{proof}
Since a convex function has left and right derivative everywhere, $\mathcal{H}$ has left and right derivative on $(s_0, \infty)$ by \Cref{th: uniform convergence to the convex functional}. Thanks to the convexity of $\mathcal{H}$ proven in \Cref{th: uniform convergence to the convex functional}, (2), outside a countable subset $\Lambda$ of $(s_0, \infty)$, $\mathcal H$ is differentiable. On $(s_0, \infty)\setminus \Lambda$, the derivative of $\mathcal H$ is decreasing as $\mathcal H$ is convex. Now a decreasing function defined on a subset of $\mathbb R$ can have only countably many points of discontinuity. So there is a countable subset of $(s_0, \infty)$ outside which $h$ is differentiable and the derivative is continuous as function on the complement of this countable set.  

(2) follows from properties of convex functions.\\

\noindent (3), (4): Recall
\[\mathcal{H}(s,s_0)=h_{M,I, J_{\bullet}}(s)/c(s)-h_{M,I, J_{\bullet}}(s_0)/c(s_0)+\int^{s}_{s_0} h_{M,I, J_{\bullet}}(t)c'(t)/c^2(t)dt.\]
Since in the context of (A) and (B) $h_{M,I,J_{\bullet}}$ is continuous on $(0, \infty)$, the part of $\mathcal H(s,s_0)$ given by the integral is always differentiable. So (3) follows from the analogous properties of $\mathcal H(s,s_0)$ in (1) by varying $s_0$.  The formulas in (4) follow from a direct computation. That $h'_+(s) \leq h'_{-}(s)$ follows from these formulas and (2).
\end{proof}

\begin{remark}\label{re: density function is differentiable outside a set of measure zero}
Trivedi asks when the Hilbert-Kunz density function of a graded pair $(R,J)$ is $\dim(R)-2$ times continuously differentiable; see \cite[Question 1]{TrivediQuadric}. In general the Hilbert-Kunz density function need not be $\dim(R)-2$ times continuously differentiable; see \cite[Example 8.3.2]{AlapanThesis}. Our work shows that the Hilbert-Kunz density function is always differentiable outside a set of measure zero. Indeed, a convex function on an interval is twice differentiable outside a set of measure zero; see \cite[Section 1.4]{Convex2006}. Thus from \Cref{th: uniform convergence to the convex functional}, it follows that outside a set of measure zero the $h$ function is twice differentiable. Now from \Cref{th: HK density when J is homogeneous but not of finite colength}, we conclude that the Hilbert-Kunz density function of a standard graded domain of dimension at least two is differentiable outside a set of measure zero. 
\end{remark}

\begin{remark}\label{re: reason for assuming condition A or B}
The conclusions of \Cref{th: uniform convergence to the convex functional} and \Cref{th: one sided differentiability of h-functions} are deduced in the context of situation (A) or (B), because we prove existence and continuity of $h_{M,I,J_{\bullet}}$ in those two contexts. So even outside the context of (A) or (B) whenever there is an $h$-function continuous on $(0, \infty)$, we have a corresponding version of \Cref{th: uniform convergence to the convex functional} and \Cref{th: one sided differentiability of h-functions}.    
\end{remark}

\subsection{Differentiability of $h$-function and existence of density function}

\noindent We return to the question of existence of $f_{M,I,J_{\bullet}}(s)$ at a given $s \in \mathbb R$. We make comparisons between the limsup and and liminf of the sequence defining $f_{M,I,J_{\bullet}}(s)$ and the corresponding $h'_+(s)$ and $h'_{-}(s)$.

\begin{lemma}\label{le: alternate limit expression of the one sided derivatives}
With the notation set in \Cref{th: uniform convergence to the convex functional}, set 
$$D_{n,t}= f_{n,M,I,J_{\bullet}}(t/p^n)=h_{n,M,I,J_{\bullet}}((t+1)/p^n)-h_{n,M,I,J_{\bullet}}(t/p^n).$$ In the context of situation (A) or (B),
\begin{enumerate}
\item $$h'_+(s)=\lim_{m \to \infty}\lim_{n \to \infty}\frac{\sum\limits^{\ceil{sp^mp^n}+p^n-1}_{t=\ceil{sp^mp^n}}D_{m+n,t}}{p^{m(d-1)}p^{nd}}.$$
\item $$h'_-(s)=\lim_{m \to \infty}\lim_{n \to \infty}\frac{\sum\limits^{\ceil{sp^mp^n}-1}_{t=\ceil{sp^mp^n}-p^n}D_{m+n,t}}{p^{m(d-1)}p^{nd}}.$$
\end{enumerate}
\end{lemma}
\begin{proof}
(1) Note
\begin{align*}
\sum\limits^{\ceil{sp^mp^n}+p^n-1}_{t=\ceil{sp^mp^n}}D_{m+n,t}\\
=\sum\limits^{\ceil{sp^mp^n}+p^n-1}_{t=\ceil{sp^mp^n}}f_{m+n,M}(t/p^mp^n)\\
=h_{m+n,M}(\ceil{sp^mp^n}/p^mp^n+1/p^m)-h_{m+n,M}(\ceil{sp^mp^n}/p^mp^n) \ .
\end{align*}
Since in the context of (A) or (B), the $h$-function exists, the right hand side of the desired equation in (1) is
\begin{align*}
\lim_{m \to \infty}\lim_{n \to \infty}\frac{h_{m+n,M}(\ceil{sp^mp^n}/p^mp^n+1/p^m)-h_{m+n,M}(\ceil{sp^mp^n}/p^mp^n)}{p^{m(d-1)}p^{nd}}\\
=\lim_{m \to \infty}\frac{h_M(s+1/p^m)-h_M(s)}{1/p^m}\\
=h'_+(s) \ .
\end{align*}
(2) Note
\begin{align*}
\sum\limits^{\ceil{sp^mp^n}-1}_{t=\ceil{sp^mp^n}-p^n}D_{m+n,t}\\
=\sum\limits^{\ceil{sp^mp^n}-1}_{t=\ceil{sp^mp^n}-p^n}f_{m+n,M}(t/p^mp^n)\\
=h_{m+n,M}(\ceil{sp^mp^n}/p^mp^n)-h_{m+n,M}(\ceil{sp^mp^n}/p^mp^n-1/p^m)
\end{align*}
Thus the right hand side of the desired equation in (1) is
\begin{align*}
\lim_{m \to \infty}\lim_{n \to \infty}\frac{h_{m+n,M}(\ceil{sp^mp^n}/p^mp^n)-h_{m+n,M}(\ceil{sp^mp^n}/p^mp^n-1/p^m)}{p^{m(d-1)}p^n}\\
=\lim_{m \to \infty}\frac{h_M(s)-h_M(s-1/p^m)}{1/p^m}\\
=h'_-(s) \ .
\end{align*}
\end{proof}

\begin{theorem}\label{th: differentiability of h function implies existence of density function}With the same notation as in \Cref{th: uniform convergence to the convex functional}, in the context of situation (A) or (B),
\begin{enumerate}
\item for any $s>0$,
$$h'_+(s) \leq \underline{\lim}_{n \to \infty}f_{n,M,I, J_\bullet}(s)/p^{n(d-1)} \leq \overline{\lim}_{n \to \infty}f_{n,M,I, J_\bullet}(s)/p^{n(d-1)} \leq h'_-(s),$$

\noindent where $\underline{\lim}$ and $\overline{\lim}$ denote liminf and limsup respectively.

\item At $s>0$, if $h_M$ is differentiable, then $f_{M,I, J_\bullet}(s)$- the density function of $(M,I,J_{\bullet})$ at $s$ exists and is equal to $h'_{M,I,J_{\bullet}}(s)$. If $h_M(s)$ is a $C^1$-function, then $f_M(s)$ is continuous. 

\item  There is a countable subset of $(0, \infty)$ outside which $f_{M,I,J_{\bullet}}(s)$ exists and is equal to $h'_{M,I,J_{\bullet}}(s)$. 
\end{enumerate}
\end{theorem}

\begin{proof}
(1) In the proof, we also use the notation set in \Cref{le: alternate limit expression of the one sided derivatives}, (1). Set 
\[\alpha_{\mu,t}= {{\mu+t-1}\choose{\mu-1}}.\]
Note $D_{n,t}= l((I^t+J_n)M/(I^{t+1}+J_n)M)$. For a fixed $n$, $D_{n,t}/\alpha_{\mu,t}$ is a decreasing function of $t$, thanks to \Cref{l: Boij's theorem strengthened}. So for $\ceil{sp^mp^n} \leq t \leq \ceil{sp^mp^n}+p^n-1$, $D_{m+n,t}/\alpha_{\mu,t} \leq D_{m+n,\ceil{sp^mp^n}}/\alpha_{\mu,\ceil{sp^mp^n}}$, so
\begin{equation*}
\begin{split}
D_{m+n,t}& \leq D_{m+n,\ceil{sp^mp^n}}\frac{\alpha_{\mu,t}}{\alpha_{\mu,\ceil{sp^mp^n}}}\\
&\leq D_{m+n,\ceil{sp^mp^n}}\frac{\alpha_{\mu,\ceil{sp^mp^n}+p^n}}{\alpha_{\mu,\ceil{sp^mp^n}}} \ .
\end{split}
\end{equation*}

\noindent Also $\alpha_{\mu,t}$ is a polynomial of degree $\mu-1$ in $t$, so
\begin{equation*}
\begin{split}
\lim_{m \to \infty}\lim_{n \to \infty}\frac{\alpha_{\mu,\ceil{sp^mp^n}+p^n}}{\alpha_{\mu,\ceil{sp^mp^n}}}
&=\lim_{m \to \infty}\lim_{n \to \infty}\frac{(\ceil{sp^mp^n}+p^n)^{\mu-1}}{\ceil{sp^mp^n}^{\mu-1}}\\
&=\lim_{m \to \infty}\frac{(sp^m+1)^{\mu-1}}{(sp^m)^{\mu-1}}\\
&=1.
\end{split}
\end{equation*}
So
\begin{equation*}
\begin{split}
h'_+(s)&=\lim_{m \to \infty}\lim_{n \to \infty}\frac{\sum\limits^{\ceil{sp^mp^n}+p^n-1}_{t=\ceil{sp^mp^n}}D_{m+n,t}}{p^{m(d-1)}p^{nd}}\\
&\leq \underline{\lim}_{m \to \infty}\underline{\lim}_{n \to \infty}\frac{p^nD_{m+n,\ceil{sp^mp^n}}}{p^{m(d-1)}p^{nd}}\frac{\alpha_{\mu,\ceil{sp^mp^n}+p^n}}{\alpha_{\mu,\ceil{sp^mp^n}}}\\
&= \underline{\lim}_{m \to \infty}\underline{\lim}_{n \to \infty}\frac{p^nD_{m+n,\ceil{sp^mp^n}}}{p^{m(d-1)}p^{nd}}\\
&= \underline{\lim}_{m \to \infty}\underline{\lim}_{n \to \infty}\frac{D_{m+n,\ceil{sp^mp^n}}}{p^{m(d-1)}p^{n(d-1)}} \ .
\end{split}
\end{equation*}
For a sequence of real numbers $\beta_n$ and any $m$, $\underline{\lim}_{n \to \infty}\beta_{m+n}=\underline{\lim}_{n \to \infty}\beta_{n}$ is independent of $m$, so $\underline{\lim}_{m \to \infty}\underline{\lim}_{n \to \infty}\frac{D_{m+n,\ceil{sp^mp^n}}}{p^{m(d-1)}p^{n(d-1)}}=
\underline{\lim}_{n \to \infty}\frac{D_{n,\ceil{sp^n}}}{p^{n(d-1)}}$. Therefore we have
$$h'_+(s) \leq \underline{\lim}_{n \to \infty}\frac{D_{n,\ceil{sp^n}}}{p^{n(d-1)}}= \underline{\lim}_{n \to \infty}\frac{f_n(s)}{p^{n(d-1)}}.$$
The proof of the last inequality is similar. First we have
If $\ceil{sp^mp^n}-p^n \leq t \leq \ceil{sp^mp^n}-1$, then $D_{m+n,t}/\alpha_{\mu,t} \geq D_{m+n,\ceil{sp^mp^n}}/\alpha_{\mu,\ceil{sp^mp^n}}$, so

\begin{equation*}
\begin{split}
D_{m+n,t}& \geq D_{m+n,\ceil{sp^mp^n}}\frac{\alpha_{\mu,t}}{\alpha_{\mu,\ceil{sp^mp^n}}}\\
& \geq D_{m+n,\ceil{sp^mp^n}}\frac{\alpha_{\mu,\ceil{sp^mp^n}-p^n}}{\alpha_{\mu,\ceil{sp^mp^n}}} \ .
\end{split}
\end{equation*}
Also $\alpha_{\mu,t}$ is a polynomial of degree $\mu-1$ in $t$, so
\begin{equation*}
\begin{split}
\lim_{m \to \infty}\lim_{n \to \infty}\frac{\alpha_{\mu,\ceil{sp^mp^n}-p^n}}{\alpha_{\mu,\ceil{sp^mp^n}}}
&=\lim_{m \to \infty}\lim_{n \to \infty}\frac{(\ceil{sp^mp^n}-p^n)^{\mu-1}}{\ceil{sp^mp^n}^{\mu-1}}\\
&=\lim_{m \to \infty}\frac{(sp^m-1)^{\mu-1}}{(sp^m)^{\mu-1}}\\
&=1.
\end{split}
\end{equation*}
So
\begin{equation*}
\begin{split}
h'_-(s)&=\lim_{m \to \infty}\lim_{n \to \infty}\frac{\sum\limits^{\ceil{sp^mp^n}-1}_{t=\ceil{sp^mp^n}-p^n}D_{m+n,t}}{p^{m(d-1)}p^{nd}}\\
&\geq \overline{\lim}_{m \to \infty}\overline{\lim}_{n \to \infty}\frac{p^nD_{m+n,\ceil{sp^mp^n}}}{p^{m(d-1)}p^{nd}}\frac{\alpha_{\mu,\ceil{sp^mp^n}-p^n}}{\alpha_{\mu,\ceil{sp^mp^n}}}\\
&= \overline{\lim}_{m \to \infty}\overline{\lim}_{n \to \infty}\frac{p^nD_{m+n,\ceil{sp^mp^n}}}{p^{m(d-1)}p^{nd}}\\
&= \overline{\lim}_{m \to \infty}\overline{\lim}_{n \to \infty}\frac{D_{m+n,\ceil{sp^mp^n}}}{p^{m(d-1)}p^{n(d-1)}} \ .
\end{split}
\end{equation*}

\noindent  For a sequence of real numbers $\beta_n$ and any $m$, $\overline{\lim}_{n \to \infty}\beta_{m+n}=\overline{\lim}_{n \to \infty}\beta_{n}$ is independent of $m$, so $\overline{\lim}_{m \to \infty}\overline{\lim}_{n \to \infty}\frac{D_{m+n,\ceil{sp^mp^n}}}{p^{m(d-1)}p^{n(d-1)}}=
\overline{\lim}_{n \to \infty}\frac{D_{n,\ceil{sp^n}}}{p^{n(d-1)}}$. Therefore we have
$$h'_-(s) \geq \overline{\lim}_{n \to \infty}\frac{D_{n,\ceil{sp^n}}}{p^{n(d-1)}}= \overline{\lim}_{n \to \infty}\frac{f_n(s)}{p^{n(d-1)}}.$$

\noindent (2) If $h_M$ is differentiable at $s$, $h'_+(s)= h'_{-}(s)$. Thus (1) implies that $f_{n,M}(s)/q^{d-1}$ exists and is equal to $h'(s)$, rest of (2) is clear.\\

\noindent (3) follows from \Cref{th: one sided differentiability of h-functions}, (3).
\end{proof}

\begin{remark}\label{re: density function in generality}
 We prove \Cref{th: differentiability of h function implies existence of density function} in the context of situation (A) or (B) defined in \Cref{th: uniform convergence to the convex functional}- which is precisely the contexts where we prove existence of $h_{M,I,J_{\bullet}}$ in this article. Thus when $(R, \mathfrak m)$ is a domain, $I, J_{\bullet}$ satisfies \textbf{Condition C}, we get a corresponding density function which is well-defined outside a countable subset of $(0, \infty)$. One particular special case, potentially important for its application to prime characteristic singularity theory, is when $J_{\bullet}$ is the ideal sequence that defines the $F$-signature of $(R,\mathfrak m)$; see \Cref{eg: when J is big}.

 When $J_n= \bpq{J}{q}$, \Cref{th: differentiability of h function implies existence of density function} yields a Hilbert-Kunz density function of $(I,J)$ well defined outside a countable subset of $(0, \infty)$.

 The function $h_{M,I, J_{\bullet}}$ need not be continuous or differentiable at zero. In \Cref{th: asymptotic behaviour of h function near zero}, we prove that for a local domain $R$, $h_{R,I,J}$ is continuous at zero if and only if $\dim R-\dim R/I \geq 1$ and differentiable at zero if and only if $\dim R-\dim R/I \geq 2$.
\end{remark}

\subsection{Properties of density functions}\label{sse: properties of local density functions}

The following consequence of \Cref{th: one sided differentiability of h-functions} and \Cref{th: differentiability of h function implies existence of density function} is used in \Cref{th: density function near zero}. For the notions of integral closure and analytic spread appearing below, we refer to \cite{HunekeSwanson}.

\begin{theorem}\label{pr: monotonicity result for the desnity function}
    Let $(R, \mathfrak{m}), I,J_{\bullet},M$ be as in \Cref{th: uniform convergence to the convex functional}. Let $r$ be an integer greater than the analytic spread of $I$. Denote the right and left hand derivatives of $h_{M,I,J_{\bullet}}$ at $s>0$ by $h^{'}_{+}(s)$ and $h^{'}_{-}(s)$ respectively. Then,
    \begin{enumerate}
        \item Both the functions $\frac{h^{'}_{+}(s)}{s^{r-1}}$ and $\frac{h^{'}_{-}(s)}{s^{r-1}}$ are decreasing on $(0, \infty)$.

        \item For positive real numbers $s_1\leq s_2$, $\frac{h^{'}_{-}(s_2)}{s_2^{r-1}}\leq \frac{h^{'}_{+}(s_1)}{s_1^{r-1}}$. Moreover, if $f_{M,I,J_\bullet}(s_1), f_{M,I,J_\bullet}(s_2)$ exist, we have

         $$\frac{f_{M, I, J_\bullet}(s_2)}{s_2^{r-1}} \leq \frac{f_{M, I, J_\bullet}(s_1)}{s_1^{r-1}}.$$
    \end{enumerate}
\end{theorem}

\begin{proof}
    Once we prove either (1) or (2) for a certain $r$, the corresponding assertion for a larger value of $r$ follows immediately. So we assume $r$ is the analytic spread of $I$. Without loss of generality we can assume that the residue field of $R$ is infinite. So there exists $r$ elements $x_1, x_2,\ldots, x_r$ of $I$ such that the integral closure of $(x_1, x_2, \ldots, x_r)$ is $I$. 
    
    We prove (1) now. By \Cref{pr: invariance under closure operations}, we can assume $I$ is generated by $r$ elements.  By \Cref{th: uniform convergence to the convex functional}, assertion (2) and \Cref{th: one sided differentiability of h-functions}, assertion (4), $\frac{h^{'}_{+}(s)}{s^{r-1}}$ and $\frac{h^{'}_{-}(s)}{s^{r-1}}$ are right and left hand derivatives of a convex function respectively. Since the right and left hand derivatives of a convex function are always decreasing, we prove (1).\\

    \noindent (2) The inequality $\frac{h^{'}_{-}(s_2)}{s_2^{r-1}}\leq \frac{h^{'}_{+}(s_1)}{s_1^{r-1}}$
    is a consequence of \Cref{th: one sided differentiability of h-functions}, (2), (4).
    This inequality combined with \Cref{th: differentiability of h function implies existence of density function}, (1), implies
    
    $$\frac{f_{I, J_\bullet}(s_2)}{s_2^{r-1}} \leq \frac{h^{'}_{-}(s_2)}{s_2^{r-1}}\leq \frac{h^{'}_{+}(s_1)}{s_1^{r-1}} \leq \frac{f_{I, J_\bullet}(s_1)}{s_1^{r-1}}.$$
    \end{proof}

As a consequence, we prove results about one sided limits of the left, right derivatives of $h_{I,J_{\bullet}}$ and the density function $f_{I, J_{\bullet}}$:
\begin{proposition}\label{pr: one sided limits exists}
    With the same notation as in \Cref{pr: monotonicity result for the desnity function},  for $s_0>0$,
    \begin{enumerate}
        \item $\underset{s \to s_0-}{\lim}h^{'}_{+}(s)= \underset{s \to s_0-}{\lim}h^{'}_{-}(s)$ and $\underset{s \to s_0+}{\lim}h^{'}_{+}(s)= \underset{s \to s_0+}{\lim}h^{'}_{-}(s)$.

        \item Given a sequence of positive real numbers $(s_n)_n$ converging to $s_0$ such that $s_n < s_0$ for all $n$, and $f_{I, J_{\bullet}}(s_n)$ exist for all $n$, we have
        $$\underset{n \to \infty}{\lim}f_{I, J_{\bullet}}(s_n)= \underset{s \to s_0-}{\lim}h^{'}_{+}(s).$$

        \item Given a sequence of positive real numbers $(s_n)_n$ converging to $s_0$ such that $s_n > s_0$ for all $n$ and $f_{I, J_{\bullet}}(s_n)$ exist for all $n$, we have
        $$\underset{n \to \infty}{\lim}f_{I, J_{\bullet}}(s_n)= \underset{s \to s_0+}{\lim}h^{'}_{+}(s).$$
    \end{enumerate}
\end{proposition}

\begin{proof}
    (1) Because of similarity in the arguments, we only prove the assertion involving left hand limits. The limit of each of $h^{'}_{+}(s)/s^{r-1}$ and $h^{'}_{-}(s)/s^{r-1}$ as $s$ approaches $s_0$ from the left exists since these functions are decreasing by \Cref{pr: monotonicity result for the desnity function}, (1). The first inequality in assertion (2), \Cref{pr: monotonicity result for the desnity function} implies 
    $$\underset{s \to s_0-}{\lim}\frac{h^{'}_{+}(s)}{s^{r-1}} \geq \underset{s \to s_0-}{\lim}\frac{h^{'}_{-}(s)}{s^{r-1}}.$$
    So $\underset{s \to s_0-}{\lim}h^{'}_{+}(s) \geq \underset{s \to s_0-}{\lim}h^{'}_{-}(s).$  Since $h^{'}_{+}(s) \leq h^{'}_{-}(s)$ for positive $s$, by \Cref{th: one sided differentiability of h-functions}, (4), the desired equality of the left hand limits in (1) follows.\\

    \noindent (2) For a positive $s$, whenever $f_{I,J_{\bullet}}$ exists, by \Cref{th: differentiability of h function implies existence of density function}, (1), we have
    $h^{'}_{+}(s) \leq f_{I,J_{\bullet}}(s) \leq h^{'}_{-}(s)$, whence it follows 
    
    $$\underset{n \to \infty}{\textup{limsup}}\,f_{I,J_{\bullet}}(s_n) \leq \underset{s \to s_0-}{\lim}h^{'}_{-}(s)= \underset{s \to s_0-}{\lim}h^{'}_{+}(s) \leq \underset{n \to \infty}{\textup{liminf}}\, f_{I,J_{\bullet}}(s_n).$$
    This proves (2). The argument for (3) is similar.
\end{proof}

\subsection{Examples}\label{sse: examples}
Our next result allows constructions of $h$-functions and density functions with a prescribed order of smoothness, as detailed in \Cref{eg: examples of differentiable density function}. In the sequel, \Cref{th: HK density when J is homogeneous but not of finite colength} shows that, when $(R, \mathfrak{m})$ is standard graded of dimension at least two, and $J$ is homogeneous, the corresponding density function $f_{\mathfrak{m},J}$ exits and is continuous. \Cref{eg: examples of differentiable density function} contrasts those results by producing continuous density functions in the not necessarily graded setup. 

\begin{theorem}\label{th: h function and integration}

Let $R$ be a local ring of characteristic $p$, $R[[t]]$ be a power series ring with indeterminate $t$. Let $M$ be a finitely generated $R$-module, $I,J$ be two $R$-ideals such that $I+J$ is $\mathfrak{m}$-primary. Let $M[[t]]=M \otimes_R R[[t]]$. Then

\begin{enumerate}
\item
$h_{M[[t]],R[[t]],(I,t^\alpha),(J,t^\beta)}(s)=\alpha\int_{s-\beta/\alpha}^s h_{M,R,I,J}(x)dx$.

\item $h_{M[[t]],R[[t]],(I,t^\alpha),J}(s)=\alpha\int_{0}^s h_{M,R,I,J}(x)dx$.

\item $h_{M[[t]],R[[t]],I,(J,t^\beta)}(s)=\beta h_{M,R,I,J}(s)$.
\end{enumerate}
\end{theorem}

\begin{proof}We will use the convention $I^s=R$ when $s \leq 0$. To prove the equality we may assume $s=s_0/q_0 \in \mathbb{Z}[1/p]$ because the functions on both sides are continuous when $s>0$. Then for $q \geq q_0$, $sq$ is an integer.
$$h_{n,M[[t]],R[[t]],(I,t^\alpha),(J,t^\beta)}=l(\frac{M[[t]]}{((I,t^\alpha)^{sq}+(J^{[q]}, t^{\beta q}))M[[t]]})$$
The above length is also equal to
$$l(\frac{M[[t]]}{(\sum_{0 \leq j \leq sq}I^{sq-j}t^{\alpha j}+(J^{[q]}, t^{\beta q}))M[[t]]})$$
But by the convention, it is also
$$l(M[[t]]/\sum_{0 \leq j \leq \infty}I^{sq-j}t^{\alpha j}+(J^{[q]}, t^{\beta q})M[[t]])$$
and because the existence of the $t^{\beta q}$-term, it is also equal to
$$l(M[[t]]/(\sum_{0 \leq j \leq \lfloor \beta q/\alpha \rfloor}I^{sq-j}t^{\alpha j}+(J^{[q]}, t^{\beta q}))M[[t]])$$
Note that the module inside is nonzero only in $t$-degree at most $\beta q-1$. So summing up over the lengths in different $t$-degrees, the above length is also equal to the following sum:
$$L=\sum_{0 \leq x \leq \beta q-1} l(M/(J^{[q]}+I^{sq-\lfloor x/\alpha \rfloor})M)$$
Let $y=\lfloor x/\alpha \rfloor$ and
$$L_1=\sum_{\alpha\lfloor \beta q/\alpha \rfloor \leq x \leq \beta q-1} l(M/(J^{[q]}+I^{sq-\lfloor x/\alpha \rfloor})M),$$
$$L_2=\sum_{0 \leq x \leq \alpha\lfloor \beta q/\alpha \rfloor-1} l(M/(J^{[q]}+I^{sq-\lfloor x/\alpha \rfloor})M).$$
Here we denote $L_1=0$ if $\beta q/\alpha \in \mathbb{Z}$. Then $L=L_1+L_2$, and $L_1$ has at most $\alpha$ terms and each term is of order $O(q^d)$, so $L_1$ has order $O(q^d)$. Now
\begin{align*}
L_2=\alpha\sum_{0 \leq y \leq \lfloor \beta q/\alpha \rfloor-1} l(M/J^{[q]}+I^{sq-y}M)\\
=\alpha\sum_{0 \leq y \leq \lfloor \beta q/\alpha \rfloor-1} h_{n,M,I,J}(s-y/q)\\
=\alpha q \int^s_{s-\lfloor \beta q/\alpha \rfloor/q}h_{n,M,I,J}(x)dx
\end{align*}
Now $\lim_{q \to \infty}L_1/q^{d+1}=0$, so
$$\lim_{q \to \infty}L/q^{d+1}=\lim_{q \to \infty}L_2/q^{d+1}=\alpha \int^s_{s-\beta/\alpha}h_{M,I,J}(x)dx.$$
Since the equation
$$h_{M[[t]],R[[t]],(I,t^\alpha),(J,t^\beta)}=\alpha\int_{s-\beta/\alpha}^s h_{M,R,I,J}(x)dx$$
is true on $\mathbb{Z}[1/p]$ and both sides are continuous with respect to $s$, they are equal on all of $\mathbb{R}$. The rest of the two equations can be obtained by taking limit as $\alpha$ or $\beta$ goes to infinity and using the $\mathfrak{m}$-adic continuity proven in \Cref{th: m-adic continuity of h function}. 
\end{proof}

\begin{example}\label{eg: examples of differentiable density function}
Given a noetherian local ring $(R, \mathfrak{m})$ of characteristic $p$, and ideals $I$, $J$ such that $I+J$ is $\mathfrak {m}$-primary, consider the power series ring $S=R[[t_1, t_2, \ldots, t_n]]$ for $n \geq 2$ and the ideals $I'= (I, t_1, \ldots, t_n)$, $J'= (J, t_1, \ldots, t_n)$. The theorem above and the continuity of $h_{R,I,J}$ on $R_{>0}$ implies that $h_{S,I',J'}$ is $(n-1)$-times differentiable on $\mathbb{R}_{>0}$ and the $(n-1)$-th derivative is continuous on $\mathbb{R}_{>0}$. Therefore the density function $f_{S,I',J'}$ exists everywhere on $\mathbb{R}_{>0}$ and is $(n-2)$-times differentiable on $\mathbb{R}_{>0}$, where the $(n-2)$-th order derivative is continuous, by \Cref{th: differentiability of h function implies existence of density function}, (2). When $R$ is a domain and $I= \mathfrak{m}$, the density function $f_{S,I',J'}$ is continuous at zero by \Cref{th: density function near zero}, (2).
\end{example}

Although the differentiability of the $h$-function at a point is sufficient for the existence of the density function at that point, the differentiability is not a necessary condition. The next few examples illustrate some related subtleties.

\begin{example}\label{eg: h function is not differentiable at a point}
    We point out that the $h$-function need not be differentiable on $(0, \infty)$. Our example of a non differentiable $h$-function comes from \cite{BlicleSchwedeTuckerFractal}. Fix a regular local domain $(R,\mathfrak m)$ of dimension $d$ and a nonzero $f \in R$. For $t \in \mathbb R$, \cite{BlicleSchwedeTuckerFractal} considers the function $t \to s(R,f^t)$:  the $F$-signature of the pair $(R,f^t)$ which is shown to be the same as
    \[s(R,f^t)= \lim_{n \to \infty} \frac{1}{q^d}l(\frac{R}{\bp{\mathfrak m}{n}:f^{\ceil{tp^n}}}).\]
    With $I=(f)$, $h_{R,I,\mathfrak m}(t)=1- s(R, f^t)$; see \cite[section 4]{BlicleSchwedeTuckerFractal}. At $t=1$, the left hand derivative of $h_{R,I, \mathfrak{m}}$ is the $F$-signature of $R/f$; see \cite[Theorem 4.6]{BlicleSchwedeTuckerFractal}, while the right hand derivative is zero since $h(s)=1$ for $s \geq 1$. So $h$ is not differentiable at one if and only if the $F$-signature of $R/f$ is nonzero, precisely when $R/f$ is strongly $F$-regular. A concrete example comes from the strongly $F$-regular ring, $\mathbb{F}_p[[x,y,z]]/(x^2+y^2+z^2)$ with $p \geq 3$. A direct calculation shows that the density function $f_{R,I, \mathfrak{m}}$ exists at $s=1$ and has value zero. So when when $R/f$ is strongly $F$-regular, at $s=1$, $f_{R,I, \mathfrak{m}}$ coincides with the right derivative of $h_{R,I,J}$ and differs from the left derivative of $h_{R,I,J}$.
\end{example}

\begin{example}\label{eg: where the limit defining the desnity function does not exist}
    We point out that the limit defining the density function at a particular $s \in \mathbb R$, i.e. of $f_{n,M,I,J}(s)/q^{\dim(M)-1}$ may not converge. For example, when $I=0$, $M=R$, then $f_{n,M,I,J}(0)=l(R/J^{[q]})$; thus $f_{n,M,I,J}(0)/q^{\dim R}=e_{HK}(J,R)$ is a nonzero real number, so $f_{n,M,I,J}(0)/q^{\dim R-1}$ goes to infinity. This example implies that some assumption is necessary to guarantee the existence of the density function at every point.
\end{example}

\begin{example}
In the definition of the density function if we replace $\ceil{sq}$ by $\flor{sq}$, then we have more examples where the density function does not exist. We recall Otha's example mentioned in \cite[sec 3]{kosuke2017function} which produces such instances. Let $R$ be the power series ring $k[[x_1,\ldots,x_{d+1}]]$, $\alpha_1 \leq \ldots \leq \alpha_{d+1}$ be a sequence of positive integers, $I=(x_1^{\alpha_1}\ldots x_{d+1}^{\alpha_{d+1}})$ be a monomial principal ideal, $J=(x_1,\ldots,x_{d+1})$ be the maximal ideal of $R$. Assume moreover that $\alpha_d < \alpha_{d+1}$, $\alpha_{d+1}$ does not divide $p$,and $\epsilon_n \in [0, \alpha_{d+1}-1]$ is the residue of $p^n$ modulo $\alpha_{d+1}$. Let $\tilde{f}$ be the density function defined using $\flor{sq}$, then $\lim_{n \to \infty} \tilde{f}_{n,R,I,J}(\frac{1}{\alpha_{d+1}})/(p^{nd}\epsilon_n)$ exists and is nonzero. So $\lim_{n \to \infty} \tilde{f}_{n,R,I,J}(\frac{1}{\alpha_{d+1}})/p^{nd}$ exists if and only if $\epsilon_n$ is a constant sequence- this is false in general. In general, $\epsilon_n$ is a periodic function and its period is the order of $p+\alpha_{n+1}\mathbb{Z}$ in the multiplicative group $(\mathbb{Z}/\alpha_{n+1}\mathbb{Z})^*$.
\end{example}

\begin{example}\label{eg: density functions exist but is not continuous}
We give an example, where the density function exists everywhere although the $h$-function is not differentiable everywhere. Note that the resulting density function is not continuous in this case; compare with \Cref{th: h is the integral of f}. Let $M=R=k[[x]]$ be the power series ring, $ I=J=(x)$. Then $h_n(s)=l(R/I^{\lceil sq \rceil}+J^{[q]})=\max\{\min\{\lceil sq \rceil, q\},0\}$. By simple calculation we get $f_{n}(s)=1$ when $-1/q<s \leq 1-1/q$ and is $0$ otherwise. So $f(s)=1$ when $0 \leq s<1$ and $f(s)=0$ otherwise.

Here $f_{n}$ converges pointwise but not uniformly. Outside an arbitrary neighborhood of $0$ and $1$ then $f_{n}$ converges uniformly.

On the other hand, $h(s)$ is $0$ when $s\leq 0$, $s$ when $0 \leq s \leq 1$, $1$ when $s \geq 1$, and is continuous. We have $f(s)=h'(s)$ when $s \neq 0,1$; when $s=0,1$ $h'(s)$ does not exist and $f(s)=h'_+(s)$. 
\end{example}

\Cref{eg: h function is not differentiable at a point}, \Cref{eg: density functions exist but is not continuous} leads us to guessing that whenever the density function exists at a positive number $s$, it coincides with the right hand derivative $h'_+(s)$; see \Cref{qe: continuity of density function}. 

\begin{remark}\label{re: analyzing monotonicity of density function}
    Assume $J_{\bullet}$ is big and $h_{M,I, J_{\bullet}}$ is differentiable everywhere. Since $h_{M,I, J_{\bullet}}$ is eventually constant (\Cref{le: compactly supported h-function}), the resulting density function $f_{M,I,J_\bullet}= h'_{M,I,J_\bullet}$ is supported on some compact interval $[0,b]$. So the density function has to increase and decrease on $[0,b]$. By \Cref{th: one sided differentiability of h-functions}, $f_{M,I,J_\bullet}= h'(s)= \mathcal{H}'(s)s^{\mu-1}/(\mu-1)!$, where $\mathcal H'$ is decreasing since $\mathcal H$ is convex; so this gives a natural way to represent $f_{M,I,J_\bullet}$ as a product of a decreasing and an explicit increasing function, namely $c(s)$. This may help analyzing the monotonicity of the density function.
\end{remark}

\section{Relation among $h$, density, and Frobenius-Poincar\'e functions}\label{se: relation among different function}

In \Cref{sse: Frob-Poincare in the local setting}, we developed a notion of Frobenius-Poincar\'e function in the local setting. Work of \Cref{se: differentiability of h function} gives a notion of Hilbert-Kunz density function in the local setting, at least outside a countable subset of $(0, \infty)$. When $(R,\mathfrak m)$ is graded, we prove that these local notions defined using the $\mathfrak m$-adic filtration coincide with the classical notion of Frobenius-Poincar\'e function and Hilbert-Kunz density function defined (see \Cref{se: background material}) using the graded structures of the underlying objects; see \Cref{th: h function differentiable when J is m primary}, \Cref{pr: the graded Frob-Poincare coincides wth the local Frob Poincare}. Moreover, using the equivalence stated in  \Cref{th: continuous differentiability vs continuity of the density function}, we prove continuous differentiability of $h$-functions of two or higher dimensional standard graded rings in \Cref{th: HK density when J is homogeneous but not of finite colength}.  
The next few lemmas are used in \Cref{th: h function differentiable when J is m primary}, \Cref{th: HK density when J is homogeneous but not of finite colength} to compare the density function  defined using the $\mathfrak{m}$-adic filtration in the graded setup and the classical one using the underlying grading. 

\begin{lemma}\label{le: degenerated in degree zero implies filtration coincides with grading}
Let $(R,\mathfrak m)$ be a standard graded ring, $M$ be a finitely generated $\mathbb Z$-graded module of dimension $d$, $J$ be a homogeneous ideal of finite colength. Set 
\[g_{n,M,J,d-1}(s)=\frac{1}{q^{d-1}}l(\frac{M}{J^{[q]}M})_{\lceil sq \rceil}, \,\, g_{n,M,J}(s)= l(\frac{M}{J^{[q]}M})_{\lceil sq \rceil} \ .\]
\begin{enumerate}
    \item When $M$ is generated in degree zero, for any graded submodule $N \subseteq M$, $(M/N)_j=\mathfrak{m}^j(M/N)/\mathfrak{m}^{j+1}(M/N)$.

    \item When $M$ is generated in degree zero, $ g_{n,M,J}(s)= l(\frac{M}{J^{[q]}M})_{\lceil sq \rceil}= f_{n,M,\mathfrak m,J}(s)$.
\end{enumerate}
\end{lemma}
\begin{proof}Let $N$ be any submodule of $M$, then $M/N$ is also generated in degree $0$, so $(M/N)_{\geq j}=\mathfrak{m}^j(M/N)$ and $(M/N)_j=\mathfrak{m}^j(M/N)/\mathfrak{m}^{j+1}(M/N)$ for any $j$. This implies $g_{n,M,J}(s)=f_{n,M,\mathfrak{m},J}(s)$.
\end{proof}

\begin{lemma}\label{le: any graded module is equivalent to one generated in degree zero}We define an equivalence relation $\sim$ on $\mathbb Z$-graded modules over a standard graded ring $R$ of positive dimension over a field: we say $M \sim N$ when there is a homogeneous map $\phi: M \to N$ such that $\dim \textup{Ker}(\phi), \dim \textup{Coker}(\phi) \leq \dim R-1$, and let $\sim$ also denote the minimal equivalence relation generated by such relations. Then $M$ is equivalent to some module generated in degree 0.
\end{lemma}

\begin{proof}Since $\dim R>0$, there is a positive integer $t$ and an element $c \in R_t$ such that $\dim R/cR<\dim R$. First, we find a sufficient large $n>0$ such that $M$ is generated in degree at most $nt$. Then we truncate at degree $nt$ to get $M_{\geq nt}:= \oplus_{j=nt}^{\infty} M_j$, which is generated in degree $nt$. The module $M/M_{\geq nt}$ is Artinian. The inclusion $M_{\geq nt} \hookrightarrow M$ shows $M_{\geq nt} \sim M$. The map $M_{\geq nt} \to M_{\geq nt}[nt]$ given by multiplication by $c^n$ has its kernel and cokernel annihilated by $c^n$. So the kernel and cokernel have dimension less than $\dim R$. Thus $M \sim M_{\geq nt} \sim M_{\geq nt}[nt]$ . Since $M_{\geq nt}[nt]$ is generated in degree zero, we are done.
\end{proof}

\noindent The next result follows directly from the lemma above and \Cref{pr: associativity formula for h function}.

\begin{lemma}\label{le: for h functions we can assume the module is generated zero}
 Let $(R, \mathfrak m)$ be standard graded, $M$ be a finitely generated $\mathbb{Z}$-graded $R$-module, $I$, $J_{\bullet}$ be homogeneous; assume that the corresponding objects obtained by localizing at $\mathfrak m$ satisfy condition (A) or (B) stated in \Cref{th: uniform convergence to the convex functional}. Then there is a finitely generated $\mathbb N$-graded $R$-module $M'$ generated in degree zero such that,
 $h_{M,I,J_{\bullet}}= h_{M',I,J_{\bullet}}$.
\end{lemma}

 In the context of (A) or (B) stated in \Cref{th: uniform convergence to the convex functional}, there is an $h$-function and an associated density function defined outside a countable subset of $(0, \infty)$. Although the limit defining the density function may not exist at every point of $(0, \infty)$, we can define the integral of $f_{M,I,J_{\bullet}}$ on any bounded measurable subset $\Sigma$ of $[0, \infty)$ by integrating the class in $L^1(\Sigma)$ represented by the density function. Indeed, being locally constant, each $f_{n,M,I,J_{\bullet}}$ is measurable. Let $S$ be the subset of $\Sigma$, where the limit of $(f_{n,M,I,J_{\bullet}}(s)/q^{\dim(M)-1})_n$ exists. Since $S$ has at most countable complement by \Cref{th: differentiability of h function implies existence of density function}, $S$ is measurable and $f_{M,I,J_{\bullet}}$ being the limit of measurable functions defines a measurable function on $S$. Since any countable set is measurable, any function $g: \Sigma \rightarrow \mathbb{R}$ extending $f_{M,I,J_{\bullet}}$ is automatically measurable. Then we define
 $\int_{\Sigma}f_{M,I,J_{\bullet}}(t)dt= \int_{\Sigma}g(t)dt$, since the later integral is independent of the choice of $g$. The next result implies that $f_{M,I,J_{\bullet}}$ is in fact integrable on $\Sigma$.

\begin{theorem}\label{th: h is the integral of f}
Let $(R,\mathfrak m)$, $I$, $J_{\bullet}$, $M$ be as in \Cref{th: uniform convergence to the convex functional}. Then in the context of situation (A) or (B) as stated in \Cref{th: uniform convergence to the convex functional}, we have for any $s>s_0>0$,
$$h_{M,I,J_{\bullet}}(s)-h_{M,I,J_{\bullet}}(s_0)=\int_{s_0}^{s}f_{M,I,J_{\bullet}}(t)dt.$$
Moreover if the density $f_{M,I,J_{\bullet}}$ exists and is continuous \footnote{See \Cref{de: continuity of density function} for the meaning of continuity of the density function at a point.} at $s>0$, then $h_{M,I,J_{\bullet}}$ is differentiable at $s$ and $f_{M}(s)=h'_{M}(s)$.
\end{theorem}

\begin{proof}
Note,

\begin{equation*}
        h_n(s)-h_n(s_0)= \sum \limits_{j=\ceil{s_0q}}^{\ceil{sq}-1}f_n(\frac{j}{q}) \ .
\end{equation*}

\noindent For any integer $j$ and $t \in (j-1/q, j/q]$, $f_n(t)= f_n(j/q)$. So the above equation yields
\begin{equation*}
    \frac{1}{q^d} h_n(s)-\frac{1}{q^d}h_n(s_0)= \int \limits_{\frac{\ceil{s_0q}}{q}-\frac{1}{q}}^{\frac{\ceil{sq}}{q}- \frac{1}{q}} \frac{f_{n}(t)}{q^{d-1}}dt \ .
\end{equation*}

\noindent Fix $0<a<b<0$ such that $[s_0,s] \subseteq (a,b)$. By \Cref{th: bound on the density function}, we can choose a constant $C$ such that for any $n\in \mathbb N$ and $t \in [a,b]$.
\[\frac{1}{q^{d-1}}f_n(t) \leq C .\]
Thus taking limit as $n$ approaches infinity and using dominated convergence, we get
\[h_{M,I,J_{\bullet}}(s)-h_{M,I,J_{\bullet}}(s_0)=\int_{s_0}^{s}f_{M,I,J_{\bullet}}(t)dt.\]

Whenever $f_M(t)$ exists at $s$ and is continuous at $s$, the differentiability of $h_M$ at $s$ and that $h'_M(s)= f_M(s)$ follows from the second fundamental theorem for integrable functions as stated in \cite[Thm 7.11]{Rudin}. Indeed, a point of continuity of $f_{M,I,J_{\bullet}}(t)$ is a Lebesgue point of the $L^1$- function $f_{M,I,J_{\bullet}}(t)$ as meant in \cite[Thm 7.11]{Rudin}.
\end{proof}

\begin{remark}\label{re: failure of first fundamental theorem for measurable function}
    We warn the reader that the above theorem is not simply a consequence of the fact that $f_{M,I, J_{\bullet}}(x)= h'_{M,I, J_{\bullet}}(x)$ almost everywhere; see \cite[7.16, Thm 7.18]{Rudin}. In fact, the theorem above implies absolute continuity of $h_{M,I, J_{\bullet}}$; see \cite[Def 7.17, Thm 7.18]{Rudin}.
\end{remark}

The continuity of the density function at a point (see \Cref{de: continuity of density function}) is in fact equivalent to the continuous differentiability of the corresponding $h$-function at that point, as interpreted in \Cref{de: continuous differentiability of a function} below.

\begin{definition}\label{de: continuous differentiability of a function}
    Let $\phi: \mathbb{R} \rightarrow \mathbb{R}$ be a continuous function. Let $\Lambda$ be the set of all points where $\phi$ is differentiable, endowed with the subspace topology from $\mathbb{R}$. Given an $s \in \Lambda$, we say that \textit{$\phi$ is continuously differentiable at $s$}, if the derivative function $\phi': \Lambda \rightarrow \mathbb{R}$ is continuous at $s$.
\end{definition}

\begin{theorem}\label{th: continuous differentiability vs continuity of the density function}
  Let $R, I, J_{\bullet}, M$ be as in \Cref{th: uniform convergence to the convex functional}. Given $s_0 \in \mathbb{R}_{>0}$, the function $h_{M,I,J_{\bullet} }$ is continuously differentiable at $s_0$, as meant in \Cref{de: continuous differentiability of a function}, if and only if the density function $f_{M,I, J_{\bullet}}$ exists at $s_0$ and is continuous at $s_0$ (see \Cref{de: continuity of density function}).
\end{theorem}

\begin{proof}
    If the density function exists at $s_0$ and is continuous at $s_0$, the continuous differentiability of $h_{M,I, J_{\bullet}}$ is proven in \Cref{th: h is the integral of f}. For the other direction: The existence of $f_{M,I, J_{\bullet}}(s_0)$ follows from \Cref{th: differentiability of h function implies existence of density function}. To prove the continuity of $f_{M,I, J_{\bullet}}$, take a sequence of positive numbers $(s_n)_n$ approaching $s_0$ from either left or right. By \Cref{pr: one sided limits exists}, (2), (3), the sequence $(f_{M,I, J_{\bullet}}(s_n))_n$ converges. So we can compute the limit point by specializing to a subsequence. Choose a subsequence $(t_n)_n$ of $(s_n)$ such that $h_{M,I, J_{\bullet}}$ is differentiable at each $t_n$. Then by \Cref{pr: one sided limits exists}, (2), (3),
    $\underset{n \to \infty}{\lim} f_{M,I, J_{\bullet}}(s_n)= \underset{n \to \infty}{\lim} h'(t_n)$. Since $h'$ is continuous at $s_0$, 
    $$\underset{n \to \infty}{\lim} f_{M,I, J_{\bullet}}(s_n)= h'(s_0)= f_{M,I, J_{\bullet}}(s_0),$$
    where the last equality follows from \Cref{th: differentiability of h function implies existence of density function}, (2). This proves the continuity of $f_{M,I, J_{\bullet}}$ at $s_0$.
\end{proof}

\begin{proposition}\label{pr: matches with Trivedi's density function when J is m-primary}
Continue with the same notation as in \Cref{le: degenerated in degree zero implies filtration coincides with grading} but $M$ not necessarily generated in degree zero. Set
$$\tilde{g}_{n,M,J,d-1}(s)=l(M/J^{[q]}M)_{\lfloor sq \rfloor}/q^{d-1}.$$
If additionally $d=\dim(M) \geq 2$, the two limits below exist for all $s \in \mathbb R$:
$$\tilde{g}_{M,J}(s)=\lim_{n \to \infty}\tilde{g}_{n,M,J,d-1}(s),\, g_{M,J}(s)=\lim_{n \to \infty}g_{n,M,J,d-1}(s).$$
Moreover $\tilde{g}_{M,J}(s)=g_{M,J}(s)$.
\end{proposition}

\begin{proof}
By \cite{TriExist}, $\tilde g_{n,M,J,d-1}(s)$ converges for all $s \in \mathbb{R}$. For $s \in \mathbb{Z}[1/p]$, $g_{n,M,J,d-1}(s)= \tilde g_{n,M,J,d-1}(s)$ for $q$ large; so we conclude convergence of $g_{n,M,J,d-1}(s)$. When $s$ is not in $\mathbb{Z}[1/p]$, 
$$g_{n,M,J,d-1}(s)= \tilde g_{n,M,J,d-1}(s+\frac{1}{q}).$$
Now for $d \geq 2$, the uniform convergence of the sequence of functions $\tilde g_{n,M,J,d-1}$ and continuity of $\tilde g_{M,J}$ imply that the sequence $\tilde g_{n,M,J,d-1}(s+\frac{1}{q})$ converges to $\tilde{g}_{M,J}(s)$.
\end{proof}

\begin{theorem}\label{th: h function differentiable when J is m primary}
Let $(R,\mathfrak m)$ be standard graded, $J$ be a homogeneous $\mathfrak{m}$-primary ideal, $M$ be an $R$-module of dimension $d \geq 2$. Then
\begin{enumerate}
    \item $h_{M,\mathfrak m, J}$ is continuously differentiable on $\mathbb R$. The density function $f_{M,\mathfrak m,J}(s)$ exists everywhere on $\mathbb R$ and is the same as $h'_{M,\mathfrak m,J}(s)$. 

    \item  Moreover $f_{M,\mathfrak m,J}$ is the same as Trivedi's Hilbert-Kunz density function $\tilde{g}_{M,J}(s)$; see \Cref{se: background material}. 
\end{enumerate}
\end{theorem}

\begin{proof}(1) It follows from \cite[Lemma 3.3]{TaylorInterpolation}, that for $s \leq 1$, $h_M(s)= e(\mathfrak m,M)s^d/d!$. Since $d \geq 2$, $h_M$ is continuously differentiable at zero and the derivative is zero. A direct computation shows that the density function at zero exists and is zero. So we can restrict to $(0, \infty)$.
Thanks to \Cref{th: differentiability of h function implies existence of density function}, (2), it is enough to show that $h_M$ is differentiable on $(0, \infty)$. By using \Cref{le: for h functions we can assume the module is generated zero}, we can assume that $M$ is generated in degree zero. Thus by \Cref{le: degenerated in degree zero implies filtration coincides with grading}
\[f_{n,M,\mathfrak m,J}(s)= g_{n,M, J}(s):= l([\frac{M}{\bp{J}{n}M}]_{\ceil{sq}}) \, \text{for all}\, s \in \mathbb R.\]

\noindent As $d \geq 2$, by \Cref{pr: matches with Trivedi's density function when J is m-primary}, $g_{n,M, J}(s)/q^{d-1}$ converges to Trivedi's density function $\tilde{g}_{M,J}(s)$ for all $s$. Since $\tilde g_{M,J}(s)$ is continuous, $f_{M,\mathfrak{m},J}(s)$ is also continuous. Now by \Cref{th: h is the integral of f}, (2), $h_{M,I, J}$ is differentiable on $(0, \infty)$.\\

\noindent (2) Fix an $M'$ which is generated in degree zero and equivalent to $M$ in the sense of \Cref{le: any graded module is equivalent to one generated in degree zero}. Thanks to \Cref{le: for h functions we can assume the module is generated zero} and part (1)
\[h_{M}= h_{M'} \, , f_M= f_{M'}.\]

\noindent The associativity formula for Trivedi's density function implies $\tilde{g}_{M,J}= \tilde g _{M',J}$; see \cite[Prop 2.14]{TriExist}. Since $M'$ is generated in degree zero and has dimension at least two, by \Cref{le: degenerated in degree zero implies filtration coincides with grading} and \Cref{pr: matches with Trivedi's density function when J is m-primary}, $\tilde g _{M',J}= f_{M',\mathfrak m, J}$. Putting together we conclude that $f_{M, \mathfrak m, J}= \tilde g_{M,J}.$ 
\end{proof}

\noindent We further strengthen the above theorem by proving it for any homogeneous $J$ which not necessarily has finite colength,

\begin{theorem}\label{th: HK density when J is homogeneous but not of finite colength} Let $(R,\mathfrak{m})$ be standard graded, $J$ be a homogeneous ideal, $s \in \mathbb{R}$, $M$ be a finitely generated $\mathbb Z$-graded module of dimension $d$. Assume $d \geq 2$. Set $\tilde{g}_{n,M,J,d-1}(s)=l(M/J^{[q]}M)_{\flor{sq} }/q^{d-1}$. Then:

\begin{enumerate}
\item The sequence $(\tilde{g}_{n,M,J,d-1}(s))_n$ converges uniformly on every compact subset of $\mathbb R$. The limiting function is continuous.

\item $h_{M,\mathfrak m, J}$ is continuously differentiable and for $s \in \mathbb {R},$ 
\[h'_{M,\mathfrak m,J}(s)= f_{M,\mathfrak m, J}(s)= \lim_{n \to \infty} \tilde{g}_{n,M,J,d-1}(s).\]
\end{enumerate}
\end{theorem}

\begin{proof}
(1) For a positive integer $N$, set $J'= J+ \mathfrak m ^{N+1}$. Then on $[0,N]$, $\tilde{g}_{n,M,J,d-1}= \tilde g_{n,M,J',d-1}$. Since $J'$ is $\mathfrak m$-primary, by \cite{TriExist}, $\tilde{g}_{n,M,J',d-1}$ converges uniformly to a continuous function. Thus on $[0,N]$, $\tilde{g}_{n,M,J,d-1}$ converges uniformly to a continuous function.\\

\noindent (2) Fix a compact interval $[a,b] \subseteq \mathbb R$. By \Cref{th: m-adic continuity}, (1) we can choose $t_0$ such that for all $t \geq t_0$, $h_{M,\mathfrak m, J}= h_{M,\mathfrak m, J+\mathfrak m^{t}}$ on $[a,b]$. Using the ideas from the argument in part(1), fix an integer $t \geq t_0$ ensuring $\tilde{g}_{n,M,J,d-1}= \tilde g_{n,M,J+\mathfrak m ^t,d-1}$ on $[a,b]$ for all $n$. By \Cref{th: h function differentiable when J is m primary}, $h_{M,\mathfrak m, J+\mathfrak m^t}$ is differentiable on $\mathbb R$ with derivative $\tilde {g}_{M, J+ \mathfrak m ^t}$. Thus on $(a,b)$, $h_{M,\mathfrak m, J}$ is differentiable with derivative being the continuous function $\tilde g _{M,J}$. Since by \Cref{th: differentiability of h function implies existence of density function} $h'_M= f_M$ on $(a,b)$, we are done.   
\end{proof}

We point out below that in the graded context the Frobenius-Poincar\'e function defined using the underlying grading and the maximal ideal adic filtration coincide.
Recall that by $\Omega$, we denote the open lower half complex plane. Let $(R,\mathfrak m)$ be standard graded, $M$ be an $\mathbb{N}$-graded $R$-module, $J$ be a homogeneous ideal. For $y \in \Omega,$

\begin{proposition}\label{pr: the graded Frob-Poincare coincides wth the local Frob Poincare}
Let $(R,\mathfrak m)$ be standard graded, $M$ a finitely generated $\mathbb Z$-graded $R$-module of dimension $d$, $J$ be a homogeneous ideal. Consider the sequence of functions on the open lower half plane

\[G_{n,M,J}(y)= \sum_{j=0}^{\infty}l([\frac{M}{\bpq{J}{q}M}]_j)e^{-iyj/q} \ .\]

\begin{enumerate}
    \item $\frac{1}{q^d}G_{n,M,J}(y)$ defines a holomorphic function on $\Omega$ for every $n$.

    \item Recall that $F_{M,\mathfrak m, J}$ denotes the Frobenius-Poincar\'e function defined in \Cref{de: Frobenius-Poincare general version}. The sequence
    \[\lim_{n \to \infty}\frac{1}{q^d}G_{n,M,J}(y)\] converges to $F_{M,\mathfrak m, J}(y)$.

    \item When $J$ is $\mathfrak m$-primary, $G_{n,M,J}(y)/q^d$ converges to $F_{M,\mathfrak m ,J}(y)$ on $\mathbb C$.
\end{enumerate}
\end{proposition}

\begin{proof}
Fix an $\mathbb N$-graded module $M'$ generated in degree zero and equivalent to $M$ in the sense of \Cref{le: any graded module is equivalent to one generated in degree zero}.

(3) Since $J$ is $\mathfrak m$-primary, $G_n$ is a sum of finitely many entire functions. So $G_n$ is entire.\\
Fix a compact subset $K$ of $\mathbb C$. By \cite[Lemma 3.2.5]{AlapanThesis}, we can find a constant $D$ such that

\[|\frac{1}{q^d}G_{n,M,J}(y)- \frac{1}{q^d}G_{n,M',J}(y)| \leq \frac{D}{q} \, \text{for all}\, n \,\text{and}\,y \in K.\]

Since $M'$ is generated in degree zero, $F_{n,M',\mathfrak m,J}= G_{n,M',J}$. Since $F_{n,M', \mathfrak m,J}/q^d$ uniformly converges to $F_{M', \mathfrak m, J}$ on $K$, the last inequality implies that $\frac{1}{q^d}G_{n,M,J}$ converges uniformly to $F_{M', \mathfrak m, J}$ on $K$; see \Cref{th: exitence of Frob-Poincare for general ideals}. Thanks to \Cref{le: for h functions we can assume the module is generated zero} and \Cref{th: exitence of Frob-Poincare for general ideals},  $F_{M', \mathfrak m, J}= F_{M, \mathfrak m,J}$ on $\mathbb C$.\\

\noindent (1) There is a polynomial $P$ of degree $d$ with nonnegative coefficients such that 
\[l([\frac{M}{\bpq{J}{q}M}]_j) \leq l(M_j) \leq P(j).\]

\noindent Fix a compact subset $K \subseteq \Omega$. Choose $\epsilon>0$ such that $\Im y< -\epsilon$ for every $y \in K$. Since 

\[\sum_{j=0}^{\infty}\frac{1}{q^d}|P(j)|e^{-j\epsilon/q}\]
is convergent, we conclude that the sequence of holomorphic functions 
\[(\frac{1}{q^d}\sum_{j=0}^{N}l([\frac{M}{\bpq{J}{q}M}]_j)e^{-iyj/q})_N\]
converges uniformly to $\frac{1}{q^d}G_{n,M,J}(y)$ on $K$. This proves the holomorphicity of $\frac{1}{q^d}G_{n,M,J}$ on $\Omega$.

\noindent (2) When $d=0$, the conclusion follows from a direct computation. Assume $d \geq 1$. Since
\[l([\frac{M}{\bpq{J}{q}M}]_j)= l([\frac{M}{\bpq{J}{q}M}]_{\leq j})- l([\frac{M}{\bpq{J}{q}M}]_{\leq {j-1}}),\]

\noindent a direct computation using the equation above shows that,
\begin{equation}\label{eq: integral form of Frob-Poincare using the truncation}
\begin{split}
  \sum \limits_{j=0}^{\infty}l([\frac{M}{\bpq{J}{q}M}]_j)e^{-iyj/p^n} &= \sum \limits_{j=0}^{\infty}l([\frac{M}{\bpq{J}{q}M}]_ {\leq j})e^{-iyj/p^n}(1-e^{-iy/p^n}) \ .
\end{split} 
\end{equation}

\noindent Since
\[l(\frac{(\mathfrak m^j+ \bpq{J}{q})M}{(\mathfrak m^{j+1}+ \bpq{J}{q})M})= l([\frac{M}{(\mathfrak m^{j+1}+ \bpq{J}{q})M}])- l([\frac{M}{(\mathfrak m^{j}+ \bpq{J}{q})M}]),\]
a direct computation shows that, 

\begin{equation}\label{eq: integral representation of adic Frob-Poincare}
    \sum \limits_{j=0}^{\infty}l(\frac{(\mathfrak m^j+ \bpq{J}{q})M}{(\mathfrak m^{j+1}+ \bpq{J}{q})M})e^{-iyj/p^n}= \sum \limits_{j=0}^{\infty}l(\frac{M}{(\mathfrak m^{j+1}+\bpq{J}{q})M})e^{-iyj/p^n}(1-e^{-iy/p^n})
\end{equation}

\noindent Choose $a$ such that as an $R$-module $M$ is generated by homogeneous elements of  degree at most $a$. Therefore

\[\mathfrak m^jM \subseteq M_{\geq j} \subseteq \mathfrak m^{j-a}M.\]

\noindent So,

\begin{equation*}
\begin{split}
   l(\frac{M}{(\mathfrak m^{j+1}+\bpq{J}{q})M})- l([\frac{M}{\bpq{J}{q}M}]_ {\leq j})&= l(\frac{M_{\geq {j+1}}+ \bpq{J}{q}M}{\mathfrak m^{j+1}M+ \bpq{J}{q}M})\\
   &\leq l(\frac{\mathfrak m^{j+1-a}M+ \bpq{J}{q}M}{\mathfrak m^{j+1}M+ \bpq{J}{q}M})\\
   & \leq l(\frac{\mathfrak m^{j+1-a}M}{\mathfrak m^{j+1}M})\\
   & \leq Cj^{d-1},
\end{split}
\end{equation*}

\noindent for some $C$, which is independent of $q$ and $j$. Using \Cref{eq: integral form of Frob-Poincare using the truncation}, \Cref{eq: integral representation of adic Frob-Poincare} and the comparison above, we get that for any $y \in \Omega$,

\begin{equation*}
\begin{split}
    |\frac{1}{q^d}G_{n,M,J}(y)-\frac{1}{q^d}F_{n, \mathfrak m, J}(y)| &\leq \sum \limits_{j=0}^{\infty} C\frac{1}{q}(\frac{j}{q})^{d-1}e^{-\Im y j/q}|1-e^{-iy/q}|\\
    & = C|1-e^{-iy/q}|\int \limits_{0}^{\infty}\flor{s}^{d-1}e^{-\Im y \flor{s}}ds \\
    & \leq C  |1-e^{-iy/q}|\int \limits_{0}^{\infty}s^{d-1}e^{-\Im y (s-1)}ds.
\end{split}
\end{equation*}

\noindent Since $\Im y<0$ for $y \in \Omega$, the last integral is convergent. It follows from the last chain of inequalities that on a compact subset of $\Omega$, 
\[ |\frac{1}{q^d}G_{n,M,J}(y)-\frac{1}{q^d}F_{n, \mathfrak m, J}(y)|\]

\noindent uniformly converges to zero. This finishes the proof of (2).
\end{proof}

\section{$h$-function and density function near boundaries} In this section, we discuss the behaviour of the $h$-function $h(s)$ near zero and $s$ large enough. The regions near zero and away from zero where the $h$-function often shows interesting behaviour are marked by two other already known invariants, namely $F$-limbus and $F$-threshold. Recall that $F$-threshold is a well-known numerical invariant in characteristic $p$ which compares the ordinary power and Frobenius power of an ideal pair. It was initially defined as a limsup in \cite{huneke2008f} and \cite{MTW}, and is shown to be a limit in \cite{de2018existence}. The $F$-limbus is defined in \cite{TaylorInterpolation}. Recall,

\begin{definition}\label{de: definition of F-limbus and F-threshold}
Let $R$ be a ring of characteristic $p > 0$ which is not necessarily local, and let $I, J$ be ideals of $R$. Define
$$c^J_I(n)=\sup\{t \in \mathbb{N}: I^t \nsubseteq J^{[p^n]}\},$$
$$c^J(I)=\lim_{n \to \infty} \frac{c^J_I(n)}{p^n},$$
$$b^J_I(n)=\inf\{t \in \mathbb{N}: J^{[p^n]} \nsubseteq I^t\},$$
$$b^J(I)=\lim_{n \to \infty} \frac{b^J_I(n)}{p^n}.$$
\end{definition}
The number $c^J(I)$ is called the $F$-threshold of $I$ with respect to $J$ and the number $b^J(I)$ is called the $F$-limbus of $I$ with respect to $J$. The following properties are proven in \cite[Lemma 3.2]{TaylorInterpolation}.

\begin{lemma}
Let $(R, \mathfrak{m})$ be a local ring of characteristic $p > 0$, and let $I, J$ be proper ideals of $R$.

\begin{enumerate}
\item For any $I,J$, any limit in \Cref{de: definition of F-limbus and F-threshold} either exists or goes to infinity.

\item If $I$ is not nilpotent, then $b^J(I) \leq c^J(I)$.

\item If $I \nsubseteq \sqrt{J}$ then $c^J(I)=\infty$.
\item If $I \subset \sqrt{J}$ then $0\leq c^J(I)<\infty$.
\item If $J \nsubseteq \sqrt{I}$ then $b^J(I)=0$.
\item If $J \subset \sqrt{I}$ then $0<b^J(I)\leq \infty$.
\item If  $I$ is not nilpotent, $I \subset \sqrt{J}$, $J \subset \sqrt{I}$, then $0<b^J(I)\leq c^J(I)<\infty$.
\end{enumerate}    
\end{lemma}
\begin{lemma}\label{le: h function is stable near boundaries}
Let $(R,\mathfrak{m})$ be a local ring of dimension $d$ and characteristic $p$, let $I,J$ be two proper ideals of $R$, and let $M$ be a finitely generated $R$-module.
\begin{enumerate}
\item If $I$ is $\mathfrak{m}$-primary, then $b^J(I)>0$ and for $s \leq b^J(I)$, $h_M(s)=\frac{s^d}{d!}e(I,M)$.
\item If $J$ is $\mathfrak{m}$-primary, then $c^J(I)<\infty$ and for $s \geq c^J(I)$, $h_M(s)=e_{HK}(J,M)$.
\end{enumerate}
\end{lemma}
\begin{proof}
The above lemma is a generalization of Lemma 3.3 of \cite{TaylorInterpolation}. The proof is identically the same since it only uses the containment relation and does not use the $\mathfrak{m}$-primariness of $I,J$. If $I$ is $\mathfrak{m}$-primary then $J \subset \sqrt{I}$, so $b^J(I)>0$; if $J$ is $\mathfrak{m}$-primary then $I \subset \sqrt{J}$, so $c^J(I)<\infty$.   
\end{proof}

\subsection{Tail: $F$-threshold, minimal stable point and maximal support}\label{sse: tail of h function}

Let $(R,\mathfrak{m})$ be a local ring of characteristic $p>0$, $I,J$ are proper $R$-ideals. Assume $J$ is $\mathfrak{m}$-primary and $I \subseteq \mathfrak{m}$. By \Cref{le: h function is stable near boundaries}, (2), when $J$ is $\mathfrak{m}$-primary, the $h_{M,I,J}(s)$-becomes the constant $e_{HK}(J,M)$ for large enough $s$. Since $h$ is increasing and $h_M(s) \leq e_{HK}(J,M)$ for any $s$, there is a smallest point after which $h_{M,I,J}(s)$ becomes a constant. We relate this smallest point  to another seemingly unrelated invariant of $(I,J)$ which we call the $F$-threshold upto tight closure; see \Cref{de: F threshold upto tight closure}. The next lemma guarantees the existence of the invariant named the $F$-threshold upto tight closure in \Cref{de: F threshold upto tight closure}. 

\begin{lemma}\label{le: F threshold upto tight closure}
   Let $(R,\mathfrak{m},k)$ be a local ring of characteristic $p>0$, $I,J$ be two $R$-ideals, $I \subset \sqrt{J}$. Let
$$r^J_I(n)=\max\{t \in \mathbb{N}|\, I^t \nsubseteq (J^{[p^n]})^*\},$$
Then $(r^J_I(n)/p^n)_n$ is a increasing sequence converging to a real number.
\end{lemma}

\begin{proof}
Given a natural number $n$, pick $x \in I^{r^J_I(n)}\setminus (J^{[q]})^*$. Note that $x^p$ cannot be in $(J^{[pq]})^*$. Indeed, in contrary say $x^p \in (J^{[pq]})^*$. Then there is a $c \in R$ not in any minimal prime of $R$ such that $cx^{p^{m+1}} \in (J^{[q]})^{[p^{m+1}]}$ for any large $m$. This implies $x \in (J^{[q]})^*$. So we conclude
$$r^J_I(n+1) \geq pr^J_I(n),$$
whence the desired increasingness follows. By \Cref{l: comparison of ususal and Frobenius powers}, $(r^J_I(n)/p^n)_n$ is bounded and hence converges to a real number.
\end{proof}

\begin{definition}\label{de: F threshold upto tight closure} Let $(r^J_I(n))_n$ be the same as in \Cref{le: F threshold upto tight closure}. The limit of $(r^J_I(n)/p^n)_n$, which exists by \Cref{le: F threshold upto tight closure} is called the \textit{$F$-threshold up to tight closure} for the ideal pair $I,J$ and is denoted by $r_{R,I,J}$.
\end{definition}

\begin{lemma}\label{le: tight closre of Frobenius power defines Hilbert-Kunz} Let $(R,\mathfrak{m},k)$ be a $d$-dimensional reduced local ring of characteristic $p>0$, $J$ be an $\mathfrak{m}$-primary $R$-ideal. Then $e_{HK}(J,R)=\lim_{n \to \infty}l(R/(J^{[q]})^*)/q^d$.
\end{lemma}

\begin{proof}It suffices to show $\lim_{n \to \infty}l((J^{[q]})^*/J^{[q]})/q^d=0$. There is a test element $c \in R$, which is in not contained in any minimal primes of $R$ such that, $c(\bpq{J}{q})^* \subseteq \bpq{J}{q}$ for all $q$; see \cite{HunekeTight}. So we have $l((J^{[q]})^*/J^{[q]}) \leq l(0:_{R/J^{[q]}}c)=l(R/cR+J^{[q]}) \leq Cq^{d-1}$ for some constant $C$, so $$\underset{n \to \infty}{\lim}l((J^{[q]})^*/J^{[q]})/q^d=0.$$
\end{proof}

We analyze how the $h$-function changes if we replace an ideal by its usual power or Fobenius power. This is used in the sequel.

\begin{proposition}\label{th: h-function of powers of ideals}
Let $n_0 \in \mathbb{N}$, then $$h_{M,I^{n_0},J}(s)=h_{M,I,J}(sn_0), h_{M,I,J^{[p^{n_0}]}}(s)=p^{n_0d}h_{M,I,J}(s/p^{n_0}).$$
\end{proposition}

\begin{proof}If $s \leq 0$ then both sides of the equations are 0 and the equality holds. Now we assume $s>0$. By definition $h_{n,M,I^{n_0},J}(s)=l(M/I^{n_0\ceil{sq}}+J^{[q]}M)$. Since $\ceil{sqn_0} \leq n_0\ceil{sq} \leq \ceil{sqn_0}+n_0-1$, $h_{n,M,I,J}(sn_0) \leq h_{n,M,I^{n_0},J}(s) \leq h_{n,M,I,J}(sn_0+(n_0-1)/q)$. We have $\lim_{n \to \infty}(h_{n,M,I,J}(sn_0+(n_0-1)/q)-h_{n,M,I,J}(sn_0))/q^d=0$ by \Cref{th: Lipschitz continuity for a family}. So 
$$\lim_{n \to \infty}h_{n,M,I^{n_0},J}(s)/q^d=\lim_{n \to \infty}h_{n,M,I,J}(sn_0)/q^d,$$ 
which means $h_{M,I^{n_0},J}(s)=h_{M,I,J}(sn_0)$. We have $h_{n,M,I,J^{[p^{n_0}]}}(s)=l(M/I^{\ceil{sq}}+J^{[qp^n_0]}M)=l(M/I^{\ceil{s/p^{n_0}\cdot qp^{n_0}}}+J^{[qp^n_0]}M)$. So
\begin{align*}
\lim_{n \to \infty}\frac{h_{n,M,I,J^{[p^{n_0}]}}(s)}{q^d}\\
=p^{n_0d}\lim_{n \to \infty}\frac{h_{n+n_0,M,I,J}(s/p^{n_0})}{q^dp^{n_0d}}\\
=p^{n_0d}h_{M,I,J}(s/p^{n_0}).
\end{align*}
\end{proof}

\begin{theorem}\label{th: maximal stable point of h function for J}
Let $(R,\mathfrak{m},k)$ be a noetherian local ring of characteristic $p>0$, $I, J$ be proper $R$-ideals and $J$ be moreover $\mathfrak{m}$-primary, $M$ be a finitely generated $R$-module. Define $$\alpha_{M,I,J}=\sup\{s|\, 
 h_{M,I,J}(s) \neq e_{HK}(J,M)\}=\min\{s|\, h_{M,I,J}(s) = e_{HK}(J,M)\}.$$
Then
$$\alpha_{R,I,J}= r_{R,I,J}.$$
\end{theorem}

\begin{proof}
For simplicity, first assume $R$ is a complete local domain. It suffices to prove: 
\begin{enumerate}
\item For $x \in \mathbb{Z}[1/p]$, if $x > r_{R,I,J}$, then $x \geq \alpha_{R,I,J}$;
\item For $x \in \mathbb{Z}[1/p]$, if $x < r_{R,I,J}$, then $x \leq \alpha_{R,I,J}$.
\end{enumerate}

\noindent (1): If $x > r_{R,I,J}$, then there is an infinite sequence $n_i$, such that $xp^{n_i}>r^J_I(n_i)$ and $xp^{n_i}$ is an integer for all $i$. By definition of $r_n$, $I^{xp^{n_i}} \subset  (J^{[p^{n_i}]})^*$. So 
$$h_{R,I,J}(x)=\lim_{i \to \infty}l(R/I^{\lceil xp^{n_i} \rceil}+(J^{[p^{n_i}]})^*)/q^d=\lim_{i \to \infty}l(R/(J^{[p^{n_i}]})^*)/q^d=e_{HK}(J,R).$$ 
The last equality in the above chain follows from \Cref{le: tight closre of Frobenius power defines Hilbert-Kunz}. So $x \geq \alpha_{R,I,J}$.\\\\
(2): If $x < r_{R,I,J}$, then there is an integer $n_0$, such that $xp^{n_0} \leq r^J_I(n_0)$ and $xp^{n_0}$ is an integer. Let $q_0=p^{n_0}$. By definition of $r^J_I(n) $, $I^{xq_0} \nsubseteq  (J^{[q_0]})^*$.  Choose $f \in I^{xq_0}\backslash (J^{[q_0]})^*$. Let $\tilde{J}=J^{[q_0]}+fR$; then $e_{HK}(\tilde{J},R)<e_{HK}(J^{[q_0]},R)$; see \cite[Theorem 5.5]{HunekeExp}, \cite[Theorem 8.17]{HH}. Now fix an $s<xq_0$, then for any $q=p^n$, $sq<xqq_0$. Since $f \in I^{xq_0}$, $f^q \in I^{xqq_0} \subseteq I^{\ceil{sq}}$. Therefore, $$I^{\ceil{sq}}+(J^{[q_0]}+fR)^{[q]}=I^{\ceil{sq}}+(J^{[q_0]})^{[q]}.$$ This means $h_{R,I,\tilde{J}}(s)=h_{R,I,J^{[q_0]}}(s)$. So for $s<xq_0$, $h_{R,I,J^{[q_0]}}(s)=h_{R,I,\tilde{J}}(s)\leq e_{HK}(\tilde{J},R)<e_{HK}(J^{[q_0]},R)$. This means $\alpha_{R,I,J^{[q_0]}} \geq xq_0$. By \Cref{th: h-function of powers of ideals}, $h_{R,I,J^{[q_0]}}(s)=q_0^dh_{R,I,J}(s/q_0)$, $\alpha_{R,I,J}=\frac{\alpha_{R,I,J^{[q_0]}}}{q_0} \geq x$.

Now we argue that without loss of generality $R$ can be taken to be a complete domain. Note,
\begin{equation}\label{eq: alpha and r can be computed mod minimal primes}
    \alpha_{R,I,J}= \underset{Q\, \textup{minimal prime of}\,R}{\textup{max}} \{\alpha_{\frac{R}{Q},I\frac{R}{Q}, J\frac{R}{Q}}\}, \,\,\,\, r_{R,I,J}= \underset{Q\, \textup{minimal prime of}\,R}{\textup{max}} \{r_{\frac{R}{Q},I\frac{R}{Q}, J\frac{R}{Q}}\}.
\end{equation}
The above description of $\alpha_{R,I,J}$ follows from \Cref{pr: associativity formula for h function}. The above description of $r_{R,I,J}$ follows from \cite[Thm 1.3]{HunekeTight}. Thanks to \Cref{eq: alpha and r can be computed mod minimal primes}, it suffices to prove the present theorem when $R$ is a domain. Assume $R$ is a domain. Since $J$ is $\mathfrak m$-primary, $r_{R,I,J}$ coincides with $r_{\hat R, I\hat R, J\hat R}.$ Indeed as $\bpq {J}{q}$ is $\mathfrak m$-primary for all $q$ and $R$ is a domain, $(\bpq {J}{q}\hat R)^*= (\bpq {J}{q})^*\hat R$; see \cite[Thm 7.16, (a)]{HHBaseChange}. On the other hand, $h_{R,I,J}= h_{\hat R, I\hat R,J\hat R}$; so  $\alpha_{R,I,J}= \alpha_{\hat R, I\hat R,J\hat R}$. So without loss of generality $R$ can be taken to be $\mathfrak m$-adically complete. We can pass to the complete domain case using \Cref{eq: alpha and r can be computed mod minimal primes}.
\end{proof}

Next we show that $[0, \alpha_{R,I,J}]$ is the support of the density function $f_{R,I,J}$ and the $h$-function is strictly increasing until it becomes a constant function. Precisely, we have:

\begin{proposition}\label{c: minimal stable point is the support}
Let $(R,\mathfrak{m},k)$ be a local ring of characteristic $p>0$, $I,J$ be proper $R$-ideals such that $I+J$ is $\mathfrak{m}$-primary. 
\begin{enumerate}
    \item When $J$ is additionally $\mathfrak{m}$-primary,  $\alpha_{R,I,J}=\sup\{s|\, f_{R,I,J}(s) \, \text{exists and is nonzero}\}$. Moreover for $s> \alpha_{R,I,J}$, $f_{R,I,J}(s)$ is zero, and for $0<s< \alpha_{R,I,J}$, it is nonzero whenever it is well-defined.

    \item Let $s_0= \sup\{s \in \mathbb{R}_{>0} \, |\,
    h_{R,I,J}\, \text{is strictly increasing on}\, (0, s] \}$; the supremum is taken to be zero when the set is empty. When $s_0$ is a real number, $h_{R,I,J}$ is constant on $(s_0, \infty)$.

    \item When $J$ is additionally $\mathfrak{m}$-primary, the quantity $s_0$ defined in (2) coincides with $\alpha_{R,I,J}$. So $h_{R,I,J}$ is strictly increasing on $(0,\alpha_{R,I,J}]$.
\end{enumerate}

\end{proposition}

\begin{proof} (1) For $s> \alpha_{R,I,J}$, $h_{I,J}(s)$ is constant. So by \Cref{th: differentiability of h function implies existence of density function}, $f_{I,J}$ exists and is zero. Since $h_{I,J}$ is the integral of the density function by \Cref{th: h is the integral of f} and $h$ is a non-constant increasing function on $(a, \alpha_{R,I,J})$ for any $0<a< \alpha_{R,I,J}$, $f_{I,J}$ has to be nonzero on a subset of $(a,\alpha_{R,I,J})$ of nonzero measure. So when $f_{I,J}(a)$ exists, it is nonzero by \Cref{pr: monotonicity result for the desnity function}.\\

\noindent (2) Assume $s_0 \in \mathbb{R}$. Suppose by contradiction, $h_{R,I,J}$ is not constant on $(s_0, \infty)$. Then there exists $s_2>s_1>s_0$ such that $h_{I,J}(s_2)>h_{I,J}(s_1)$. Since by \Cref{th: h is the integral of f}, 
$$h_{R,I,J}(s_2)-h_{R,I,J}(s_1)= \int_{s_1}^{s_2}f_{R,I,J}(t)dt,$$
the integral above is positive. So the function $f_{R,I,J}$ must be defined and positive on a subset of $(s_1, s_2)$ with positive measure. So by \Cref{pr: monotonicity result for the desnity function}, (2), on $(0, s_1)$ wherever $f_{R,I,J}$ exists, it takes a positive value. Since by \Cref{th: differentiability of h function implies existence of density function}, (3),  $f_{R,I,J}$ exists outside a countable subset of $(0, s_1)$, $h_{I,J}$ is strictly increasing on $(0, s_1]$ by \Cref{th: h is the integral of f}. This contradicts the supremality of $s_0$.

For (3), it is enough to prove that the $s_0$ in (2) coincides with $\alpha_{R,I,J}$. To that end, note, $h_{R,I,J}$ must be the constant function $e_{HK}(J,R)$ on $(s_0, \infty)$. So $\alpha_{R,I,J} \leq s_0$. Since $h_{R,I,J}$ is strictly increasing on $(0, s_0)$, $s_0$ cannot be strictly larger than $\alpha_{R,I,J}$. 
\end{proof} 

Moreover, in \Cref{co: h is strictly increasing when J is not m primary} we prove that $h_{I,J}$ eventually becomes constant if and only if $J$ is $\mathfrak{m}$-primary. In other words, the $s_0$ defined in \Cref{c: minimal stable point is the support} is finite if and only if $J$ is $\mathfrak{m}$-primary. These statements follow from our next result, where we compute the limit of $h_{I,J}(s)/s^{\dim(R/J)}$ as $s$ approaches infinity. 

\begin{theorem}\label{th: behaviour of h at infinity}
    Let $R$ be a domain, $I,J$ be ideals such that $I+J$ is $\mathfrak{m}$-primary. Then
    \[\underset{s \to \infty}{\lim}\frac{h_{R,I,J}(s)}{s^{\dim(R/J)}}= \frac{\underset{Q \in \textup{Assh}(R/J)}{\sum}e(I,R/Q)e_{HK}(JR_Q, R_Q)}{\dim(R/J)!},\]
  where $\textup{Assh}(R/J)$ is the set of associated primes of $R/J$ of dimension $\dim(R/J)$.  
\end{theorem}

\begin{proof}
Without loss of generality, we assume that the residue field of $R$ is perfect. Set $d= \dim(R), d'= \dim(R/J)$, $h_n(s)= l(R/(I^{\ceil{sq}}+\bpq{J}{q})R)$. We first aim to show that \\
\textbf{Claim:} For $s_0>0$, the sequence of functions $(h_n(s)/s^{d'}q^d)_n$ is uniformly Cauchy on $(s_0, \infty)$.
The key ingredient in the proof of the claim presented further below is:
\begin{lemma}\label{le: convergence estimate of h}
    There is an $R$ module $N$ of dimension at most $d-1$, where $N$ depends only on $R$ and $I$, such that
    \[|\frac{h_{n+1}(s)}{p^{(n+1)d}}- \frac{h_{n}(s)}{p^{nd}}| \leq \frac{1}{q^d}l(\frac{N}{(I^{\ceil{sq}}+\bpq{J}{q})N})\]
\end{lemma}
\begin{proof}
Note that 
\[\frac{h_{n+1}(s)}{p^{(n+1)d}}- \frac{h_{n}(s)}{p^{nd}} \leq \frac{1}{(qp)^d}l(\frac{R}{I^{\ceil{sq}p}+ \bpq{J}{qp}})-\frac{h_n(s)}{q^d} \leq \frac{1}{(qp)^d}l(\frac{R}{I^{\ceil{sq}[p]}+ \bpq{J}{qp}})-\frac{h_n(s)}{q^d} .\]
Choose an $R$-module $N_1$ of dimension at most $d-1$, fitting into an exact sequence
$$0 \rightarrow R^{\oplus p^d} \rightarrow F_*R \rightarrow N_1 \rightarrow 0.$$
Tensoring this exact sequence with $R/(I^{\ceil{sq}}+\bpq{J}{q})$, we obtain a right exact sequence
\[ (\frac{R}{I^{\ceil{sq}}+\bpq{J}{q}})^{\oplus p^d} \rightarrow F_*(\frac{R}{I^{\ceil{sq}[p]}+ \bpq{J}{qp}}) \rightarrow \frac{N_1}{(I^{\ceil{sq}}+\bpq{J}{q})N_1} \rightarrow 0.\]
Thus we conclude,
\begin{equation} \label{eq: inequality of the difference}
    \frac{h_{n+1}(s)}{p^{(n+1)d}}- \frac{h_{n}(s)}{p^{nd}} \leq \frac{1}{q^d}l(\frac{R}{I^{\ceil{sq}[p]}+ \bpq{J}{qp}})-\frac{h_n(s)}{q^d} \leq \frac{1}{(qp)^d} l(\frac{N_1}{(I^{\ceil{sq}}+\bpq{J}{q})N_1}) \ .
\end{equation}
Next choose an $R$-linear injection
$F_*R \xrightarrow{\phi} R^{\oplus p^d}$
where $\dim \textup{Coker}(\phi)<\dim R$. Let $\mu$ be the minimal number of generators of $I$. Choose a nonzero $c \in I$ and let $\psi=c^\mu\phi$. Note that 
\begin{equation*}
    \begin{split}
        \psi(F_*(I^{\ceil{spq}}+ \bpq{J}{pq})) \subseteq \psi(F_*(I^{\ceil{sq}p}+ \bpq{J}{pq})) & \subseteq \psi(F_*(I^{(\ceil{sq}-\mu)[p]}+ \bpq{J}{pq}))\\
        & \subseteq c^\mu(I^{\ceil{sq}-\mu}+ \bpq{J}{q}) \subseteq I^{\ceil{sq}}+ \bpq{J}{q}.
    \end{split}
\end{equation*}
The second containment above follows from \Cref{l: comparison of ususal and Frobenius powers}. Thus $\psi$ induces a map $$F_*(R/I^{\ceil{spq}}+ \bpq{J}{pq}) \rightarrow R/(I^{\ceil{sq}}+ \bpq{J}{q}),$$ whose cokernel is a quotient of $\text{coker}(\phi)/(I^{\ceil{sq}}+ \bpq{J}{q})\text{coker}(\phi)$. This implies

\begin{equation}\label{eq: the other inequality}
   \frac{h_{n}(s)}{p^{nd}}-  \frac{h_{n+1}(s)}{p^{(n+1)d}} \leq \frac{1}{(qp)^d}l(\frac{\text{coker}(\phi)}{(I^{\ceil{sq}}+ \bpq{J}{q})\text{coker}(\phi)}) \ . 
\end{equation}
From \Cref{eq: inequality of the difference}, \Cref{eq: the other inequality}, it follows that taking $N= N_1 \oplus \text{coker}(\phi)$ suffices for the desired statement.
\end{proof}
\begin{proof}[Proof of \textbf{Claim}] Let $N$ be as in \Cref{le: convergence estimate of h}. Note,
\begin{equation*}
    \begin{split}
        \frac{1}{q^d}l(\frac{N}{(I^{\ceil{sq}}+\bpq{J}{q})N}) &\leq \frac{1}{q^{d}}l(\frac{N}{(I^{\ceil{s}q}+\bpq{J}{q})N}) \leq \frac{1}{q^d}l(\frac{N}{(I^{\ceil{s}[q]}+\bpq{J}{q})N})\\
        & \leq \frac{1}{q^d}l(\frac{F^n_*N}{(I^{\ceil{s}}+J)F^n_*N}) \leq \frac{1}{q}\frac{\mu(F^n_*N)}{q^{d-1}}l(\frac{R}{I^{\ceil{s}}+J}) \ ,
    \end{split}
\end{equation*}
where $\mu(-)$ is the minimal number of the generators of the corresponding $R$-module. Since $\dim(N)$ is at most $d-1$, we can choose $C_1$ such that $\frac{\mu(F^n_*N)}{q^{d-1}}\leq C_1$ for all $n>0$. Choose $C_2$ such that, $l(\frac{R}{I^{\ceil{s}}+J}) \leq C_2(\ceil{s})^{d'}$ for all $s \in \mathbb{R}$. Using the above chain inequalities with \Cref{le: convergence estimate of h}, conclude
\[|\frac{h_{n+1}(s)}{s^{d'}p^{(n+1)d}}- \frac{h_{n}(s)}{s^{d'}p^{nd}}| \leq \frac{C_1C_2\ceil{s}^{d'}}{qs^{d'}}\leq \frac{C_1C_2}{q}(1+\frac{1}{s})^{d'}.\]
Since $1/s$ is bounded on $(s_0, \infty)$, the claim is established.
\end{proof}
\noindent \textbf{Proof of \Cref{th: behaviour of h at infinity} continued:} The \textbf{Claim} above asserting uniform Cauchyness implies that we can commute the order of the limits in $\lim_{s \to \infty}\lim_{n \to \infty}h_n(s)/s^{d'}q^d$: this means, $\lim_{s \to \infty}h_{I,J}(s)/s^{d'}$ exists if and only if $\lim_{n \to \infty}\lim_{s \to \infty}h_n(s)/s^{d'}q^d$ exists and the limits coincide once existence is proven. We finish the proof by noting,
\begin{equation*}
\begin{split}
   \lim_{n \to \infty}\lim_{s \to \infty}h_n(s)/s^{d'}q^d &= \frac{1}{d'!}\lim_{n \to \infty} \frac{1}{q^{d-d'}}e(I, \frac{R}{\bpq{J}{q}}) \\
   & = \frac{1}{d'!}\lim_{n \to \infty} \frac{1}{q^{d-d'}}\underset{Q \in \text{Assh}(R/J)}{\sum}e(I, R/Q)l(R_Q/\bpq{J}{q}R_Q)\\
  & = \frac{1}{d'!}\underset{Q \in \text{Assh}(R/J)}{\sum}e(I, R/Q)e_{HK}(JR_Q, R_Q) \ .
\end{split}   
\end{equation*}

\end{proof}
\begin{remark}
    Using the associativity formula \Cref{pr: associativity formula for h function}, one can get an analogue of \Cref{th: behaviour of h at infinity}, for any finitely generated module over an $F$-finite ring.
\end{remark}

\begin{corollary}\label{co: h is strictly increasing when J is not m primary}
    Let $I,J$ be ideals of a domain $R$ such that $I+J$ is $\mathfrak{m}$-primary.
    \begin{enumerate}
        \item Then $h_{I,J}$ is eventually constant if and only if $J$ is $\mathfrak{m}$-primary.

        \item If $J$ is not $\mathfrak{m}$-primary, $h_{I,J}$ is a strictly increasing function on $\mathbb{R}_{>0}$.
    \end{enumerate}
\end{corollary}

\begin{proof}
  (1)  When $J$ is $\mathfrak{m}$-primary, $h_{I,J}(s)$ is $e_{HK}(J,R)$ for large values of $s$; see \Cref{le: compactly supported h-function}. Assume $J$ is not $\mathfrak{m}$-primary. So $\dim(R/J)>0$. Thus $h_{I,J}(s)$ approaches infinity as $s$ approaches infinity.

  (2) Since $J$ is not $\mathfrak{m}$-primary, $h_{I,J}(s)$ is not eventually constant. Thus the $s_0$ defined in \Cref{c: minimal stable point is the support}, (2) is not finite. Thus $h_{I,J}(s)$ must be strictly increasing on $\mathbb{R}_{>0}$.
\end{proof}
The following immediate consequence of \Cref{c: minimal stable point is the support} proves that the Hilbert-Kunz density function as defined in \cite{TriExist} is strictly positive on the interior of its support.

\begin{corollary}\label{co: HK density is strictly positive in the interior}
    Let $(R, \mathfrak{m})$ be a standard graded ring of dimension at least two, $J$ be an $\mathfrak{m}$-primary homogeneous ideal. Then the Hilbert-Kunz density function $\tilde{g}_{R,J}$ as defined in \cite{TriExist} (see \Cref{de: HK desnity after Trivedi}) is positive in the interior of the support of $\tilde{g}_{R,J}$.
\end{corollary}

\begin{proof}
    By \Cref{th: h function differentiable when J is m primary}, $f_{R, \mathfrak{m},J}$ exists at all points of $\mathbb{R}$ and coincides with $\tilde{g}_{R,J}$. Now by \Cref{c: minimal stable point is the support}, (1) the support of $\tilde{g}_{R,J}$ is $[0, \alpha_{R, \mathfrak{m},J}]$ and $\tilde{g}_{R,J}$ is positive on $(0, \alpha_{R, \mathfrak{m},J})$.
\end{proof}

\begin{remark}\label{re: maximal support of h and relation to Trivedi;s result}
Recall from \Cref{th: HK density when J is homogeneous but not of finite colength} that for standard graded $(R,\mathfrak m)$ of Krull dimension at least two and a finite colength homogeneous ideal $J$, Trivedi's density function $\tilde{g}_{R,J}$ coincides with $f_{R,\mathfrak m,J}$ and both are continuous. So \Cref{th: maximal stable point of h function for J} gives a precise description of the support of $\tilde{g}_{R,J}$. Thus \Cref{th: maximal stable point of h function for J} and the theorem below extends \cite[Theorem 4.9]{TrivediWatanabeDomain}, where $\alpha_{R,J}$ is shown to coincide with the $F$-threshold $c^J(\mathfrak m)$ under suitable hypothesis. 
\end{remark}

\begin{remark}\label{re: Trivedi's question regarding the support}
    When $(R, \mathfrak{m})$ is standard graded of dimension at least two and $J$ is an $\mathfrak{m}$-primary homogeneous ideal, Trivedi-Watanabe asks if $\alpha_{R,\mathfrak{m},J}= c^J(\mathfrak{m})$\footnote{note the change in notation.}; see \cite[pn 537]{TriFthreshold}. In this context, since Trivedi's density function coincides with $f_{R, \mathfrak{m},J}$, \Cref{c: minimal stable point is the support} translates Trivedi- Watanabe's question to the question asking if the $F$-threshold and the $F$-threshold upto tight closure coincide for the pair $\mathfrak{m},J$; see \Cref{th: maximal stable point of h function for J}, \Cref{c: minimal stable point is the support}.
\end{remark}

Next we compare $F$-threshold and the $F$-threshold up to tight closure. The following immediate observation is used later in \Cref{se: applications}; see \Cref{co: lower bound on F threshold} for example.

\begin{lemma}\label{le: F threshold upto tight closure is larger}
    Let $I, J $ be two ideals in a prime characteristics $F$-finite noetherian ring $R$ such that $I \subseteq \sqrt{J}$. Then $r^J(I) \leq c^J(I)$.
\end{lemma}

Below, we isolate a few cases where the $F$-threshold and the $F$-threshold up to tight closure coincides. We do not know whether in general these two invariants coincide; see \Cref{qe: minnimal stable point vs F-threshold}.

\begin{theorem}\label{th: when the maximal support coincides with F-threshold}Let $(R,\mathfrak{m},k)$ be a local ring of characteristic $p>0$, $I$ be an $R$-ideal, $J$ be an $\mathfrak{m}$-primary $R$-ideal. Then $c^J(I)=r_{R,I,J}$ is true under either of the assumptions below:
\begin{enumerate}
\item There exists a sequence of positive numbers $r'_n$ such that $I^{r'_n} \subset J^{[q]}:(J^{[q]})^*$ for infinitely many $q \gg 0$ and $\lim_n r'_n/p^n \to 0$.
\item There exists a constant $n_0$ such that $I^{n_0} \subset J^{[q]}:(J^{[q]})^*$ for infinitely many $q \gg 0$.
\item  $R$ is $F$-rational\footnote{ see \cite{FedderWatanabe89}, \cite{KarenRationalSing}.}, i.e. the tight closure of every parameter ideal coincides with the ideal and $J$ is a parameter ideal.
\item $I \subset \sqrt{\tau(R)}$, where $\tau(R)=\cap_{\mathfrak a \subset R}\mathfrak a:\mathfrak a ^*$ is the test ideal of $R$. See \cite[Definition 8.22, Proposition 8.23]{HH} for details on the test ideal.
\item (Theorem 4.9, \cite{TriFthreshold})$R$ is strongly $F$-regular on the punctured spectrum.
\end{enumerate}
\end{theorem}
\begin{proof}
\begin{enumerate}
\item By definition $r^J_I(n) \leq c^J_I(n)$, and the condition implies $c^J_I(n) \leq r^J_I(n)+r_n$, so $\lim_n (c^J_I(n)-r^J_I(n))/p^n=0$ and $c^J(I)=r^J(I)$.
\item By (1) and the fact that $\lim_n n_0/n=0$.
\item If $J$ is a parameter ideal, so is $\bpq{J}{q}$. Since $R$ is $F$-rational, $J^{[q]}:(J^{[q]})^*=R$ for any $q$, so $n_0=1$ satisfies the assumption of (2). 
\item There exist an $n_0$ such that $I^{n_0} \subset \tau(R) \subset \cap_q J^{[q]}:(J^{[q]})^*$, and this $n_0$ satisfies the assumption of (2).
\item In this case $\tau(R)$ is either $\mathfrak{m}$-primary or is the unit ideal, so $I \subset \sqrt{\tau(R)}$ always holds.
\end{enumerate}
\end{proof}
\subsection{Head: Order of vanishing at 0 and Hilbert-Kunz multiplicity of quotient rings}
So far we have proven continuity of the $h$-function on $\mathbb R_{>0}$; see \Cref{th: Lipschitz continuity for a family} and \Cref{pr: cotinuity of $h$ functions}. In this section we determine when $h_{M,I,J}$ is continuous at $s=0$; see \Cref{th: continuity of h general case}. In \Cref{th: asymptotic behaviour of h function near zero}, we determine the order of vanishing of $h$-functions near the origin and show that the asymptotic behaviour of $h_{I,J}$ near the origin captures other numerical invariants of $(R,I,J)$. A major intermediate step involved in proving \Cref{th: asymptotic behaviour of h function near zero} is \Cref{th: commutation of double limit}, which boils down to proving commutation of the order of a double limit. We lay the groundwork for that.

Let $(R,\mathfrak{m},k)$ be a local of characteristic $p>0$, $I,J$ be two $R$-ideals such that $I+J$ is $\mathfrak{m}$-primary. Let $d=\dim R$, $d'=\dim R/I$. For a positive integer $s_0$, consider the sequence of real numbers:
$$\Gamma_{s_0,m,n}=\frac{l(R/I^{s_0p^n}+J^{[p^np^m]})}{p^{nd}p^{md'}s^{d-d'}_0}.$$

\begin{equation}\label{eq: as n goes to infinity}
    \lim_{n \to \infty}\Gamma_{s_0,m,n}=\frac{h_R(s_0/p^m)}{(s_0/p^m)^{d-d'}} .
\end{equation}

\begin{equation*}
\begin{split}
\lim_{m \to \infty}\Gamma_{s_0,m,n}&=\frac{e_{HK}(J^{[p^n]},R/I^{s_0p^n})}{p^{nd}s_0^{d-d'}}\\
&=\frac{e_{HK}(J,R/I^{s_0p^n})}{(s_0p^n)^{d-d'}}\\
&=\frac{1}{(s_0p^n)^{d-d'}}\sum_{P \in \Assh(R/I)}e_{HK}(J,R/P)l_{R_P}(R_P/I^{s_0p^n}R_P) \ ,
\end{split}
\end{equation*}
where $\Assh(R/I)$ is the set of associated primes of $R/I$ in $R$  of dimension $\dim(R/I)$; see \Cref{de: Assh}.

\noindent For $P \in \Assh(R/I)$, we have $\textup{ht}(P) \leq \dim R-\dim R/P=\dim R-\dim R/I=d-d'$. So
\begin{equation*}
\begin{split}
\lim_{n \to \infty}\lim_{m \to \infty}\Gamma_{s_0,m,n}
&=\lim_{n \to \infty}\frac{1}{(s_0p^n)^{d-d'}}\sum_{P \in \Assh(R/I)}e_{HK}(J,R/P)l_{R_P}(R_P/I^{s_0p^n}R_P)\\
&=\frac{1}{(d-d')!}\sum_{P \in \Assh(R/I), \textup{ht}P=d-d'}e_{HK}(J,R/P)e(I,R_P) \ .
\end{split}
\end{equation*}

\noindent When $R$ is an $F$-finite domain and hence an excellent domain (see \cite{Kunz76}), for all $P \in \Assh(R/I)$, $\textup{ht}(P)=d-d'$. So the above quantity is $$\frac{1}{(d-d')!}\sum_{P \in \Assh(R/I)}e_{HK}(J,R/P)e(IR_P, R_P).$$ When $R$ is a Cohen-Macaulay domain and $I$ is part of a system of parameters, the above quantity recovers the Hilbert-Kunz multiplicity $e_{HK}(J, R/I)$ as,
\begin{equation*}
\begin{split}
&\sum_{P \in \Assh(R/I)}e_{HK}(J,R/P)e(IR_P, R_P)\\
&= \sum_{P \in \Assh(R/I)}e_{HK}(J,R/P)l(R_P/IR_P)\\&= e_{HK}(J,R/I) \ .
\end{split}
\end{equation*}
\begin{theorem}\label{th: commutation of double limit}
Assume $R$ is a domain and $I \neq 0$ and $J$ be such that $I+J$ is $\mathfrak m$-primary. Fix a positive integer $s_0$. Set $\dim(R/I)=d'$. Then 
\begin{equation*}
    \lim_{m \to \infty}\frac{h(s_0/p^m)}{(s_0/p^m)^{d-d'}}= \frac{1}{(d-d')!}\sum_{P \in \Assh(R/I)} e_{HK}(J,R/P)e(I,R_P) \ ,
\end{equation*}
where $\Assh(R/I)$ is the set of associated primes of $R/I$ in $R$  of dimension $\dim(R/I)$.
\end{theorem}

\begin{proof}
We use the notation set above in this subsection. It follows from \Cref{eq: as n goes to infinity} and above that we need to show
\begin{equation*}\label{eq: commutation of double limit}
\lim_{m \to \infty}\lim_{n \to \infty}\Gamma_{s_0,m,n}=\lim_{n \to \infty}\lim_{m \to \infty}\Gamma_{s_0,m,n}.
\end{equation*}

We already see that $\lim_{n \to \infty}\Gamma_{s_0,m,n}$ and $\lim_{n \to \infty}\lim_{m \to \infty}\Gamma_{s_0,m,n}$ exist. We claim that the sequence $n \to \Gamma_{s_0,m,n}$ is uniformly convergent in terms of $m$; then, by argument of analysis, we get $\lim_{m \to \infty}\lim_{n \to \infty}\Gamma_{s_0,m,n}$ exists, and is equal to $\lim_{n \to \infty}\lim_{m \to \infty}\Gamma_{s_0,m,n}$.

To this end, we prove that there exist a constant $C$  such that $|\Gamma_{s_0,m,n+1}-\Gamma_{s_0,m,n}| \leq C/p^n$ for all $m$, which implies that $|\lim_{n \to \infty}\Gamma_{s_0,m,n}-\Gamma_{s_0,m,n}| \leq 2C/p^n$ for all $m$. We can prove it in two steps: we first prove there is a constant $C_1$ such that $\Gamma_{s_0,m,n+1}-\Gamma_{s_0,m,n} \leq C_1/p^n$, then we prove there is a constant $C_2$ such that $\Gamma_{s_0,m,n}-\Gamma_{s_0,m,n+1} \leq C_2/p^n$, then $C=\max\{|C_1|,|C_2|\}$ satisfies the statement of the claim. Without loss of generality we assume $R/\mathfrak m$ is a perfect field; see \Cref{re: changing the residue field}.

Choice of $C_1$: since $\dim R=d$, there is an exact sequence
$$0 \to R^{\oplus p^d} \to F_*R \to N \to 0$$
where $N$ is an $R$-module with $\dim N<d$. Then we have
$$(R/I^{s_0p^n}+J^{[p^np^m]})^{\oplus p^d} \to F_*R/(I^{s_0p^n}+J^{[p^np^m]})F_*R \to N/(I^{s_0p^n}+J^{[p^np^m]})N \to 0.$$
This means
\begin{equation*}
    \begin{split}
       l(\frac{R}{I^{s_0p^{n+1}}+J^{[p^{n+1}p^m]}})&\leq l(\frac{R}{I^{s_0p^n[p]}+J^{[p^{n+1}p^m]}})  \\
       & \leq p^dl(\frac{R}{I^{s_0p^n}+J^{[p^np^m]}})+l(\frac{N}{(I^{s_0p^n}+J^{[p^np^m]})N}).
    \end{split}
\end{equation*}

So dividing $p^{(n+1)d}p^{md'}s^{d-d'}_0$, we get
$$\Gamma_{s_0,m,n+1} \leq \Gamma_{s_0,m,n}+l(N/(I^{s_0p^n}+J^{[p^np^m]})N)/p^{(n+1)d}p^{md'}s^{d-d'}_0.$$
Now we claim that there is a constant $C_1>0$ that depends on $N,I,J$ and $s_0$ but is independent of $m,n$ such that $l(N/I^{s_0p^n}+J^{[p^np^m]}N)/p^{n(d-1)+d}p^{md'}s^{d-d'}_0 \leq C_1$. We have

\begin{equation*}
\begin{split}
l(N/(I^{s_0p^n}+J^{[p^np^m]})N)& \leq l(N/(I^{s_0[p^n]}+J^{[p^np^m]})N) \\
&=l(F^n_*N/(I^{s_0}+J^{[p^m]})F^n_*N)\\
&\leq \mu_R(F^n_*N)l(R/I^{s_0}+J^{[p^m]}).
\end{split}
\end{equation*}
Since $\dim N \leq d-1$ and $\dim R/I=d'$, $\mu_R(F^n_*N)/p^{n(d-1)}$ and $l(R/I^{s_0}+J^{[p^m]})/p^{md'}$ are both bounded. And $p^{-d}s^{d-d'}_0$ is independent of $m,n$. This means there is a constant $C_1>0$ that depends on $N,I,J$ and $s_0$ but is independent of $m,n$ such that $l(N/I^{s_0p^n}+J^{[p^np^m]}N)/p^{n(d-1)+d}p^{md'}s^{d-d'}_0 \leq C_1$. Thus we have
$$\Gamma_{s_0,m,n+1} \leq \Gamma_{s_0,m,n}+C_1/p^n.$$

Choice of $C_2$: since $\dim R=d$, there is an injection
$F_*R \xrightarrow{\phi} R^{\oplus p^d}$
where $\dim \textup{Coker}\phi<\dim R$. Let $\mu$ be the minimal number of generators of $I$. Choose $0 \neq c \in I$ and let $\psi=c^\mu\phi$. Since $R$ is a domain, $\psi$ is still an injection, and we have a short exact sequence
$$0 \to F_*R \xrightarrow{\psi} R^{\oplus p^d} \to N' \to 0$$
and we have $\dim N'<\dim R$. Therefore, we get an exact sequence: 
$$F_*R/(I^{s_0p^n}+J^{[p^np^m]})F_*R \xrightarrow{\bar{\psi}} (R/I^{s_0p^n}+J^{[p^np^m]})^{\oplus p^d} \to N'/(I^{s_0p^n}+J^{[p^np^m]})N' \to 0.$$
We claim that $\bar{\psi}$ induces an $R$-linear map $\Phi:F_*(R/(I^{s_0p^{n+1}}+J^{[p^{n+1}p^m]})) \xrightarrow{} (R/I^{s_0p^n}+J^{[p^np^m]})^{\oplus p^d}$. It suffices to show $\psi(F_*(I^{s_0p^{n+1}}+J^{[p^{n+1}p^m]})) \subset (I^{s_0p^n}+J^{[p^np^m]})^{\oplus p^d}$. We have $I^{s_0p^{n+1}}=I^{s_0p^np} \subset I^{(s_0p^n-\mu)[p]}$. So

\begin{align*}
\psi(F_*(I^{s_0p^{n+1}}+J^{[p^{n+1}p^m]}))\\
\subset \psi(F_*(I^{(s_0p^n-\mu)[p]}+J^{[p^{n+1}p^m]}))\\
\subset (I^{(s_0p^n-\mu)}+J^{[p^np^m]})\psi(F_*R)\\
\subset c^\mu(I^{(s_0p^n-\mu)}+J^{[p^np^m]})\phi(F_*R)\\
\subset (I^{s_0p^n}+J^{[p^np^m]})\phi(F_*R)\\
\subset (I^{s_0p^n}+J^{[p^np^m]})^{\oplus p^d}.
\end{align*}
Since $\text{coker}(\Phi)$ is a quotient of $N'/(I^{s_0p^n}+J^{[p^np^m]})N'$,
$$p^dl(R/I^{s_0p^n}+J^{[p^np^m]}) \leq l(R/I^{s_0p^{n+1}}+J^{[p^{n+1}p^m]})+l(N'/(I^{s_0p^n}+J^{[p^np^m]})N')$$
So dividing $p^{(n+1)d}p^{md'}s^{d-d'}_0$, we get
$$\Gamma_{s_0,m,n+1} \leq \Gamma_{s_0,m,n}+l(N'/(I^{s_0p^n}+J^{[p^np^m]})N')/p^{(n+1)d}p^{md'}s^{d-d'}_0$$
Since $\dim N'<\dim R$, we can use the same proof in the previous step to show that there is a constant $C_2>0$ that depends on $N',I,J$ and $s_0$ but independent of $m,n$ such that $l(N'/(I^{s_0p^n}+J^{[p^np^m]})N')/p^{n(d-1)+d}p^{md'}s^{d-d'}_0 \leq C_2$, so
$$\Gamma_{s_0,m,n} \leq \Gamma_{s_0,m,n+1}+C_2/p^n.$$
\end{proof}

\begin{theorem}\label{th: asymptotic behaviour of h function near zero}
Let $(R,\mathfrak{m},k)$ be a local domain, $I,J$ be two $R$-ideals, 
$I \neq 0$, $I+J$ is $\mathfrak{m}$-primary. Let $d=\dim R$, $d'=\dim R/I$. Then:
\begin{enumerate}
\item $\lim_{s \to 0+}h(s)/s^{d-d'}=\frac{1}{(d-d')!}\sum_{P \in \Assh(R/I)}e_{HK}(J,R/P)e(I,R_P)$.
\item The order of vanishing of $h(s)$ at $s=0$ is exactly $d-d'$.
\item $h(s)$ is continuous at $0$.
\end{enumerate}
\end{theorem}
\begin{proof}
(1) Let $\frac{1}{(d-d')!}\sum_{P \in \Assh(R/I)}e_{HK}(J,R/P)e(I,R_P)=c=c_{I,J}$, which is a constant that only depends on $I,J$. The last theorem implies for any fixed $s_0$,
$$\lim_{m \to \infty}h(s_0/p^m)/(s_0/p^m)^{d-d'}=c$$

Choose a sequence $\{s_i\}_i \subset (0,\infty)$ such that $\lim_{i \to \infty}s_i=0$ and $\lim_{i \to \infty}h(s_i)/s^{d-d'}_i$ exists. Below we argue that $\lim_{i \to \infty}h(s_i)/s^{d-d'}_i=c$; then (1) follows. Fix any $n_0 \in \mathbb{N}$. There exists an integer $\alpha_i$ for each $s_i$ such that $s_ip^{\alpha_i} \in (p^{n_0-1},p^{n_0}]$. Since $h(s)$ is an increasing function,
\begin{gather*}
\frac{h(\flor{s_iq^{\alpha_i}}/q^{\alpha_i})}{((\flor{s_iq^{\alpha_i}}+1)/q^{\alpha_i})^{d-d'}}\leq \frac{h(s_i)}{s_i^{d-d'}}
\leq \frac{h(\ceil{s_iq^{\alpha_i}}/q^{\alpha_i})}{((\ceil{s_iq^{\alpha_i}}-1)/q^{\alpha_i})^{d-d'}}\\
\implies (\frac{\flor{s_iq^{\alpha_i}}}{\flor{s_iq^{\alpha_i}}+1})^{d-d'}\frac{h(\flor{s_iq^{\alpha_i}}/q^{\alpha_i})}{(\flor{s_iq^{\alpha_i}}/q^{\alpha_i})^{d-d'}}  \leq \frac{h(s_i)}{s_i^{d-d'}}\leq(\frac{\ceil{s_iq^{\alpha_i}}}{\ceil{s_iq^{\alpha_i}}-1})^{d-d'}\frac{h(\ceil{s_iq^{\alpha_i}}/q^{\alpha_i})}{(\ceil{s_iq^{\alpha_i}}/q^{\alpha_i})^{d-d'}}\\
\implies (\frac{p^{n_0-1}}{p^{n_0-1}+1})^{d-d'}\frac{h(\flor{s_iq^{\alpha_i}}/q^{\alpha_i})}{(\flor{s_iq^{\alpha_i}}/q^{\alpha_i})^{d-d'}}  \leq \frac{h(s_i)}{s_i^{d-d'}} \leq (\frac{p^{n_0-1}}{p^{n_0-1}-1})^{d-d'}\frac{h(\ceil{s_iq^{\alpha_i}}/q^{\alpha_i})}{(\ceil{s_iq^{\alpha_i}}/q^{\alpha_i})^{d-d'}}.\\
\end{gather*}
Let $i \to \infty$, then $s_i \to 0$, $\alpha_i \to \infty$. Since $\flor{s_iq^{\alpha_i}}, \ceil{s_iq^{\alpha_i}}$ lies in $[p^{n_0-1},p^{n_0}]$, so there are only finitely many possible values of $\flor{s_iq^{\alpha_i}}, \ceil{s_iq^{\alpha_i}}$. So by \Cref{th: commutation of double limit},
$$\lim_{i \to \infty}\frac{h(\flor{s_iq^{\alpha_i}}/q^{\alpha_i})}{(\flor{s_iq^{\alpha_i}}/q^{\alpha_i})^{d-d'}}=\lim_{i \to \infty}\frac{h(\ceil{s_iq^{\alpha_i}}/q^{\alpha_i})}{(\ceil{s_iq^{\alpha_i}}/q^{\alpha_i})^{d-d'}}=c.$$
This means
$$(\frac{p^{n_0-1}}{p^{n_0-1}+1})^{d-d'} c \leq \lim_{i \to \infty}h(s_i)/s_{i}^{d-d'} \leq (\frac{p^{n_0-1}}{p^{n_0-1}-1})^{d-d'}c.$$
Since this is true for arbitrary $n_0$, we get
$$\lim_{i \to \infty}h(s_i)/s_{i}^{d-d'}=c.$$
This finishes the proof of (1).

(2) follows from (1).

(3) Since $R$ is a domain and $I \neq 0$, $d'=\dim R/I<\dim R=d$, $d-d'\geq 1$. So the order of $h(s)$ at $0$ is at least 1; in particular, $\lim_{s \to 0+}h(s)=0=h(0)$.
\end{proof}

\begin{lemma}\label{le: continuity of h at zero domain case}
Let $(R,\mathfrak{m})$ be a noetherian local domain, $I,J$ be two $R$-ideal such that $I+J$ is $\mathfrak{m}$-primary. Then $h_{R,I,J}(s)$ is continuous at $0$ if and only if $I \neq 0$.
\end{lemma}
\begin{proof}
If $I \neq 0$ then by previous theorem it is continuous at $0$. If $I=0$, then $h_R(s)=e_{HK}(J,R) \neq 0=h_R(0)$ for $s>0$, so it is discontinuous at $0$.   
\end{proof}

\begin{theorem}\label{th: continuity of h general case}
Let $(R,\mathfrak{m})$ be a noetherian local ring, $I,J$ be two $R$-ideals such that $I+J$ is $\mathfrak{m}$-primary, $M$ be a finitely generated $R$-module. Then $h_{M,I,J}(s)$ is continuous at $0$ if and only if $I \nsubseteq P$ for any $P \in \Supp(M)$ with $\dim R/P=\dim M$. In particular, $h_{R,I,J}(s)$ is continuous at $0$ if and only if $\dim R>\dim R/I$. If $h_M$ is discontinuous at $0$ then we have 
$$\lim_{s \to 0+}h_M(s)=\sum_{P\in \Supp(M), I \subset P,\dim R/P=\dim M}l_{R_P}(M_P)e_{HK}(J,R/P).$$
\end{theorem}

\begin{proof}
By the associativity formula for $h$-function in \Cref{pr: associativity formula for h function}, $$h_M(s)=\sum_{P \in \Supp(M),\dim R/P=\dim M}l_{R_P}(M_P)h_{R/P}(s).$$
For any $P \in \Supp(M)$, $\underset{s \to 0+}{\lim}h_{R/P, I,J}(s)$ is always nonnegative; the limit is positive if and only if $I\subseteq P$, in which case the limit is $e_{HK}(J, R/P)$; see \Cref{le: continuity of h at zero domain case}. Thus taking limit as $s$ approaches zero from the right, we get the expression of the right hand limit of $h_M$. Since $h_M$ is continuous at 0 if and only if $\lim_{s \to 0^+}h_{R/P}(s)=0$ for any $P \in \Supp(M)$ with $\dim R/P=\dim M$, the continuity of $h_M$ at zero is equivalent to asking $I \nsubseteq P$ for any $P \in \Supp(M)$ with $\dim R/P=\dim M$. If $M=R$, then this means $I \nsubseteq P$ for any $P \in \Assh(R)$ which means $\dim R>\dim R/I$.
\end{proof}

Next we show that, when $J$ is $\mathfrak{m}$-primary, under mild hypothesis, the integral of the density function $f_{M,I,J}$ is the Hilbert-Kunz multiplicity $e_{HK}(J,M)$.

\begin{proposition}\label{pr: integral of the density function}
  In the local ring $(R, \mathfrak{m})$,  let $I,J_{\bullet}$ satisfy \textbf{Condition C}. Let $M$ be a finitely generated $R$-module such that $M,I, J_{\bullet}$ satisfy hypothesis (A) or (B) of \Cref{th: uniform convergence to the convex functional}. We denote the left and right derivative of $h_{M,I,J_{\bullet}}$ by $h'_{M,+}$ and $h'_{M,-}$ respectively.

    \begin{enumerate}

        \item For $0<a<b< \infty$ The functions $f_{M,I,J_\bullet}, (h_{M})'_{+}, (h_{M})_{-}'$ are integrable on $[a,b]$. Moreover the integrals of all three functions on $[a,b]$ coincide.
        
        \item Assume that $J_{\bullet}$ is $\bp{J}{n}$ for some $\mathfrak{m}$-primary ideal $J$. Then
        $$\underset{s_0 \to 0+}{\lim}\int \limits_{s_0}^{\infty} f_{M,I,J}(t)dt= e_{HK}(J,M)-\sum_{P\in \Supp(M), I \subset P,\dim R/P=\dim M}l_{R_P}(M_P)e_{HK}(J,R/P).$$
        Moreover, if $\dim(M/IM) < \dim (M)$, the integral is $e_{HK}(J,M)$.
 
        \item Let $J_\bullet$ be the same as in (2). Then $\underset{s_0 \to 0+}{\lim}\int \limits_{s_0}^{\infty} f_{M,I,J}(t)dt= \underset{s_0 \to 0+}{\lim}\int \limits_{s_0}^{\infty} (h_{M})'_{+}(t)dt= \underset{s_0 \to 0+}{\lim} \int \limits_{s_0}^{\infty} (h_{M})_{-}'(t)dt.$ In particular $(h_{M})'_{+},(h_{M})_{-}'$ are integrable on $\mathbb{R}_{>0}$.
    \end{enumerate}
\end{proposition}

\begin{proof}
(1) The integrability of $f_{M,I,J_{\bullet}}$ follows from \Cref{th: h is the integral of f} and the paragraph before it. Being continuous, $h_M$ is measurable on $\mathbb{R}_{>0}$. On $[a,b]$, $(h_{M})'_{+}, (h_{M})_{-}'$, being limits of the sequences of measurable functions $n(h_{M}(s+\frac{1}{n})- h_M(s)), -n(h_{M}(s-\frac{1}{n})- h_M(s))$ respectively, are measurable. Since $f_{M,I,J_{\bullet}}, (h_{M})'_{+}, (h_{M})_{-}'$ coincide outside a countable subset of $[a,b]$, their integrals on $[a,b]$ coincide; see \Cref{th: differentiability of h function implies existence of density function}. Since the integral of $f_{M,I,J_{\bullet}}$ on $[a,b]$ is finite by \Cref{th: h is the integral of f}, $(h_{M})'_{+}, (h_{M})_{-}'$ are also integrable on $[a,b]$.

 Now the first assertion of (2) follows directly from \Cref{th: h is the integral of f} and \Cref{th: continuity of h general case}. The assertion in (2) follows directly from the first one.
Since outside a countable subset of $\mathbb{R}$, $f_{M,I,J_{\bullet}}$, $(h_M)'_{+}$ and $(h_M)'_{-}$ coincide, the equations in (3) are immediate. We now argue the integrability of $(h_M)'_{+}$. For $s_0>0$, let $\chi_{s_{0}}$ be the characteristic function of $[s_0, \infty)$. Then,
$$\int_{\mathbb{R}_{>0}}(h_M)'_{+}(t)dt= \underset{s_0 \to 0+}{\lim}\int_{\mathbb{R}_{>0}}(h_M)'_{+}(t)\chi_{s_{0}}(t)dt= \underset{s_0 \to 0+}{\lim} \int \limits_{s_0}^{\infty} (h_{M})'_{+}(t)dt,$$
where the first equality follows from Lebesgue Monotone Convergence Theorem. Thus $(h_{M})'_{+}$ is integrable on $\mathbb{R}_{>0}$. The integrability of $(h_{M})'_{-}$ follows by a similar argument.
\end{proof}

Now we analyze the behaviour of the density function $f_{R,I,J}$ near the origin. Our argument below uses the monotonicity property proven in \Cref{pr: monotonicity result for the desnity function}. This forces us to make a simplifying assumption about $I$.

\begin{theorem}\label{th: density function near zero} Let $(R,\mathfrak{m})$ be a noetherian local domain. Let $x_1, x_2, \ldots, x_r$  be part of system of parameters of $R$, where $r \geq 1$. Let $I \supseteq (x_1, \ldots, x_r)$ be an ideal contained in the integral closure of $(x_1, \ldots, x_r)$. Let $J$ be an ideal such that $I+J$ is $\mathfrak{m}$-primary. Denote the left and right hand derivatives of $h_{R,I,J}$ at $s$ by $h'_{+}(s)$ and $h'_-(s)$ respectively. Then,

\begin{enumerate}
    \item Both $\lim_{s \to 0+} \frac{h'_{+}(s)}{s^{r-1}}$ and $\lim_{s \to 0+} \frac{h'_{-}(s)}{s^{r-1}}$ exist and coincide with
         $$r \lim_{s \to 0+} \frac{h_{R,I,J}(s)}{s^{r}}= r \frac{1}{r!}\sum_{P\, \textup{is minimal over} I}e_{HK}(J\frac{R}{P}, \frac{R}{P})e(IR_P,R_P).$$
         
     \item Let $f_{R,I,J}(s)$ be the density function associated to the pair $(I,J)$ at $s$, when it exists. Then
     $$\lim_{s \to 0+} \frac{f_{R,I,J}(s)}{s^{r-1}}= r \lim_{s \to 0+} \frac{h_{R,I,J}(s)}{s^{r}}= r \frac{1}{r!}\sum_{P\, \textup{is minimal over} I}e_{HK}(J\frac{R}{P}, \frac{R}{P})e(IR_P,R_P).$$
\end{enumerate}
\end{theorem}

\begin{proof}
The proof relies on the following comparisons, which are again used in \Cref{th: decreasingness of normalized h function}.

\begin{lemma}\label{le: bounds on higher slopes of h from convexity}
    Let $R, I$ be as in \Cref{th: density function near zero}. Let $J_{\bullet}$ be a family of ideal such that $I, J_{\bullet}$ satisfy \textbf{Condition C}. In this lemma, denote the left and right derivatives of $h_{R,I,J_{\bullet}}$ at $s$ by $h'_{+}(s)$ and $h'_-(s)$ respectively. Then, given real numbers $s>s_0$,
 \begin{equation}\label{eq: bound on higher slopes of h}
        \frac{h'_{+}(s_0)}{s_{0}^{r-1}} \geq r\frac{h_{R,I,J_{\bullet}}(s)- h_{R,I,J_{\bullet}}(s_0)}{s^r-s_0^r} \geq \frac{h'_{+}(s)}{s^{r-1}}.
    \end{equation}
 Moreover, we have similar inequalities where $h'_{+}$ is replaced by $h'_{-}$.   
\end{lemma}

\begin{proof}[Proof of lemma]
We just prove the part for $h'_{+}$ since the part for $h'_{-}$ follows from a similar argument. By \Cref{pr: invariance under closure operations}, we may assume $I= (x_1, \ldots, x_r)$. By \Cref{pr: monotonicity result for the desnity function}, $h'_{+}(s)/s^{r-1}$ is a decreasing function on the positive real line. So for positive real numbers $s_0 \leq t\leq s$,
   
    \begin{equation}\label{eq: inequality coming from decreasingness of the normlalized density function}
      h'_{+}(s_0)\frac{t^{r-1}}{s_{0}^{r-1}} \geq h'_{+}(t) \geq h'_{+}(s)\frac{t^{r-1}}{s^{r-1}}. 
    \end{equation}
    
    Thanks to \Cref{pr: integral of the density function}, \Cref{th: h is the integral of f}, the integral of $h'_{+}$ on $[s_0, s]$ is $h(s)-h(s_0)$. So taking integral of the functions in the variable $t$ on $[s_0,s]$, the above inequalities imply
    \[ \frac{h'_{+}(s_0)}{rs_{0}^{r-1}}(s^r- s_0^r) \geq h(s)-h(s_0) \geq \frac{h'_{+}(s)}{rs^{r-1}}(s^r- s_0^r),\]
    
    whence the desired inequalities in \Cref{eq: bound on higher slopes of h} follow.
\end{proof}

\noindent\textbf{Proof of theorem continued:} For assertion (1), we claim that we can assume $I= (\underline{x}):= (x_1, \ldots, x_r)$ without loss of generality. Indeed, by \Cref{pr: invariance under closure operations}, $h_{R,I,J}(s)= h_{R,(\underline{x}),J}(s)$. The right side of the equation in (1) remains unchanged after the replacement of $I$ by $(\underline{x})$ since the minimal primes over $I$ and $(\underline{x})$ coincide; see \cite[lemma 8.1.10]{HunekeSwanson}. So we assume $I= (\underline{x})$. First consider the case of $\lim_{s \to 0+} \frac{h'_{+}(s)}{s^{r-1}}$. We use \Cref{le: bounds on higher slopes of h from convexity} when $J_{\bullet}$ is $(J^{[p^{n}]})_n$. Taking $J_{\bullet}= (J^{[p^{n}]})_n$ in \Cref{eq: bound on higher slopes of h}, and then taking limit, we get the chain of inequalities,

$$\liminf_{s_0 \to 0+}\frac{h'_{+}(s_0)}{s_{0}^{r-1}} \geq \lim_{s \to 0+}\lim_{s_0 \to 0+}r\frac{h_{R,I,J}(s)- h_{R,I,J}(s_0)}{s^r-s_0^r}= \lim_{s \to 0+} r\frac{h_{R,I,J}(s)}{s^r} \geq \limsup_{s \to 0+}\frac{h'_{+}(s)}{s^{r-1}}.$$
    The claimed equality in the last chain follows as $h(t)$ approaches zero as $t$ approaches zero from right since we assume $r \geq 1$; see \Cref{th: asymptotic behaviour of h function near zero}. So we conclude that $h'_{+}(t)/t^{r-1}$ and $rh_{R,I,J}(s)/s^{r}$ have the same limit as $t$ approaches zero from the right\footnote{This assertion is true when $J^{[q]}$ is replaced by $J_{\bullet}$ where $I, J_{\bullet}$ satisfy \textbf{Condition C}.}. Rest follows from \Cref{th: asymptotic behaviour of h function near zero} once we note that all the minimal primes over $(\underline{x})$ have the same dimension.

    The case of $\lim_{s \to 0+} \frac{h'_{-}(s)}{s^{r-1}}$ follows by a similar argument.\\

    \noindent (2) Whenever $f_{R,I,J}(s)$ exists at some positive $s$, we have
    $$h^{'}_{+}(s) \leq f_{R,I,J}(s) \leq h'_{-}(s); \,\,\,\, \textup{see \Cref{th: differentiability of h function implies existence of density function}}.$$
    Rest follows from these comparisons and part (1).
\end{proof}

\begin{remark}
\begin{enumerate}
    \item Note that \Cref{th: density function near zero} includes the case when $I$ is $\mathfrak m$-primary. Indeed we can always assume that the residue field is infinite without loss of generality. So we can assume that $I$ is integral over an ideal generated by system of parameters; see \cite[Ch. 8]{HunekeSwanson}.
    
    \item When $(R,\mathfrak m)$ is not necessarily a domain, one can obtain analogues of \Cref{th: asymptotic behaviour of h function near zero} and \Cref{th: density function near zero} using the associativity formula \Cref{pr: associativity formula for h function}.
\end{enumerate}    
\end{remark}

\section{Applications}\label{se: applications}
\subsection{Comparison between Hilbert-Kunz and Hilbert-Samuel multiplicity} 
\label{sse: Watanabe's question}Let $(R, \mathfrak m)$ be an $F$-finite local ring of dimension $d$. It is well known that for an $\mathfrak m$-primary ideal $I$,
$$e(I,R) \geq e_{HK}(I,R) \geq \frac{e(I,R)}{d!}.$$
When $d$ is at least 2, Watanabe and Yoshida asked whether the rightmost inequality is always strict; see \cite[Question 2.9]{WY00}\footnote{The original question is restricted to the case $I=\mathfrak m$.}. Watanabe-Yoshida's question was affirmatively answered by Hanes by approximating an appropriate length function; see \cite[Thm 2.2, 2.4]{Hanes}. We show that Watanabe-Yoshida's question is equivalent to a question of containment of ideals, namely \Cref{pr: algebraic version of Watanabe's question}, (3), which is a priori weaker. In \Cref{th: comparison of powers and powers of tight closure}, we provide a proof of the statement appearing in \Cref{pr: algebraic version of Watanabe's question}, (3) and thus another explanation to an affirmative answer to Watanabe-Yoshida's question, which is different from Hanes'. The translation of the question regarding multiplicities to containment of ideals is facilitated by the appropriate $h$-function.

\begin{proposition}\label{pr: algebraic version of Watanabe's question}
Let $(R, \mathfrak m)$ be an $F$-finite noetherian ring of dimension $d$. For an $\mathfrak m$-primary ideal $I$, the following statements are equivalent.

\begin{enumerate}
    \item $e_{HK}(I,R) > \frac{e(I,R)}{d!}$.
    
    \item The minimal stable point $\alpha_{R,I,I}$ of $h_{R,I,I}$ is strictly larger than 1; see \Cref{th: maximal stable point of h function for J}.

    \item There exists a $q= p^e$ such that $I^{q+1}$ is not contained in $(\bpq{I}{q})^*$.
\end{enumerate}
\end{proposition}

\begin{remark}\label{re: inequality implies the non-containment}
The fact that $e_{HK}(I,R)$ is strictly greater than $e(I,R)/d!$ implies that $I^{q+1}$ can be contained in $(\bpq{I}{q})^*$ only for finitely many $q$'s. Indeed, otherwise $l(R/I^{q+1})/q^d$ is at least $l(R/(\bpq{I}{q})^*)/q^d$ for infinitely many $q$'s. Taking limit as $q$ approaches infinity, this implies $e(I,R)/d! \geq e_{HK}(I,R)$; see \Cref{le: tight closre of Frobenius power defines Hilbert-Kunz}. One implication of the previous proposition is that assertion (3), which is a priori weaker than assertion (1), implies assertion (1).
\end{remark}

\begin{proof}
The value of $h_{R,I,I}(s)$ at $1$ and $\alpha_{R,I,I}$ are $e(I,R)/d!$ and $e_{HK}(I,R)$ respectively; see \Cref{le: h function is stable near boundaries}, \Cref{th: maximal stable point of h function for J}. If (1) holds, $h_{R,I,I}(s)$ cannot be a constant on $[1, \alpha_{R,I,I}]$, so (2) follows. Now (2) implies that $h_{R,I,I}$ is a non-constant increasing function on $[1, \alpha_{R,I,I}]$; see \Cref{th: maximal stable point of h function for J}. So (1) follows.

Now we argue that (2) and (3) are equivalent. Let $r^I_I(n)$ be as in \Cref{le: F threshold upto tight closure}. Then $(r^I_I(n)/p^n)_n$ is a nondecreasing sequence converging to $\alpha_{R,I,I}$; see \Cref{le: F threshold upto tight closure}, \Cref{th: maximal stable point of h function for J}. If (2) holds, then $r^I_I(e)$ is strictly greater than $p^e$ for some $e$, so (3) follows. Conversely if (3) holds, $r^I_I(e)$ must be strictly greater than $p^e$. So $\alpha_{R,I,I}$, being the limit of the nondecreasing sequence $(r^I_I(n)/p^n)_n$ must be strictly greater than 1.
\end{proof}

The line of argument in the above proposition shows that: 

\begin{corollary}
    Suppose $J \subset I$ are two $\mathfrak m$-primary ideals in a local ring $(R,\mathfrak m)$. Suppose there exists some $q=p^e$ such that $I^{q+1} \nsubseteq (\bpq{J}{q})^*$. Then
    \[e_{HK}(J,R) > \frac{e(I,R)}{\dim(R)!}.\]
\end{corollary}

\begin{remark}\label{re: continuity of density function and Watanabe's inequality}
    We do not know what motivated Watanabe-Yoshida to formulate the question mentioned above. But from the point of view of $h$-functions, this inequality seems probable. Indeed assume additionally that the density function $f_{R,I,I}$ is continuous at 1. Since the value of the density function at 1 is $\frac{1}{(\dim R-1)!}e(I,R)>0$, $f_{R,I,I}$ remains positive in a neighborhood of $1$. This implies $\alpha_{R,I,I}>1$ and hence the inequality sought for by Watanabe-Yoshida follows; see \Cref{pr: algebraic version of Watanabe's question}. Although we do not know whether $f_{R,I,I}$ is continuous when $\textup{ht}(I)$ is at least 2, we expect that to be the case; see \Cref{qe: continuity of density function}.
\end{remark}

\begin{theorem}\label{th: comparison of powers and powers of tight closure}
Let $(R,\mathfrak{m})$ be an $F$-finite noetherian local ring of dimension $d \geq 2$ and $I$ be an $\mathfrak{m}$-primary $R$-ideal. Then there exists $q_0$ such that for all $q \geq q_0$, $I^{q+1} \nsubseteq (I^{[q]})^*$.   
\end{theorem}

\begin{proof}
Let $r^I_I(n)$ be as in \Cref{le: F threshold upto tight closure}. Since $(r^I_I(n)/p^n)_n$ is a nondecreasing sequence, once $I^{q_0+1} \nsubseteq (I^{[q_0]})^*$ for some $q_0= p^{n_0}$, $I^{q+1} \nsubseteq (I^{[q]})^*$ for all $q \geq q_0$. So we just prove that $I^{q+1} \nsubseteq (I^{[q]})^*$ for some $q$. We prove this by contradiction.

 Assume on the contrary that we have $I^{q+1} \subset (I^{[q]})^*$ for all $q$. Then for any $P \in \textup{Min}(R)$, $(I+P/P)^{q+1} \subset ((I+P/P)^{[q]})^*$ where the tight closure is taken in $R/P$. Thus we may replace $R$ with $R/P$ for some $P$ satisfying $\dim R=\dim R/P$, so we may assume $R$ is a domain. Since $R$ is $F$-finite, we can choose a test element $0 \neq c \in R$. Then $cI^{q+1} \subset I^{[q]}$ for all $q \gg 0$. We may replace $c$ by a multiple to assume $cI^{q+1} \subset I^{[q]}$ for all $q\geq 1$. We fix such an element $c$.

Now by Artin-Rees lemma, there is a number $M_0 \in \mathbb{N}$ that only depends on $c,I$ such that for all $n \geq 0$,
$$cR \cap I^{n+M_0}=I^n(cR\cap I^{M_0}).$$
Let $M=M_0+1$. We use $\mu$ to denote the number of minimal generators of an $R$-module or an $R$-ideal, and denote $r=\mu(I)$. Denote
$$N_1=\frac{cI^{q+1}+I^{[q]}I^M}{cII^{[q]}+I^{[q]}I^M},N_2=\frac{I^{[q]}}{cII^{[q]}+I^{[q]}I^M}.$$
By the assumption $N_1$ is a submodule of $N_2$. We see $\mu(N_2) \leq \mu(I)=r$ and $I^MN_2=0$, thus
$$l(N_1) \leq l(N_2) \leq rl(R/I^M).$$
Consider the short exact sequence
$$0 \to cII^{[q]}+I^{[q]}I^M \to cI^{q+1}+I^{[q]}I^M \to N_1 \to 0.$$
Tensoring with $R/I$ we get an exact sequence
$$(cII^{[q]}+I^{[q]}I^M)\otimes_RR/I \to (cI^{q+1}+I^{[q]}I^M)\otimes_R R/I \to N_1/IN_1 \to 0.$$
Thus we have
\begin{align*}
l((cI^{q+1}+I^{[q]}I^M)\otimes_R R/I) \leq l(N_1/IN_1)+l((cII^{[q]}+I^{[q]}I^M)\otimes_RR/I)\\
\leq l(N_1)+\mu(cII^{[q]}+I^{[q]}I^M)l(R/I)\leq rl(I^M)+(r^2+r\mu(I^M))l(R/I)=C
\end{align*}
Here $C$ is a constant that only depends on $I,M$ and is independent of $q$.
On the other hand, we prove that the map induced by natural inclusion $cI^{q+1}\to cI^{q+1}+I^{[q]}I^M$,
$$\varphi:cI^{q+1}\otimes_R R/I \to (cI^{q+1}+I^{[q]}I^M)\otimes_R R/I$$
is an injection. We rewrite this map as
$$\varphi:\frac{cI^{q+1}}{cI^{q+2}}\to \frac{cI^{q+1}+I^{[q]}I^M}{cI^{q+2}+I^{[q]}I^{M+1}}$$
Then we see
$$ker\varphi=\frac{cI^{q+1}\cap(cI^{q+2}+I^{[q]}I^{M+1})}{cI^{q+2}}=\frac{cI^{q+2}+cI^{q+1}\cap I^{[q]}I^{M+1}}{cI^{q+2}}$$
By our choice of $M$ we have inclusion
$$cI^{q+1}\cap I^{[q]}I^{M+1} \subset cR\cap I^{q+M+1}=cR\cap I^{q+M_0+2}\subset cI^{q+2}$$
Thus $ker\varphi=0$, so $\varphi$ is an injection. Since $c$ is a nonzero divisor, this means
$$l(I^{q+1}/I^{q+2})=l(cI^{q+1}/cI^{q+2}) \leq l((cI^{q+1}+I^{[q]}I^M)\otimes_R R/I)\leq C$$
where $C$ is independent of $q$. But it is well-known that
$$l(I^{q+1}/I^{q+2})=\frac{e(I,R)}{(d-1)!}q^{d-1}+O(q^{d-2})$$
and $d \geq 2$ implies $d-1 \geq 1$, so $\lim_{q \to \infty}l(I^{q+1}/I^{q+2})=\infty$, which is a contradiction. Thus we cannot have $I^{q+1} \subset (I^{[q]})^*$ for $q \gg 0$, so $I^{q+1} \nsubseteq (I^{[q]})^*$ for infinitely many $q$.
\end{proof}

\begin{corollary}\label{co: lower bound on F threshold}
    Let $R$ be an $F$-finite ring and $I$ be an ideal.
    \begin{enumerate}
        \item If one of the minimal primes over $I$ has height at least one, then $1 \leq r^I(I) \leq c^{I}(I)$.

        \item If one of the minimal primes over $I$ has height at least two \footnote{This assumption holds when $\text{ht}{I} \geq 2$ or $I$ is $\mathfrak{m}$-primary and $\dim(R) \geq 2$.}, then $1 < r^I(I) \leq c^I(I)$.
    \end{enumerate}   
\end{corollary}

\begin{proof}
In \Cref{le: F threshold upto tight closure is larger} we note that $r^I(I) \leq c^I(I)$. So we just prove the assertions for $r^I(I)$. For a minimal prime $Q$ over $I$ observe that $r^{IR_Q}(IR_Q) \leq r^I(I)$. So without loss of generality we replace $R$ and $I$ by $R_Q$ and $IR_Q$ respectively and assume $I$ is $\mathfrak{m}$-primary. For (1), the height assumption implies that we can assume $R$ is dimension at least one. For $s \in [0,1]$, $h_{R,I,I}(s)= e(I,R)s^{\dim(R)}/\dim(R)!$; see \Cref{le: h function is stable near boundaries}. Thus the minimal stable point of $h_{R,I,I}$: $\alpha_{R,I,I}= r^I(I)$, is at least one; see \Cref{le: F threshold upto tight closure}, \Cref{th: maximal stable point of h function for J}..

For (2), we can assume $(R, \mathfrak{m})$ has dimension at least two and $I$ is $\mathfrak{m}$-primary. Since $\dim(R)$ is at least two, $\alpha_{R,I,I}$ is strictly bigger than one by \Cref{pr: algebraic version of Watanabe's question}, \Cref{th: comparison of powers and powers of tight closure}. Since $r^I(I)= \alpha_{R,I,I}$, we are done.
\end{proof}

\begin{remark}
One can obtain an explicit $q_0$ depending on $R$ and $I$, following the arguments of \Cref{th: comparison of powers and powers of tight closure}. Such a bound should improve \Cref{co: lower bound on F threshold} and the inequality in \Cref{pr: algebraic version of Watanabe's question}, (1).
\end{remark}
We obtain another comparison of Hilbert-Samuel and Hilbert-Kunz multiplicities using $h$-functions.

\begin{proposition}\label{pr: comparison coming from behaviour of h function at infinity}
Let $(R, \mathfrak{m})$ be a local domain, $J$ be an ideal of $R$;  $f_1, f_2, \ldots, f_r$ be 
 elements in $R$, whose images in $R/J$ form a system of parameters of $R/J$. Set $I= (f_1, \ldots, f_r)R$. Then the following inequality  
 \begin{equation*}
    \sum_{P \in \Assh(R/I)}e_{HK}(J,R/P)e(IR_P,R_P) \geq \underset{Q \in \textup{Assh}(R/J)}{\sum}e(I,R/Q)e_{HK}(JR_Q, R_Q) \ .
\end{equation*}     
holds if one of the following is satisfied:
 \begin{enumerate}
     \item The height of $I$ is $r$.

    \item There exist $g_1, \ldots, g_h$ in $R$, where $h$ is the height of $J$ such that $\sqrt{J}= \sqrt{(g_1, \ldots, g_h)R}$.\footnote{i.e. $R/J$ is a set theoretic complete intersection, when $R$ is Cohen-Macaulay.} 
 \end{enumerate}
 \end{proposition}

\begin{proof}
We argue that if (2) holds then (1) also holds. Indeed, such $g_1, \ldots, g_h$ must be part of a system of parameters of $R$ as these generate an ideal of height $h$. Note that the images of $f_1, \ldots, f_r$ must be a system of parameters of $R/(g_1, \ldots, g_h)R$, since a set of elements form a part of a system of parameters if and only if their images modulo the nilradical is a part of a system of parameters. Now as  $g_1, \ldots, g_h$ is a part of a system of parameters of $R$ and images of $f_1, \ldots, f_r$ is a system of parameters of $R/(g_1, \ldots, g_h)R$, $(g_1, \ldots, g_h, f_1, \ldots, f_r)$ is a full system of parameter of $R$. So $f_1, \ldots, f_r$ must be a part of a system of parameters of $R$ proving $\text{ht}(I)=r$.

    Now we establish the desired inequality assuming (1) holds. By \Cref{th: decreasingness of normalized h function}, which we prove further below, $h_{I,J}(s)/s^r$ is a decreasing function on $\mathbb{R}_{>0}$. Thus 
    \[r!\lim_{s \to 0+}\frac{h_{I,J}(s)}{s^r} \geq r!\lim_{s \to \infty}\frac{h_{I,J}(s)}{s^r}.\]
    Since the left and right side of the above inequality evaluates to the left and right side of the desired inequality respectively, we are done; see \Cref{th: asymptotic behaviour of h function near zero}, \Cref{th: behaviour of h at infinity}.
\end{proof}

\begin{remark}\label{re: context of the inequalities}
    In the context of \Cref{pr: comparison coming from behaviour of h function at infinity}, it was known that $\underset{Q \in \text{Assh}(R/J)}{\sum}e(I,R/Q)e(JR_Q, R_Q)$ is bounded above by $e(J, R/I)$; see \cite[pn 119-122]{NormalFlatness}. These bounds appear crucially in the equimultiplicity theory; see \cite{EquimultiplicityLipman},\cite{SmirnovEquimultiplicity}.  
\end{remark}

The next two corollaries can also be deduced from Lech's associativity formula; see \cite{LechAssociative}. We provide a proof using $h$-functions.

\begin{corollary}\label{co: consequence to equimultiplicity}
Let $(R, \mathfrak{m})$ be a local domain, $(x_1, \ldots, x_d)$ be a full system of parameters of $R$. For an integer $1 \leq r \leq d$, set $I=(x_1, \ldots, x_r)$ and $J= (x_{r+1}, \ldots, x_d)$. Then 

\begin{equation}\label{eq: equality in equimultiplicity}
 \sum_{P \in \Assh(R/I)}e_{HK}(J,R/P)e(IR_P,R_P) = \underset{Q \in \text{Assh}(R/J)}{\sum}e(I,R/Q)e_{HK}(JR_Q, R_Q) \ .
\end{equation} 
\end{corollary}

\begin{proof}
The left side of \Cref{eq: equality in equimultiplicity} is at least the right side by \Cref{pr: comparison coming from behaviour of h function at infinity}. Switching the roles of $I$ and $J$ in \Cref{pr: comparison coming from behaviour of h function at infinity}, we conclude that the right side is at least the left side. 
\end{proof}

\begin{corollary}\label{co: consequence to equimultiplicity when the residue field is infinite}
    Suppose the residue field of the local domain $(R, \mathfrak{m})$ is infinite. Let $I$ be an ideal whose analytic spread is the same as the height of $I$. Let $x_1, \ldots, x_t$ be a set of elements whose image in $R/I$ forms a full system of parameters. Set $J=(x_1, \ldots, x_t)R$. Then

\begin{equation}\label{eq: equality in equimultiplicity with infinite residue field}
 \sum_{P \in \Assh(R/I)}e_{HK}(J,R/P)e(IR_P,R_P) = \underset{Q \in \text{Assh}(R/J)}{\sum}e(I,R/Q)e_{HK}(JR_Q, R_Q) \ .
\end{equation}    
\end{corollary}

\begin{proof}
    Since the analytic spread coincides with the height and the residue field is infinite, we can find part of a system of parameters $y_1, \ldots, y_h$ of $R$ such that if $I'= (y_1, \ldots, y_h)$, then $I' \subset I \subset \overline{I'}$. The elements $x_1, \ldots, x_t$ form a full system of parameters of $R/I'$ since $I, I'$ have the same radical; see the first paragraph of the proof of \Cref{pr: comparison coming from behaviour of h function at infinity}. Recall $\text{Assh}(R/I)= \text{Assh}(R/I')$; see \cite[lemma 8.1.10]{HunekeSwanson}. Moreover for any $P \in \text{Assh}(R/I)$ or $Q \in \text{Assh}(R/J)$, $e(IR_P, R_P)= e(I'R_P, R_P)$ and $e(I, R/Q)= e(I', R/Q)$, it is enough to establish \Cref{eq: equality in equimultiplicity with infinite residue field} after replacing $I$ by $I'$. Since $y_1, \ldots, y_h, x_1, \ldots, x_t$ form a full system of parameters of $R$, we can use \Cref{co: consequence to equimultiplicity} to conclude the proof. 
\end{proof}

 As a consequence of \Cref{pr: comparison coming from behaviour of h function at infinity}, we show that $h_{I,J}$ takes an especially simple form when $I,J$ `splits a full system of parameters'.

\begin{proposition}\label{pr: h functions when the ideals split a full system of parameters}
Let $(R, \mathfrak{m})$ be a local domain, $(x_1, \ldots, x_d)$ be a full system of parameters of $R$. For an integer $0 \leq r \leq d$, set $I=(x_1, \ldots, x_r)$ and $J= (x_{r+1}, \ldots, x_d)$, with the understanding that $I=0$, when $r=0$. Then for any $s \in \mathbb{R}_{>0}$,
\[h_{I,J}(s)= \frac{e(I+J, R)}{r!}s^r. \]
\end{proposition}

\begin{proof}
   When $I=0$, on the positive real axis, $h_{I,J}(s)$ is the constant function $e_{HK}(J,R)=e(J,R)$. So we assume $r>0$. In this case, by \Cref{th: decreasingness of normalized h function}, $h_{I,J}/s^r$ is a decreasing function on $\mathbb{R}_{>0}$. Our hypothesis on $I, J$ implies 
    \[\lim_{s \to 0+}\frac{h_{I,J}(s)}{s^r} = \lim_{s \to \infty}\frac{h_{I,J}(s)}{s^r};\]
    see \Cref{th: asymptotic behaviour of h function near zero}, \Cref{th: behaviour of h at infinity}, \Cref{co: consequence to equimultiplicity}.  So on  $\mathbb{R}_{>0}$, $h_{I,J}(s)/s^r$ must be the constant function with value 
    $$\lim_{s \to \infty}\frac{h_{I,J}(s)}{s^r}= \frac{1}{r!}\underset{Q \in \text{Assh}(R/J)}{\sum}e(I,R/Q)e_{HK}(JR_Q, R_Q); \,\, \text{see} \, \text{\Cref{th: behaviour of h at infinity}}.$$
    Since the analytic spread of $I$ is the same as the height of $I$,
    $$\underset{Q \in \text{Assh}(R/J)}{\sum}e(I,R/Q)e_{HK}(JR_Q, R_Q)= e(I+J, R),$$
    by \cite{LechAssociative}. So $h_{I,J}(s)= e(I+J, R)s^r/r!$ for any positive $s$.
\end{proof}

\subsection{$F$-threshold and multiplicity}\label{sse: F-threshold and multiplicity}Comparisons among Hilbert-Samuel multiplicity, Hilbert-Kunz multiplicity, $F$-threshold are abound in the literature. We show that general properties of the $h$ function combined with a very coarse approximation of it recover some of these.\\

Motivated by the comparison between Hilbert-Samuel multiplicity and log canonical threshold in \cite{EFM} and the analogy between $F$-threshold and log canonical threshold (see \cite[Thm 3.3, 3.4]{mustata2004f}) the following was conjectured:

\begin{conjecture}\label{co: F-threshold and multiplicity}(see \cite[Conj 5.1]{HMTW}) Let $(R,\mathfrak m)$ be a noetherian local ring of dimension $d$ containing a positive characteristic field. Let $J$ be an ideal generated by a full system of parameters and $I$ be an $\mathfrak m$-primary ideal. Then 
$$e(I,R) \geq \frac{d^d}{c^J(I)^d}e(J,R).$$
Here $e(-,R)$ denotes the Hilbert-Samuel multiplicity of the corresponding ideal.
\end{conjecture}
\begin{remark}
    We can assume $I$ is generated by a system of parameters without loss of generality, in \Cref{co: F-threshold and multiplicity}. Indeed in \Cref{co: F-threshold and multiplicity} one can first assume that the residue field is infinite by making standard constructions. Recall $e(I,R)= e(\overline{I},R)$ and $c^J(I)= c^J(\overline{I})$, where $\overline{I}$ is the integral closure of $I$ (see \cite[Prop 2.2, (2)]{HMTW}). When the residue field is infinite, we can choose a system of parameters $f_1, \ldots, f_d$ such that $I= \overline{(f_1, \ldots, f_d)}$.
\end{remark}

The above conjecture is settled when $(R,\mathfrak m)$ is graded; see \cite{HunekeTakagi}, \cite{HMTW}. Drawing motivations from \cite[Prop 4.5]{TW} which confirms a special case of the above conjecture the following conjecture was made:

\begin{conjecture}\label{co: F-threshold and Hilbert-Kunz}(see \cite[Conj 1.1]{BetancortSmirnov})
Let $f_1, f_2, \ldots, f_r$ be part of a system of parameters of a noetherian local ring $(R, \mathfrak m)$ of prime characteristic. Let $J$ be an $\mathfrak m$-primary ideal. Set $I= (f_1, \ldots, f_r)R$. Then
$$e_{HK}(J,R) \leq (\frac{c^J(I)}{r})^re_{HK}(J\frac{R}{I}, \frac{R}{I}).$$
\end{conjecture}
We next point out that \Cref{co: F-threshold and Hilbert-Kunz}, as stated, is false even when $R$ is regular.

\begin{proposition}\label{pr: Smirnov Betancourt conjecture is false}
Take $(R, \mathfrak m)$ to be the localization of a polynomial ring in $d$ variables over a prime characteristic field, where $d$ is at least two. Take $J=\mathfrak m^t$ and $I= \mathfrak m$. Then for large enough $t$,
$$e_{HK}(J,R) > (\frac{c^J(I)}{d})^de_{HK}(J\frac{R}{I}, \frac{R}{I}).$$
Thus for large $t$, \Cref{co: F-threshold and Hilbert-Kunz} fails.
\end{proposition}

\begin{proof}
    Since $R$ is regular, $e_{HK}(J,R)$ is the same as $l(R/J)$ which is just $\binom{d+t-1}{d}$. The $F$-threshold is $t+d-1$; see \cite[Example 2.7, (iii)]{HMTW}. When $d$ is at least two, for large enough $t$,
    $$\binom{t+d-1}{d} > (\frac{t+d-1}{d})^d.$$
    So we are done.
\end{proof}

We now relate \Cref{co: F-threshold and multiplicity} to \Cref{co: F-threshold and Hilbert-Kunz}. 

\begin{proposition}\label{pr: equivalence of conjectures}
Let $(R, \mathfrak m)$ be a noetherian local ring of prime characteristics. The following are equivalent:
\begin{enumerate}
    \item For every pair of $\mathfrak{m}$-primary ideals $I,J$ generated by system of parameters
    $$e(I,R) \geq \frac{d^d}{c^J(I)^d}e(J,R).$$
    \item For every pair of $\mathfrak{m}$-primary ideals $I,J$ generated by system of parameters
    $$e_{HK}(J,R) \leq (\frac{c^J(I)}{d})^de_{HK}(J\frac{R}{I}, \frac{R}{I}).$$
  That is the restricted case of \Cref{co: F-threshold and Hilbert-Kunz}, where both $I,J$ are generated by system of parameters, is equivalent to \Cref{co: F-threshold and multiplicity}.     
\end{enumerate}
\end{proposition}

\begin{proof}
    Assume (1). For an ideal generated by system of parameters the Hilbert-Kunz and Hilbert-Samuel multiplicities are the same; see \cite[Thm 2]{LechAssociative}. Since $I$ is generated by a system of parameters 
    $$e_{HK}(J\frac{R}{I}, \frac{R}{I})= l(\frac{R}{I}) \geq e(I,R).$$
   This implies (2).\\

    Now assume (2). Choose a system of parameters $f_1, f_2, \ldots, f_d$ so that $I=(f_1, f_2, \ldots, f_d)$. For any positive integer $n$, (2) yields,
    $$e_{HK}(J,R) \leq (\frac{c^J((f_1^n, f_2^n, \ldots, f_d^n))}{d})^dl(\frac{R}{(f_1^n, f_2^n, \ldots, f_d^n)}).$$
    Since $I^n$ is in the integral closure of $(f_1^n, f_2^n, \ldots, f_d^n)$, $c^J((f_1^n, f_2^n, \ldots, f_d^n))= c^J(I^n)$. Moreover $c^J(I^n)=c^J(I)/n$; see \cite[Prop 2.2, (3)]{HMTW}. So the last inequality gives
    $$e_{HK}(J,R) \leq (\frac{c^J(I)}{d})^d\frac{l(\frac{R}{(f_1^n, f_2^n, \ldots, f_d^n)})}{n^d},$$
    for all $n$. Taking limit as $n$ approaches infinity in the last inequality, we obtain (1).
\end{proof}

In view of the above proposition, we believe that the corrected version of \Cref{co: F-threshold and Hilbert-Kunz} should be as follows:

\begin{conjecture}\label{co: corrected F-threshold and Hilbert-Kunz}
Let $f_1, f_2, \ldots, f_r$ be part of a system of parameters of a noetherian local ring $(R, \mathfrak m)$ of prime characteristic. Let $J$ be an ideal generated by a (full) system of parameters of $R$. Set $I= (f_1, \ldots, f_r)R$. Then
$$e_{HK}(J,R) \leq (\frac{c^J(I)}{r})^re_{HK}(J\frac{R}{I}, \frac{R}{I}).$$    
\end{conjecture}

The comparison between Hilbert-Kunz multiplicity and the $F$-threshold proven in  \cite[Prop 2.1]{BetancortSmirnov} verifies \Cref{co: corrected F-threshold and Hilbert-Kunz} when $r=1$. In assertion (2), of \Cref{th: Hilbert-Kunz multiplicity and F-thresold}, we strengthen the comparison in 
\cite[Prop 2.1]{BetancortSmirnov} by using the property of $h$-function proven below.

\begin{theorem}\label{th: decreasingness of normalized h function}
Let $(R,\mathfrak m)$ be a noetherian local domain. Let $x_1, x_2, \ldots, x_r$ be elements of $R$, with $r$ being a nonnegative integer. Let $I \supseteq (x_1, \ldots, x_r)R$ be an ideal contained in the integral closure of $(x_1, \ldots, x_r)$. When $r$ is zero $I$ is understood to be the zero ideal. Let $J_\bullet$ be a family of ideals such that $I, J_\bullet$ satisfies \textbf{Condition C}. Then
$h_{R, I, J_{\bullet}}(s)/s^{r}$ is a decreasing function on $(0, \infty)$.
\end{theorem}

\begin{proof}
When $r$ is zero, $h_{I,J}$ is a constant function. So $r$ is assumed to be positive for the rest of the proof. By \Cref{pr: invariance under closure operations}, we can assume $I= (x_1, \ldots, x_r)R$. Thanks to \Cref{le: bounds on higher slopes of h from convexity}, for real numbers $s_0<s_1<s_2$, we have
$$\frac{h_{R,I, J_\bullet}(s_1)- h_{R,I, J_\bullet}(s_0)}{s_1^r- s_0^r} \geq \frac{h'_+(s_1)}{rs_1^{r-1}} \geq \frac{h_{R,I, J_\bullet}(s_2)- h_{R,I, J_\bullet}(s_1)}{s_2^r- s_1^r}.$$

Take limit as $s_0$ approaches zero from right in the above chain of inequalities. Since $r \geq 1$, $h_{I, J_{\bullet}}(s_0)$ approaches zero (see \Cref{th: asymptotic behaviour of h function near zero}). So,
$$\frac{h_{R,I, J_\bullet}(s_1)}{s_1^r} \geq \frac{h_{R,I, J_\bullet}(s_2)- h_{R,I, J_\bullet}(s_1)}{s_2^r- s_1^r}.$$
The last inequality implies
$h_{R,I, J_\bullet}(s_1)/s_1^r \geq h_{R,I, J_\bullet}(s_2)/s_2^r.$
\end{proof}

\begin{theorem}\label{th: Hilbert-Kunz multiplicity and F-thresold}
Let  $f_1, f_2, \ldots, f_r$ be part of a system of parameters of a noetherian local domain $(R, \mathfrak m)$ of prime characteristics with $r$ positive; set $I=(f_1, \ldots, f_r)$. Let $J$ be an $\mathfrak m$-primary ideal. Let $\alpha_{R,I,J}$ be the minimal stable point of $h_{R,I,J}$ as defined in \Cref{sse: tail of h function}. Then

\begin{enumerate}
    \item $$e_{HK}(J,R) \leq \frac{\alpha_{R,I,J}^r}{r!}\sum_{P\, \textup{is a minimal over prime} \,I}e_{HK}(J\frac{R}{P}, \frac{R}{P})e(IR_P,R_P).$$
    \item $$e_{HK}(J,R) \leq \frac{c^J(I)^r}{r!}e_{HK}(J\frac{R}{I},\frac{R}{I}).$$
\end{enumerate}
\end{theorem}

\begin{proof}
    First we point out (1) implies (2). We know $\alpha_{R,I,J} \leq c^J(I)$ (see \Cref{le: h function is stable near boundaries}, assertion (2)) and $e(IR_P,R_P) \leq l_{R_p}(R_P/IR_P)$. The last inequality holds as $IR_P$ is generated by a system of parameters. Using these two comparisons in (1), we get  
 $$e_{HK}(J,R) \leq \frac{c^J(I)^r}{r!}\sum_{P\, \textup{minimal prime over} I}e_{HK}(J\frac{R}{P}, \frac{R}{P})l_{R_p}(\frac{R_P}{IR_P}).$$
 The right hand side of the above inequality is $e_{HK}(J\frac{R}{I}, \frac{R}{I})$; see \cite[Thm 3.14]{HunekeExp}. So (2) follows.\\

For (1), note, since $h_{R,I,J}(s)/s^r$ is decreasing on $(0, \infty)$ by \Cref{th: decreasingness of normalized h function}, we have
$$\lim_{s \to 0+}\frac{h_{R,I,J}(s)}{s^r} \geq \frac{h_{R,I,J}(\alpha_{R,I,J})}{\alpha_{R,I,J}^r}= \frac{e_{HK}(J,R)}{\alpha_{R,I,J}^r}.$$
The above one sided limit is
$$\frac{1}{r!}\underset{P\, \textup{minimal prime over}\,I}{\sum}e_{HK}(J\frac{R}{P}, \frac{R}{P})e(IR_P,R_P),$$ 
by \Cref{th: asymptotic behaviour of h function near zero}. So assertion (1) follows.

\end{proof}

Now we point out an equivalent formulation of \Cref{co: F-threshold and multiplicity} phrased in terms of $h$ functions.

\begin{proposition}\label{pr: HMTW conjecture in terms of the h function}
Let $(R,\mathfrak{m})$ be a local domain of Krull dimension $d \geq 1$, $J$ be an ideal generated by system of parameters, $I$ be an $\mathfrak m$-primary ideal. The following are equivalent.

\begin{enumerate}
    \item $$e(I,R) \geq \frac{d^d}{c^J(I)^d}e(J,R).$$
    
    \item There exists $x_0 \in [0, c^J(I)]$ such that
$$\frac{e(I,R)}{d^d} \geq \frac{h_{R,I,J}(x_0)}{x_0^d}.$$
\end{enumerate}
\end{proposition}

\begin{proof}
Assume (1) holds. Take $x_0= c^J(I)$. By \Cref{le: h function is stable near boundaries}, we have $h_{R,I,J}(x_0)= e_{HK}(J,R)= e(J,R)$- the last equality follows because $J$ is generated by a full system of parameters. Thus $(1)$ implies the desired inequality when $x_0= c^J(I)$.

Now assume (2) holds. Using \Cref{re: changing the residue field}, we can assume, without loss of generality that the residue field of $R$ is infinite. Since the residue field is infinite, $I$ admits a minimal reduction generated by a full system of parameters. So by \Cref{th: decreasingness of normalized h function},  $h_{R,I,J}(s)/s^d$ is decreasing on $\mathbb{R}_{\geq 0}$. So
$$\frac{h_{R,I,J}(x_0)}{x_0^d} \geq \frac{h_{R,I,J}(c^J(I))}{c^J(I)^d}= \frac{e(J,R)}{c^J(I)^d}.$$
Thus (1) follows from (2).
\end{proof}

\section{Questions}\label{se: questions}

\begin{question}\label{qe: minnimal stable point vs F-threshold}
    Let $(R, \mathfrak m)$ be an $F$-finite ring; $J$ be an $\mathfrak m$-primary ideal, $I$ be any ideal. Is the minimal stable point $\alpha_{R,I,J}$ of $h_{R,I,J}$ the same as the $F$-threshold $c^J(I)$? 
\end{question}

In view of \Cref{th: maximal stable point of h function for J}, the above question is a question about asymptotic comparisons of $\bpq{J}{q}$ and $(\bpq{J}{q})^*$. Moreover, in view of \Cref{th: Hilbert-Kunz multiplicity and F-thresold}, one may hope to replace $c^J(I)$ by potentially the smaller number $\alpha_{R,I,J}$ in \Cref{co: corrected F-threshold and Hilbert-Kunz} or \Cref{co: F-threshold and multiplicity}. So this question tests the veracity of this naive hope. See \Cref{re: Trivedi's question regarding the support} for the connection of this question with a question of Trivedi-Watanabe.

\begin{question}\label{qe: continuity of density function}
Let $(R, \mathfrak m)$ be an $F$-finite ring; $I, J$ be ideals such that $I+J$ is $\mathfrak m$-primary.
\begin{enumerate}
    \item Find conditions on $(R,I,J)$ such that the limit defining the Hilbert-Kunz density function $f_{R,I,J}(s)$:
    $$\lim_{q \to \infty}\frac{1}{q^{d-1}}l(\frac{
    I^{\ceil{sq}} +\bpq{J}{q}}{I^{\ceil{sq}+1}+\bpq{J}{q}}),$$
    exists at all $s \in \mathbb R$. For $s>0$, does the above limit coincide with the right derivative of $h_{I,J}$ at $s$? 

    \item Find conditions on $(R,I,J)$ such the Hilbert-Kunz density function $f_{R,I,J}(s)$ is continuous on $(0, \infty)$.
\end{enumerate}

Recall that continuity of the density function is equivalent to the corresponding $h$-function being continuously differentiable; see \Cref{th: continuous differentiability vs continuity of the density function}. Our result suggests that a larger value of $\textup{ht}(I)$ may imply a better smoothness property of the $h$-function; see \Cref{th: asymptotic behaviour of h function near zero}, \Cref{th: h function and integration}. So we wonder if $f_{R,I,J}$ exists and is continuous on $\mathbb{R}$ when $\textup{ht}(I)$ is at least two; see \Cref{re: continuity of density function and Watanabe's inequality} for a consequence of affirmative answers. Recall that when $R$ is standard graded of dimension at least two and $I= \mathfrak m$, for any homogeneous ideal $J$, the answer to both the questions are affirmative; see \Cref{th: HK density when J is homogeneous but not of finite colength}. 

When $\text{ht}(I)=1$, the right hand derivative of $h_{R,I,J}$ at zero can differ from $f_{R,I,J}(0)$. Indeed, for $R= \mathbb{F}_p[[x,y]]/(x^2- y^3); I=J= (x,y)$, a direct computation shows $f_{R,I,J}(0)= 1$. But by \Cref{th: asymptotic behaviour of h function near zero}, the right derivative of $h_{R,I,J}$ at zero is two. If we additionally assume $R$ is normal, one can show that the right derivative of $h_{R,I,J}$ at zero coincides with $f_{R,I,J}(0)$ by using the fact that $R$ is regular in codimension one. \\

\end{question}

Inspired by Trivedi's question \cite[Question 2]{TrivediQuadric}, we ask

\begin{question}\label{qe: piecewise polynomial}
    Let $I, J$ be $\mathfrak m$-primary ideals of a noetherian local ring $R$ of dimension at least two. Is $h_{R,I,J}$ a piecewise polynomial? In other words, does there exist a countable subset $S$ of $\mathbb R$ and a partition $\mathbb R \setminus S= \underset{n \in \mathbb N}{\coprod} (a_n,b_n)$ such that on each $(a_n, b_n)$, $h_{R,I,J}$ is given by a polynomial function? 
\end{question}
We point out that, in the context of the question, $h_{R,I,J}(s)$ is $e_{HK}(J,R)$ for large $s$, $e(I,R)s^{\dim(R)}/\dim(R)!$ on some interval $(0, a]$ and zero for $s$ nonpositive.
\section*{Acknowledgements}
The first author was supported in part by NSF-FRG grant DMS-1952366. The second author thanks support of NSF DMS \# 1952399 and \# 2101075. We thank Linquan Ma for supporting our collaboration.

This material is based upon work supported by the National Science Foundation under Grant No. DMS-1928930 and by the Alfred P. Sloan Foundation under grant G-2021-16778, while the first author was in residence at the Simons Laufer Mathematical Sciences Institute (formerly MSRI) in Berkeley, California, during the Spring 2024 semester.

\printbibliography[title= {References}]
\end{document}